\documentclass[twocolumn]{autart}

\usepackage{amssymb,epsfig}
\usepackage[noadjust]{cite}
\usepackage[
]{amsmath}
\usepackage{graphicx}
\usepackage{amsfonts}
\usepackage{multirow}
\usepackage{caption}
\usepackage{enumitem}
\usepackage{algorithm,algorithmic} 

\usepackage{makecell}
\usepackage{kbordermatrix}

\DeclareMathOperator{\rank}{rank}
\DeclareMathOperator{\diag}{diag}
\DeclareMathOperator{\offdiag}{offdiag}
\DeclareMathOperator{\ind}{ind}

\def\RPAPDSF{\emph{RPA-PDSF}}
\renewcommand\Re{\mathrm{Re}}
\renewcommand\Im{\mathrm{Im}}
\def\step{\!:\!}
\newenvironment{bgapmatrix}[1][\,]{\left[#1\def\bgaptmp{#1}\begin{matrix}}{\end{matrix}\bgaptmp\right]\let\bgaptmp\undefined}
\newenvironment{row}[1][\,]{[#1\def\bgaptmp{#1}\begin{matrix}}{\end{matrix}\bgaptmp]\let\bgaptmp\undefined}
\newenvironment{bsmallmatrix}{\left[\begin{smallmatrix}}{\end{smallmatrix}\right]}

\begin{document}

\begin{frontmatter}

\title{Schur Method for Robust Pole Assignment of Descriptor Systems via Proportional Plus Derivative State Feedback}

\author[Guo]{Zhen-Chen Guo}\ead{guozhch06@gmail.com}

\address[Guo]{Shing-Tung Yau Center,  College of Sciences, National Chiao Tung University, Hsinchu 30010, Taiwan, R.O.C.}


\begin{abstract}                          
The pole assignment problem for descriptor systems is a classical
inverse algebraic eigenvalue problem, which has attracted attention for decades.
In this paper, we propose a direct method to solve
the problem with the application of the proportional plus derivative
state feedback,  intending to  obtain a robust closed-loop system.
Theorems on the feasibility of our method will be presented.
Numerical examples show that our method yields poles of high
relative accuracy.
\end{abstract}

\end{frontmatter}




\section{Introduction}\label{section1}
Consider the linear time-invariant dynamical system 
\begin{equation}\label{descriptor system}
E\dot{x}(t)=Ax(t)+Bu(t),
\end{equation}
where $E, A \in\mathbb{R}^{n\times n}$, $B\in\mathbb{R}^{n\times m}$,
$x(t)\in\mathbb{R}^n$ and $u(t)\in\mathbb{R}^m$
are the state and the input, respectively.
The system is  referred to as a descriptor system and has
found  a variety of applications,  such as
chemical processes and electrical network control \cite{Dai, KD,  Ria}. 
With a generally singular $E$, \eqref{descriptor system} is an algebraic-differential equation which attracts much recent interests.     Introduced in \cite{Lu}, studies of \eqref{descriptor system} include several meaningful mathematical problems, motivated intrinsically from the associated engineering design, such as
its controllability, regularization,  pole assignment  and so on.
Please refer to \cite{BBMN, BMN, BGM, CBS, CCH, CH1, Duan, DP1, DP3, DP2, FLE, FKN, Hou, IT, KNC, KuZ, KD,
Le1, LC, Lu, MI, OL, RZ, Sha, SL,  Van, Var2, YS, ZB, ZB1} and the references therein for more information.


When the infinity index $\ind_{\infty}(A,E)$ (or the maximal size of the Jordan blocks in the Weierstrass canonical form of the matrix pair $(A,E)$ or equivalently the matrix pencil $A - \lambda E$ corresponding to its infinite generalized eigenvalue) is no greater than 1 and $(A,E)$ is regular,
the algebraic part (or the associated redundant variables) in \eqref{descriptor system} can be  eliminated, resulting in a standard linear system (with a nonsingular $E$) of reduced order.
In contrary, systems with  $\ind_{\infty}(A,E)>1$ might lost causality for some insufficiently smooth inputs.
Consequently, one hopes to obtain a regular  closed-loop system
with an infinity index being no greater than $1$ after applying feedback.
Fortunately, authors in \cite{BMN} have shown that,
if $(E, A, B)$ is strongly controllable (or S-controllable),
a proportional plus derivative state feedback (PD-SF)
exists for such a design goal.

Regarding the  pole assignment problem,
which is of some importance for system design, the dynamical behaviour
of  \eqref{descriptor system} fundamentally depends  on
the eigen-structure of $(A, E)$, especially  the  eigenvalues \cite{BMN}.
When only the state is available,
the proportional state feedback  will be adopted \cite{KNC}; if  the derivative of the state
can be measured, we may apply the derivative state feedback \cite{ATF}.
When both are procurable, we may employ a PD-SF \cite{DP2}.
All these state feedback designs are also applicable for output feedback.

It is worthwhile to point out that a state feedback involving derivatives
has some advantages over one without.
More specifically, by modifying $E$ to $E+BG$  for some $G\in\mathbb{R}^{m\times n}$,
we could regularize the closed-loop descriptor system,
assigning $\rank(\begin{row}E& B\end{row})$ finite poles meanwhile
shifting some infinite ones.  Consequently, \eqref{descriptor system}  %
may be converted into a standard one of reduced order, under  certain conditions,
eliminating the algebraic part.

For the multi-input system (i.e., $m>1$), 
many different PD-SFs lead to regular  closed-loop systems, which has
an infinity index no higher than $1$ and the finite eigenvalues are $r$ specified complex numbers
(closed under complex conjugation),  with $\rank(E)\leq r\leq \rank(\begin{row}E& B\end{row})$.
In applications,  one may prefer a PD-SF which produces a robust closed-loop system.
Applying the regularization results in \cite{BMN},
we will focus on the \emph{robust pole assignment problem via
the proportional plus derivative state feedback  (RPA-PDSF)},
which is stated as follows:

\noindent{{\em RPA-PDSF:}}
For given $E, A\in\mathbb{R}^{n\times n}$, $B\in\mathbb{R}^{n\times m}$ with $(E, A, B)$ being S-controllable,
and an arbitrarily set $\mathfrak{L}=\{(\alpha_{1}, 0), \ldots, (\alpha_{n-r},0); (\alpha_{n-r+1}, \beta_{n-r+1}),
\ldots, (\alpha_n, \beta_n)\}$, closed under complex conjugation,
where  $\beta_j\neq0$  for  $j=n-r+1, \ldots, n$, with
$q-m \leq r\leq q$, $q\triangleq \rank(\begin{row}E&B\end{row})$,
find a pair of matrices $G, F\in\mathbb{R}^{m\times n}$,
such that $(A+BF, E+BG)$ is regular, $\ind_{\infty}(A+BF,E+BG) \leq 1$,  the spectrum $\lambda(A+BF, E+BG)=\mathfrak{L}$
and  the obtained  closed-loop system is robust, that is, the eigenvalues of
$(A+BF, E+BG)$  are as insensitive to perturbations on $(A+BF, E+BG)$ as possible.

Here we represent an eigenvalue $\lambda_j = \alpha_j/\beta_j$ by the ordered-pair $(\alpha_j, \beta_j)$, eliminating the distinction between finite and infinite eigenvalues. Note that $(\alpha_j, \beta_j)$ is a representative of an equivalence class defined by the relation $\sim$, where  $(\alpha_i, \beta_i) \sim (\alpha_j, \beta_j) \Leftrightarrow \alpha_i \beta_j = \alpha_j \beta_i$. Without loss of generality, we restrict $\alpha_1, \cdots, \alpha_{n-r}$ to be real.

Adopting different measures of robustness, different methods were
proposed to solve the  \RPAPDSF.
Two frequently used measures are the condition number of the eigenvectors matrix \cite{Var2} and the
departure from normality \cite{LC}.
When \eqref{descriptor system} is  completely controllable (C-controllable), adopting
the condition numbers of the left and right  eigenvector matrices of $(A+BF, E+BG)$ as the measure,
\cite{Var2} solved  the \RPAPDSF\ through a series of generalized Sylvester equations and the Weierstrass canonical form of $(A+BF, E+BG)$.
Arbitrary pole placement were permitted, under the harsh assumption that
the sizes of all the Jordan blocks (for both finite and infinite eigenvalues) are known {\it a priori}.
Computing the  Weierstrass canonical form  would also cause some numerical instability in general.
Furthermore, the accuracy in solving the  generalized Sylvester equations
relies on a wide separation between $\lambda(A, E)$ and $\lambda(A+BF, E+BG)$, thus placing an unreasonable
demand. (After all, why should some well-behaved poles not allowed to remain?)
Recently, a Schur-Newton method was proposed in \cite{LC}, minimizing  the departure from normality of $(A+BF, E+BG)$.
With the generalized Schur form $(A+BF, E+BG)=(XSP,XTP)$, where $X, P\in\mathbb{R}^{n\times n}$ are
nonsingular, $S, T\in\mathbb{R}^{n\times n}$ are upper quasi-triangular and
all  finite poles are real,  the method in  \cite{LC}
generates an orthogonal $P$. For complex conjugate poles,
the acquired $P$ is usually not orthogonal, implying that it virtually does  not optimize
the departure from normality of  $(A+BF, E+BG)$.
Both methods are iterative and convergence are not proven.

In  \cite{DP2}, Duan and Patton employed the proportional plus partial derivative state feedback, with the closed-loop system in the form $(A+BF, E+BGC)$ for $C\in\mathbb{R}^{l\times n}$ being the output matrix.
Adopting a sum  of the condition numbers of individual eigenvalues
as the robust measure, the method assigns $n$ distinct finite poles, requiring the existence of $G\in\mathbb{R}^{m\times l}$ with $\rank(E+BGC)=n$.
However, no sufficient and necessary condition is offered to guarantee such existence.
Besides, the method essentially computes the right coprime polynomial matrices $N(s)$ and $D(s)$
such that $(A-sE)N(s)+BD(s)=0$ \cite{Duan0}, which is theoretically elegant
yet numerically difficult to implement.

Inspired by the algorithms  \verb|schur| \cite{Chu} and  \verb|Schur-rob|  \cite{GCQX},
 we  propose a direct method to solve the \RPAPDSF, utilizing the generalized real  Schur form of $(A+BF, E+BG)$ in this paper. We shall adopt
a robustness measure which is closely related to the departure from normality.
All poles will be placed in turn, and
in each step (which assigns  an infinite pole, a real pole or a pair of complex conjugate poles),
we minimize the robust measure in an optimization subproblem. 
When assigning an infinite pole, we merely need to solve some linear equations;
while assigning a real pole, only a singular value decomposition (SVD)
is required. When  assigning a pair of complex conjugate poles,
similarly to \verb|Schur-rob|, an efficient solution of
the corresponding optimization sub-problem is proposed.
In addition, theorems will be proved to guarantee the feasibility of our method.
Abundant amount of numerical results will show the feasibility and efficiency of our method,
producing robust closed-loop systems with highly accurate finite poles.

The paper is organized as follows. In Section~\ref{section2}, we present some preliminaries. Our method is developed in Section~\ref{section3},
for the assignment of infinite poles, real finite poles and complex conjugate finite poles.
Numerical results  are reported in Section~\ref{section4}.
Some concluding remarks are then made in Section~\ref{section5}.

\section{Notations and Preliminaries}\label{section2}

Throughout this paper, for an arbitrary  matrix $M$, we denote the null space,
the range space and the submatrix comprised
by rows $k$ to $l$ and columns $s$ to $t$ by $\mathcal{N}(M)$, $\mathcal{R}(M)$ and
 $M(k\step l, s\step t)$, respectively. For any arbitrary $\lambda\in\mathbb{C}$,
 define $D_0(\lambda) \equiv \begin{bmatrix}\Re(\lambda)&\Im(\lambda)\\
-\Im(\lambda)&\Re(\lambda)\end{bmatrix}$, $D_{\delta}(\lambda) \equiv \begin{bmatrix}\Re(\lambda)&\delta\Im(\lambda)\\
-\delta^{-1}\Im(\lambda)&\Re(\lambda)\end{bmatrix}$ with some $0\neq \delta\in\mathbb{R}$.
We adopt ``$\offdiag(T)$" for the strictly upper triangular part of a matrix  $T$.

%

\begin{lem}\label{Theorem2.1}
For any regular matrix pencil $(A, E)$, $A, E\in\mathbb{R}^{n\times n}$, there exist a
nonsingular matrix $X\in\mathbb{R}^{n\times n}$ and an orthogonal matrix $P\in\mathbb{R}^{n\times n}$
such that $X^{-1}AP=S$ and $X^{-1}EP=T$ are both upper
quasi-triangular with $1\times 1$ or $2\times 2$ diagonal blocks.
Moreover, writing  the block diagonal parts of $S$ and $T$ as
$\Phi=\diag(\Phi_{11}, \ldots, \Phi_{kk})$ and
$\Psi=\diag(\Psi_{11}, \ldots, \Psi_{kk})$, respectively,
the $1\times 1$ diagonal blocks are
$\Phi_{jj}=\frac{\alpha}{\sqrt{\alpha^2+\beta^2}}$ and
$\Psi_{jj}=\frac{\beta}{\sqrt{\alpha^2+\beta^2}}$,
 corresponding to some real  eigenvalue $(\alpha, \beta)$ (with
 $\beta = 0$ indicating a classical infinite eigenvalue). The
$2\times 2$ diagonal blocks  corresponds to some complex conjugate eigenvalues
$\{(\alpha, \beta), (\bar{\alpha}, \bar{\beta})\}$ $(\alpha\beta\neq0)$ 
with $\Phi_{jj}=I_2$,
$\Psi_{jj}=D_{\delta}(\sigma+i\tau)$ if $|\alpha|\geq|\beta|$, or
$\Phi_{jj}=D_{\delta}(\tilde{\sigma}+i\tilde{\tau})$, $\Psi_{jj}=I_2$
 if $|\alpha|<|\beta|$,  where
$\sigma=\frac{\Re(\bar{\alpha}\beta)}{|\alpha|^2}$, $\tau=\frac{\Im(\bar{\alpha}\beta)}{|\alpha|^2}$,
$\tilde{\sigma}=\frac{\Re(\bar{\beta}\alpha)}{|\beta|^2}$, $\tilde{\tau}=\frac{\Im(\bar{\beta}\alpha)}{|\beta|^2}$,
$0\neq \delta\in\mathbb{R}$.
\end{lem}

\begin{pf}
For the infinite pole $(1, 0)$, let $v\in\mathbb{R}^n$ be the vector satisfying $Ev=0$ and $\|v\|_2=1$.
Then it holds that $Av\neq0$ since  $(A, E)$ is regular. Define $u\equiv\frac{1}{\|Av\|_2^2}Av$ and
denote  the orthonormal basis of the orthogonal complement subspaces  of
$\mathcal{R}(v)$ and $\mathcal{R}(u)$ by $V_{\bot}$ and $U_{\bot}$, respectively. With $U\equiv\begin{row}u& U_{\bot}\end{row}$,
$V\equiv\begin{row}v& V_{\bot}\end{row}$,
we obviously have
\[
	U^\top AV=\begin{bmatrix}1& *\\ 0 & A_1\end{bmatrix},\quad
U^\top EV=\begin{bmatrix}0& *\\ 0 & E_1\end{bmatrix}.
\]

For the real normalized eigenvector $v$
corresponding to the finite real eigenvalue $(\alpha, \beta)$ with  $\beta\neq 0$,
we have $\beta Av=\alpha Ev$ with $Ev\neq0$.  With $u\equiv\frac{\beta}{\sqrt{\alpha^2+\beta^2}\|Ev\|_2^2}Ev$, construct $V_{\bot}, U_{\bot}$ and the orthogonal $V,U$ similarly.
It then follows that 
\[
	U^\top AV=\begin{bmatrix}\frac{\alpha}{\sqrt{\alpha^2+\beta^2}} & *\\ 0&A_2 \end{bmatrix},\quad
U^\top EV=\begin{bmatrix}\frac{\beta }{\sqrt{\alpha^2+\beta^2}} & *\\ 0&E_2 \end{bmatrix}.
\]

For the complex eigenvector $v=x+iy$ $(x, y\in\mathbb{R}^n)$ corresponding to the eigenvalue
$(\alpha, \beta)\in\mathbb{C}\times \mathbb{C}$ $(\alpha\beta\neq 0)$, it follows from $\beta Av=\alpha Ev$
that $\begin{row}Ax& Ay\end{row}D_0(\beta)=\begin{row}Ex& Ey\end{row}D_0(\alpha)$.
Clearly, $x$ and $y$ are linearly independent, else $(\alpha, \beta)\not\in\mathbb{C}\times \mathbb{C}$.
Let $\tilde{x}$ and $\tilde{y}$ be the result of the Jacobi orthogonal transformation on $x$ and $y$ \cite{GCQX},
i.e., $\begin{row}\tilde{x}& \tilde{y} \end{row}=
\begin{row}x& y \end{row}\left[\begin{smallmatrix}c&s\\-s&c\end{smallmatrix}\right]$ (with $c$ and $s$ being the cosine and sine of some angle).
Denote $v_1\equiv\tilde{x}/\|\tilde{x}\|_2$, $v_2\equiv\tilde{y}/\|\tilde{y}\|_2$, then we have
$A\begin{row}v_1& v_2\end{row}D_{\delta}(\beta)=
E\begin{row}v_1& v_2\end{row}D_{\delta}(\alpha)$
with $\delta=\|\tilde{x}\|_2/\|\tilde{y}\|_2$.  Let $V_{\bot}\in\mathbb{R}^{n\times (n-2)}$ be the matrix
satisfying $V_{\bot}^{\top}v_1=V_{\bot}^{\top}v_2=0$ and $V_{\bot}^{\top}V_{\bot}=I_{n-2}$,
then $V=\begin{row}v_1& v_2&  V_{\bot}\end{row}$ is  orthogonal.
Obviously, $1\leq\dim(\mathcal{R}(\begin{row}Av_1& Av_2\end{row}))\leq2$.
Define $U_{\bot}\in\mathbb{R}^{n\times (n-2)}$ as the matrix
which satisfies $U_{\bot}^\top\begin{row}Av_1& Av_2\end{row}=0$
and $U_{\bot}^{\top}U_{\bot}=I_{n-2}$, and let  $\{w_1, w_2\}$ be a  real orthonormal basis of
the orthogonal complement subspace of $\mathcal{R}(U_{\bot})$.
If $|\alpha|\geq|\beta|$, define
\[
\begin{bmatrix}u_1& u_2\end{bmatrix}\equiv\begin{bmatrix}w_1& w_2\end{bmatrix}
\begin{bmatrix}w_1^\top Av_1& w_2^\top Av_1\\ w_1^\top Av_2& w_2^\top Av_2\end{bmatrix}^{-1},
\]
otherwise
\[
\begin{bmatrix}u_1& u_2\end{bmatrix}\equiv\begin{bmatrix}w_1& w_2\end{bmatrix}
\begin{bmatrix}w_1^\top Ev_1& w_2^\top Ev_1\\ w_1^\top Ev_2&  w_2^\top Ev_2\end{bmatrix}^{-1}.
\]
Denote the nonsingular $U=\begin{row}u_1& u_2&  U_{\bot}\end{row}$, we have
\[
	U^\top AV=\begin{bmatrix}I_2&*\\0&A_3\end{bmatrix},\quad
U^\top EV=\begin{bmatrix}D_{\delta}(\sigma+i\tau)&*\\0&E_3\end{bmatrix}
\]
if $|\alpha|\geq|\beta|$,
or
\[
	U^\top AV=\begin{bmatrix}D_{\delta}(\tilde{\sigma}+i\tilde{\tau})&*\\0&A_3\end{bmatrix},\quad
U^\top EV=\begin{bmatrix}I_2&*\\0&E_3\end{bmatrix}
\]
if $|\alpha|<|\beta|$.

Repeat the above process on $(A_1, E_1)$, $(A_2, E_2)$ or $(A_3, E_3)$, we eventually  obtain the result.
\end{pf}

\begin{note}
Suppose that  $\ind_{\infty}(A,E)\leq 1$, or the infinite eigenvalues of $(A, E)$ are
semi-simple, then when all the diagonal block parts corresponding to the infinite eigenvalues are collected together,
we have $T_{pq}=0$ for $p=j, \ldots, j+l$ and $q=p+1, \ldots, j+l$.
\end{note}

When $\Psi_{jj}=D_{\delta_j}(\sigma_j+i\tau_j)$ and
$\Phi_{ll}=D_{\delta_l}(\tilde{\sigma}_l+i\tilde{\tau}_l)$, let
$\Psi_{jj}=W^{*}_{jj}\begin{bsmallmatrix}\sigma_j+i\tau_j & \zeta_j\\ & \sigma_j-i\tau_j\end{bsmallmatrix}W_{jj}$
and  $\Phi_{ll}=\widetilde{W}^{*}_{ll}\begin{bsmallmatrix}\tilde{\sigma}_l+i\tilde{\tau}_l & \tilde{\zeta}_l\\ & \tilde{\sigma}_l-i\tilde{\tau}_l\end{bsmallmatrix}\widetilde{W}_{ll}$
be the Schur decompositions of $\Phi_{jj}$ and $\Psi_{ll}$, respectively.
Direct calculations show that $\zeta_j=(\delta_j-\frac{1}{\delta_j})\tau_j$ and
$\tilde{\zeta}_l=(\delta_l-\frac{1}{\delta_l})\tilde{\tau}_l$.
Now define $Z\equiv\diag(Z_{11}, \ldots, Z_{kk})\in\mathbb{C}^{n\times n}$,
$D\equiv\diag(D_{11}, \ldots, D_{kk})\in\mathbb{C}^{n\times n}$,
where (i) $Z_{jj}=1$, $D_{jj}=1$ if the size of $\Phi_{jj}$ equals to $1$;
(ii) $Z_{jj}=W_{jj}$, $D_{jj}=(1+\sigma_j^2+\tau_j^2)^{-1/2}I_2$ if $\Phi_{jj}=I_2$; or (iii)
$Z_{jj}=\widetilde{W}_{jj}$, $D_{jj}=(1+\tilde{\sigma}_j^2+\tilde{\tau}_j^2)^{-1/2}I_2$ if $\Psi_{jj}=I_2$.
It then holds that  $DZSZ^*$ and $DZTZ^*$ are both upper triangular with the diagonal elements satisfying
$|(DZSZ^*)_{jj}|^2+ |(DZTZ^*)_{jj}|^2=1$ for $j=1, \ldots, n$.

By writing $S$ and $T$ in partitioned form, i.e.,
\begin{align*}
S=\begin{bmatrix}\Phi_{11}&S_{12}&\cdots&S_{1k}\\&\Phi_{22}&\cdots&S_{2k}\\
&&\ddots&\vdots\\&&&\Phi_{kk} \end{bmatrix}, & \quad
T=\begin{bmatrix}\Psi_{11}&T_{12}&\cdots&T_{1k}\\&\Psi_{22}&\cdots&T_{2k}\\
&&\ddots&\vdots\\&&&\Psi_{kk} \end{bmatrix},
\end{align*}
it follows from the definition of $Z$ and $D$ that
\begin{align*}
&\quad \|\offdiag(DZSZ^*)\|_F^2 + \|\offdiag(DZTZ^*)\|_F^2 \\
&=\sideset{}{_{1\times 1 \text{blocks}}}\sum
\sum_{l>j}( |S_{jl}|^2 + |T_{jl}|^2) \ + \\
& 
\sideset{}{_{2\times 2 \text{blocks}}^{|\alpha_j|\geq|\beta_j|}}\sum
 \frac{\sum_{l>j}(\|S_{jl}\|_F^2+\|T_{jl}\|_F^2) +(\delta_j^2-1)^2\delta_j^{-2}\tau_j^2}{1+\sigma_j^2+\tau_j^2}+  \\
&  
\sideset{}{_{2\times 2 \text{blocks}}^{|\alpha_j|<|\beta_j|}}\sum
  \frac{\sum_{l>j}(\|S_{jl}\|_F^2 + \|T_{jl}\|_F^2 )+(\delta_j^2-1)^2\delta_j^{-2}\tilde{\tau}_j^2}{1+\tilde{\sigma}_j^2+\tilde{\tau}_j^2}.
\end{align*}
Furthermore, for those $2\times 2$ blocks $D_{\delta}(\sigma_j+i\tau_j)$
corresponding to $|\alpha_j|\geq |\beta_j|$, it holds that
\[
	\sigma_j^2+\tau_j^2=\frac{|\bar{\alpha}_j\beta_j|^2}{|\alpha_j|^4}=\frac{|\beta_j|^2}{|\alpha_j|^2}\leq1,
\]
 leading  to $1\leq 1+\sigma_j^2+\tau_j^2\leq2$. Analogously, we also have
$1\leq 1+\tilde{\sigma}_j^2+\tilde{\tau}_j^2\leq2$.
Now denote
\begin{multline}\label{dep}
 \Delta_F^2(A,E)
\equiv \|S-\Phi\|_F^2+\|T-\Psi\|_F^2\; + \\
\sideset{}{_{2\times 2 \text{blocks}}^{|\alpha_j|\ge|\beta_j|}}\sum
\tau_j^2\left(\delta_j-\frac{1}{\delta_j}\right)^2+ \\
\sideset{}{_{2\times 2 \text{blocks}}^{|\alpha_j|<|\beta_j|}}\sum
\tilde{\tau}_j^2\left(\delta_j-\frac{1}{\delta_j}\right)^2,
\end{multline}
we have
\begin{multline*}
	\frac{1}{2}\Delta_F^2(A,E)\\
	\leq \|\offdiag(DZSZ^*)\|_F^2 + \|\offdiag(DZTZ^*)\|_F^2\\
	\leq \Delta_F^2(A,E).
\end{multline*}
From the Henrici-type theorem \cite{ESun} for the sensitivity of generalized eigenvalues, we know from Lemma~\ref{Theorem2.1} that
$\Delta_F^2(A+BF, E+BG)$ will be an appropriate robust measure for the closed-loop system or
the corresponding \RPAPDSF.

Next we quote the following solvability result for the \RPAPDSF.
\begin{lem}(\cite{BMN})\label{Theorem2.2}
For the S-controllable descriptor system $(E, A, B)$, a positive integer  $r$ satisfying
$q-m\leq r\leq q \triangleq \rank(\begin{row}E& B\end{row})$, and
an  arbitrary set
$\mathfrak{L}=\{(\alpha_{1}, 0), \ldots, (\alpha_{n-r}, 0), (\alpha_{n-r+1}, \beta_{n-r+1}), \ldots, (\alpha_n, \beta_n)\}$
with  $\{(\alpha_{n-r+1}, \beta_{n-r+1}), \ldots, (\alpha_n, \beta_n)\}$ being self-conjugate and
$\beta_{j}\neq 0$ ($j=n-r+1, \ldots, n$),
there exist $G,F\in\mathbb{R}^{m\times n}$ which
solve the \RPAPDSF.
\end{lem}

For  such $G, F\in\mathbb{R}^{m\times n}$ in  Lemma~\ref{Theorem2.2},
it follows from Lemma~\ref{Theorem2.1}  that
there always exist a nonsingular matrix $X_{G,F}\in\mathbb{R}^{n\times n}$ and
an orthogonal matrix $P_{G,F}\in\mathbb{R}^{n\times n}$ such that
\begin{align}\label{eq-A+BF}
X_{G,F}^{-1}(A+BF)P_{G,F}=S=&
\kbordermatrix{&n-r &r\\
n-r &	S_{11} & S_{12} \\
r   &		   & S_{22} },
\\
\label{eq-E+BG}
X_{G,F}^{-1}(E+BG)P_{G,F}=T=&
\kbordermatrix{&n-r &r\\
n-r &	0	   & T_{12} \\
r   &		   & T_{22} },
\end{align}
where all the diagonal elements in $S_{11}$  are $1$ and
$\lambda(S_{22}, T_{22})=\{ (\alpha_{n-r+1}, \beta_{n-r+1}), \ldots, (\alpha_n, \beta_n)\}$.
The choice of $T_{11} = 0$ in \eqref{eq-E+BG} is  justified  in the Note above.
We shall utilize the decomposition in
\eqref{eq-A+BF} and \eqref{eq-E+BG} to solve the  \RPAPDSF.

Let  $P_{G,F}=\begin{row}p_1 & \cdots &  p_n\end{row}$,
$X_{G,F}=\begin{row}x_1 &  \cdots & x_n\end{row}$, and
define  $N_S\equiv\begin{row}\breve{v}_{1,S}& \cdots&  \breve{v}_{n,S}\end{row}$ and
$N_{T}\equiv\begin{row}\breve{v}_{1,T}& \cdots&  \breve{v}_{n,T}\end{row}$  as the strictly
upper quasi-triangular parts of $S$ and $T$, respectively. In other words, we have
$\breve{v}_{j,S}=\begin{row}v_{j,S}^\top & 0\end{row}^\top$ with $v_{j,S}\in\mathbb{R}^{j-1}$ or
$\mathbb{R}^{j-2}$ ($j=2,  \cdots, n$), and
$\breve{v}_{j,T}=\begin{row}v_{j,T}^\top & 0\end{row}^\top$ with $v_{j,T}\in\mathbb{R}^{j-1}$ or
$\mathbb{R}^{j-2}$ ($j=n-r+1,  \cdots, n$).
We  collect the first $j$ columns of $P_{G,F}$ and $X_{G,F}$ in $P_{j}=\begin{row}p_1& \cdots &  p_{j}\end{row}$ and $X_{j}=\begin{row}x_1& \cdots &  x_{j}\end{row}$, respectively.

\section{Solving the RPA-PDSF via the generalized real Schur form}\label{section3}

\subsection{The parametric solution}\label{subsection3.1}
Without loss of generality,  assume that {\em $B$ is of full column rank}.
Let $B=Q\begin{row}R^\top &  0\end{row}^\top
=Q_1R$
be the QR factorization  of $B$,  where $Q=\begin{row}Q_1& Q_2\end{row}\in\mathbb{R}^{n\times n}$ is orthogonal,
$Q_1\in \mathbb{R}^{n\times m}$, and $R\in\mathbb{R}^{m\times m}$ is nonsingular and upper triangular.
Substituting the  QR factorization  of $B$ into \eqref{eq-A+BF} and \eqref{eq-E+BG},
we deduce that
\begin{equation}\label{GF}
\left\{\begin{aligned}
F&=R^{-1}Q_1^\top(X_{G,F}SP_{G,F}^{\top}-A),\\
G&=R^{-1}Q_1^\top(X_{G,F}TP_{G,F}^{\top}-E),\\
\end{aligned}\right.
\end{equation}
where
\begin{equation}\label{condition}
\left\{\begin{aligned}
Q_2^\top(AP_{G,F}-X_{G,F}S)&=0,\\
Q_2^\top(EP_{G,F}-X_{G,F}T)&=0.
\end{aligned}\right.
\end{equation}
Consequently, once we obtain an orthogonal $P_{G,F}$, a nonsingular $X_{G,F}$
and a pair of upper  quasi-triangular $S$ and $T$ satisfying \eqref{condition},
a solution $(G,F)$ to the pole assignment problem can be computed  through \eqref{GF}.

\subsection{Assigning the infinite pole $(\alpha_{j}, 0)$}\label{subsection3.2}

Provided that there exist some infinite poles in $\mathfrak{L}$,
we shall show how to assign all infinite poles $(\alpha_j, 0)$ for $j=1, \ldots, n-r$.
Suppose  $j-1$ infinite poles ($j\geq 1$) have already been placed, suggesting that
$P_{j-1}$, $Q_2^\top X_{j-1}$  and  the $(j-1)\times (j-1)$
leading principal submatrix $S_{j-1}$  of $S$ have been acquired.
We are going to compute the $j$-th column of $P_{G,F}$, $Q_2^\top X_{G,F}$ and $S$
when assigning the infinite pole $(\alpha_{j}, 0)$.
We emphasize that $Q_2^\top X_{j-1}$ (not $X_{j-1}$) is known, in the computation of $Q_2^\top X_j$ (not $X_j$).

It is simple  to show that $\rank(\begin{row}E& B\end{row})=m+\rank(Q_2^\top E)$, indicating that
$l\triangleq\dim(\mathcal{N}(Q_2^\top E)) \geq n-r$. Let the columns of $Z\in\mathbb{R}^{n \times l}$ be an orthonormal
basis of $\mathcal{N}(Q_2^\top E)$, and the second equation in   \eqref{condition} implies that $Q_2^\top EP_{n-r}=0$.
It follows that $P_{j-1}=ZW_{j-1}$ for some $W_{j-1}\in\mathbb{R}^{l\times (j-1)}$
with $W_{j-1}^\top W_{j-1}=I_{j-1}$  and $p_{j}=Zw_{j}$ with a normalized $w_j\in\mathbb{R}^{l}$ to be specified.
From the $j$-th column of the first equation in \eqref{condition}, noting that the $j$-th diagonal
element of $S$ is $1$ and $P$ is orthogonal,  we have
\[
Q_2^\top AZw_{j}=Q_2^\top X_{j-1}v_{j,S} + Q_2^\top x_{j}, \quad W_{j-1}^\top w_j=0.
\]
Recalling the definition of $\Delta_F^2(A,E)$ in \eqref{dep},  we should solve the following optimization problem:
\begin{subequations}\label{eqinfinite-opt}
\begin{align}
\min&\quad \|v_{j,S}\|_2^2\label{eqinfinite-obj}\\
\text{s.t.}&\quad  M_{j-1}\begin{bmatrix}w_{j}^\top&\ x_{j}^\top&\ v_{j,S}^\top\end{bmatrix}^\top=0,\label{eqinfinite-conrtrain}
\end{align}
\end{subequations}
where
\[
M_{j-1}=\begin{bmatrix}
Q_2^\top AZ & -Q_2^\top & -Q_2^\top X_{j-1}\\
W_{j-1}^\top & 0& 0 \end{bmatrix}.
\]

Theorem \ref{theorem-infinite} in Section \ref{subsection3.6} demonstrates that $M_{j-1}$ is of full row rank and there
exists some vector $\begin{bmatrix}w^\top & x^\top&  v_S^\top\end{bmatrix}^\top\in\mathcal{N}(M_{j-1})$
with $w\in\mathbb{R}^l$, $x\in\mathbb{R}^n$, $v_S\in\mathbb{R}^{j-1}$,  such that $w\neq 0$.
Denote $W_{\bot}\in\mathbb{R}^{l\times (l-j+1)}$ satisfying $W_{j-1}^\top W_{\bot}=0$ and
$W_{\bot}^\top W_{\bot}=I_{l-j+1}$, then the columns of
\[
\begin{bmatrix}W_{\bot} &  0 & 0\\ Q_2Q_2^\top AZW_{\bot} &  -Q_2Q_2^\top X_{j-1}& Q_1\\  0& I_{j-1} & 0\end{bmatrix}
\]
form a basis of $\mathcal{N}(M_{j-1})$, suggesting
$w_{j}=W_{\bot}u_j$ for some normalized $u_j\in\mathbb{R}^{l-j+1}$ 
and $x_{j}=Q_1y_j+Q_2Q_2^\top AZw_j-Q_2(Q_2^\top X_{j-1})v_{j,S}$ for some $y_j\in\mathbb{R}^m$.
Accordingly, $p_j=ZW_{\bot}u_j$, where $u_j$ would take any arbitrary unit vector.
Apparently, the optimization sub-problem  \eqref{eqinfinite-opt} obtains its minimum when $v_{j,S}=0$,
leading to $x_j=Q_1y_j+Q_2Q_2^\top Ap_j$, with $y_j$ to be specified in Section~\ref{subsection3.5}.
Consequently,  we obtain $Q_2^\top x_j=Q_2^\top Ap_j$ which is sufficient
for the assigning process to continue. Note in the definition of $M_{j-1}$ and the assigning procedure for
finite poles later that only $Q_2^\top X_{j-1}$ (not $X_{j-1}$) is required.

The procedure for assigning all infinite poles  is summarized in  Algorithm \ref{algorithm1},
where we use $\Xi_{n-r} \equiv Q_2^\top X_{n-r}$.

\begin{algorithm}
\caption{\ Assigning all infinite poles $\{(\alpha_1, 0), \ldots, (\alpha_{n-r},0)\}$}
\begin{algorithmic}[1]\label{algorithm1}
\REQUIRE ~~\\
$A$, $E$ and $Q_2$.
\ENSURE ~~\\
Orthogonal $P_{n-r}\in\mathbb{R}^{n\times (n-r)}$, upper triangular $S_{n-r}\in\mathbb{R}^{(n-r)\times (n-r)}$
and $\Xi_{n-r}\in\mathbb{R}^{(n-m)\times (n-r)}$.
\STATE Find $Z\in\mathbb{R}^{n\times l}$, whose columns form an orthonormal basis of $\mathcal{N}(Q_2^\top E)$.
\STATE Set $P_{n-r}=ZW_{n-r}$ with $W_{n-r}\in\mathbb{R}^{l\times (n-r)}$ an arbitrary orthogonal matrix.
\STATE Set $S_{n-r}=I_{n-r}$.
\STATE Set $\Xi_{n-r}=Q_2^\top AP_{n-r}$.
\end{algorithmic}
\end{algorithm}

\smallskip

Provided that all infinite poles and some  finite poles have already been  assigned,
where the complex conjugate poles are placed together.
Consider the already acquired  $P_j$, $Q_2^\top X_j$  and the  $j\times j$ leading
principal submatrix $S_j$ and $T_j$ of $S$ and $T$, respectively, they satisfy
\[
Q_2^\top AP_j=Q_2^\top X_jS_j,   \quad Q_2^\top EP_j=Q_2^\top X_jT_j \quad (j\geq n-r).
\]
The details of the pole assignment for the
finite real pole $(\alpha_{j+1}, \beta_{j+1})$ and the finite
complex conjugate poles $\{(\alpha_{j+1}, \beta_{j+1}), (\bar{\alpha}_{j+1}, \bar{\beta}_{j+1})\}$
will be presented in Sections~\ref{subsection-real} and \ref{subsection-complex}, respectively.
The $(j+1)$-th column, or the $(j+1)$-th  and $(j+2)$-th columns, of $P_{G,F}$, $Q_2^\top X_{G,F}$, $S$ and $T$
are computed in the assignment process.

\subsection{Assigning the  finite real pole $(\alpha_{j+1}, \beta_{j+1})$}\label{subsection-real}

Let $(\alpha_{j+1}, \beta_{j+1})\in\mathbb{R}\times \mathbb{R}$, then the $(j+1)$-th diagonal elements of $S$ and
$T$ are $\frac{\alpha_{j+1}}{\sqrt{|\alpha_{j+1}|^2+|\beta_{j+1}|^2}}$ and
$\frac{\beta_{j+1}}{\sqrt{|\alpha_{j+1}|^2+|\beta_{j+1}|^2}}$, respectively.
The $(j+1)$-th columns in \eqref{condition} are
\[
	\left\{\begin{aligned}
Q_2^\top Ap_{j+1} -Q_2^\top X_{j}v_{j+1,S}- \frac{\alpha_{j+1}Q_2^\top x_{j+1}}{\sqrt{|\alpha_{j+1}|^2+|\beta_{j+1}|^2}}&=0,\\
Q_2^\top Ep_{j+1} -Q_2^\top X_{j}v_{j+1,T}- \frac{\beta_{j+1}Q_2^\top x_{j+1}}{\sqrt{|\alpha_{j+1}|^2+|\beta_{j+1}|^2}}&=0,
\end{aligned}
\right.
\]
which are the conditions $p_{j+1}$, $x_{j+1}$, $v_{j+1,S}$ and $v_{j+1,T}$ have to meet.
From the definition of $\Delta_F^2(A,E)$ in \eqref{dep} and the orthogonality
of $P_{G,F}$,
it is then natural to consider  the optimization sub-problem:
\begin{subequations}\label{eqreal-opt}
\begin{align}
\min_{\|p_{j+1}\|_2=1}&\|v_{j+1,S}\|_2^2+\|v_{j+1,T}\|_2^2\label{eqreal-obj}\\
\text{s.t.}\quad & M_j\begin{bmatrix}p_{j+1}^\top& x_{j+1}^\top& v_{j+1,S}^\top& v_{j+1,T}^\top\end{bmatrix}^\top=0,\label{eqreal-conrtrain}
\end{align}
\end{subequations}
where
\begin{align}\label{realM}
M_j=\begin{bmatrix}Q_2^\top A&-\frac{\alpha_{j+1}Q_2^\top}{\sqrt{|\alpha_{j+1}|^2+|\beta_{j+1}|^2}}&-Q_2^\top X_j&0\\
Q_2^\top E&-\frac{\beta_{j+1}Q_2^\top}{\sqrt{|\alpha_{j+1}|^2+|\beta_{j+1}|^2}}&0&-Q_2^\top X_j\\
P_j^\top&0&0&0\end{bmatrix}.
\end{align}

Theorem \ref{theorem-real} in Section~\ref{subsection3.6} shows that $\dim(\mathcal{N}(M_j))=2m+j$
and  there exists
$\begin{bmatrix}p_{j+1}^\top& x_{j+1}^\top& v_{j+1,S}^\top& v_{j+1,T}^\top\end{bmatrix}^\top\in\mathcal{N}(M_j)$
such that $p_{j+1}\neq0$, which guarantees the solvability of \eqref{eqreal-opt}. 
Next, we shall consider the solution of the optimization sub-problem \eqref{eqreal-opt}.

\begin{enumerate}[label={\bf Case \roman*},leftmargin=0em,itemindent=5em]
\item ($|\alpha_{j+1}| \geq |\beta_{j+1}|$) 
  Define
\[
\widetilde{M}_j=\begin{bgapmatrix}-\frac{\alpha_{j+1}Q_2^\top}{\sqrt{|\alpha_{j+1}|^2+|\beta_{j+1}|^2}} &\vrule&
	Q_2^\top A  &   -Q_2^\top X_j &   0  \\[1.5ex]
\hline
0 & \vrule& &\widetilde{M}_{2,j}&\end{bgapmatrix}
\]
with
\begin{align}\label{M_2j-real}
	\widetilde{M}_{2,j}=\begin{bmatrix}
Q_2^\top \left(E-\frac{\beta_{j+1}}{\alpha_{j+1}}A\right)&  \frac{\beta_{j+1}}{\alpha_{j+1}}Q_2^\top X_j&  -Q_2^\top X_j\\
P_j^\top& 0& 0\end{bmatrix},
\end{align}
then \eqref{eqreal-conrtrain} is equivalent to $\widetilde{M}_j\begin{bmatrix}x_{j+1}^\top&p_{j+1}^\top&v_{j+1,S}^\top&v_{j+1,T}^\top\end{bmatrix}^\top=0$.
Equivalently, we have
\[
	\frac{\alpha_{j+1} Q_2^\top x_{j+1}}{\sqrt{|\alpha_{j+1}|^2+|\beta_{j+1}|^2}} =Q_2^\top Ap_{j+1}-Q_2^\top X_jv_{j+1,S}
\]
and $\widetilde{M}_{2,j}\begin{bmatrix}p_{j+1}^\top&v_{j+1,S}^\top&v_{j+1,T}^\top\end{bmatrix}^\top=0$.
Evidently, $\widetilde{M}_{2,j}$ is of full row rank, implying that $\dim(\mathcal{N}(\widetilde{M}_{2,j}))=m+j$.
Let the columns of $\begin{bmatrix}Z_1^\top&Z_3^\top&Z_4^\top\end{bmatrix}^\top$
be  an orthonormal basis of $\mathcal{N}(\widetilde{M}_{2,j})$, where $Z_1\in\mathbb{R}^{n\times (m+j)}$,
$Z_3\in\mathbb{R}^{j\times (m+j)}$, $Z_4\in\mathbb{R}^{j\times (m+j)}$,
then the columns of
\[
\begin{bmatrix}
	0 & Q_1^\top & 0& 0\\
Z_1^\top & \frac{\sqrt{\alpha_{j+1}^2+\beta_{j+1}^2}}{\alpha_{j+1}}(AZ_1-X_jZ_3)^\top Q_2Q_2^\top&Z_3^\top&Z_4^\top\\
\end{bmatrix}^\top
\]
form a basis of $\mathcal{N}(M_j)$.

Consequently, the constrained optimization sub-problem \eqref{eqreal-opt} is reduced
to the following quadratic optimization problem:
\begin{align}\label{eqreal-opt1}
\min_{\|Z_1 u\|_2=1} u^\top(Z_3^\top Z_3+ Z_4^\top Z_4)u,
\end{align}
where $u\in\mathbb{R}^{m+j}$. Furthermore, since $Z_1^\top Z_1 + Z_3^\top Z_3 + Z_4^\top Z_4=I_{m+j}$,
\eqref{eqreal-opt1} is further reduced to 
\begin{align}\label{eqreal-opt2}
\min_{\|Z_1 u\|_2=1} u^\top u,
\end{align}
which is attained by $u$ being an eigenvector of $Z_1^\top Z_1$
corresponding to its greatest eigenvalue with $u^\top Z_{1}^\top Z_{1}u=1$.
Once such $u$ is obtained, $p_{j+1}$, $v_{j+1,S}$ and $v_{j+1,T}$ can be retrieved by
\begin{align}\label{pvv-real}
p_{j+1}=Z_1u, \quad v_{j+1,S}=Z_3u, \quad v_{j+1,T}=Z_4u.
\end{align}

We also have
\[ x_{j+1}=Q_1y_{j+1}+ \tfrac{\sqrt{\alpha_{j+1}^2+\beta_{j+1}^2}}{\alpha_{j+1}}
Q_2 \left[ Q_2^\top Ap_{j+1} - (Q_2^\top X_j)v_{j+1,S}\right] \]
for some $y_{j+1}\in\mathbb{R}^m$ to be determined. 
Clearly,
\[
Q_2^\top x_{j+1}=\tfrac{\sqrt{\alpha_{j+1}^2+\beta_{j+1}^2}}{\alpha_{j+1}}
 \left[ Q_2^\top Ap_{j+1} - (Q_2^\top X_j)v_{j+1,S}\right],
\]
which can be computed and added
to  $Q_2^\top X_{j+1}=\begin{bmatrix}Q_2^\top X_j&Q_2^\top x_{j+1}\end{bmatrix}$.
Recall from the definition of $M_{j}$ in \eqref{realM},
it is $Q_2^\top X_{j+1}$, rather than $X_{j+1}$, that is required when assigning the finite real
pole $(\alpha_{j+2}, \beta_{j+2})$, without requiring $y_{j+1}$. Similar comments hold
for the case of $(\alpha_{j+2}, \beta_{j+2})\in\mathbb{C}\times \mathbb{C}$, which will be discussed later.
The choice of $y_{j+1}$ will be discussed in Section~\ref{subsection3.5}.

\item ($|\alpha_{j+1}| < |\beta_{j+1}|$) 
Analogously to {\bf Case \romannumeral1}, let  the columns of
$\begin{bmatrix}Z_1^\top &Z_3^\top &Z_4^\top\end{bmatrix}^\top$,
 where $Z_1\in\mathbb{R}^{n\times (m+j)}$,
$Z_3, Z_4 \in\mathbb{R}^{j\times (m+j)}$, 
form an orthonormal  basis of
\begin{align}\label{M_2j-real2}
	\widetilde{M}_{2,j}=\begin{bmatrix}
Q_2^\top \left(A-\frac{\alpha_{j+1}}{\beta_{j+1}}E \right)& -Q_2^\top X_j & \frac{\alpha_{j+1}}{\beta_{j+1}}Q_2^\top X_j \\
P_j^\top& 0& 0\end{bmatrix},
\end{align}
where $\dim(\mathcal{R}(\begin{row}Z_1^\top &Z_3^\top &Z_4^\top\end{row}^\top))=m+j$ is
guaranteed by $\rank(\widetilde{M}_{2,j})=(n-m+j)$.
 Besides, the columns of
 \[
\begin{bmatrix} 0& Q_1^\top  & 0& 0\\
 Z_1^\top &\frac{\sqrt{\alpha_{j+1}^2+\beta_{j+1}^2}}{\beta_{j+1}}(EZ_1-X_jZ_4)^\top Q_2Q_2^\top  &  Z_3^\top&  Z_4^\top\\
\end{bmatrix}^\top
\]
form a  basis of $\mathcal{N}(M_j)$,
 leading to $p_{j+1}=Z_1u$,  $v_{j+1,S}=Z_3u$, $v_{j+1,T}=Z_4u$ and
 \[
 x_{j+1}=Q_1y_{j+1}+\tfrac{\sqrt{\alpha_{j+1}^2+\beta_{j+1}^2}}{\beta_{j+1}}
Q_2 \left[ Q_2^\top Ep_{j+1} - (Q_2^\top X_j)v_{j+1,T} \right],
\]
for some $u\in\mathbb{R}^{m+j}$ and $y_{j+1}\in\mathbb{R}^m$.
The constrained optimization sub-problem \eqref{eqreal-opt} can be solved in the same
manner as in  {\bf Case  \romannumeral1}, with $y_{j+1}$ and $x_{j+1}$ to be specified in Section~\ref{subsection3.5}.
\end{enumerate}

\begin{note}
When there is no  infinite pole, i.e., $r=n$, some minor modifications to our method will be required.
We first place the real finite pole $(\alpha_1, \beta_1)$.
More specifically, there is no contribution
from the first columns of $S$ and $T$ to $\Delta_F^2(A+BF,E+BG)$ and
the optimization sub-problem \eqref{eqreal-opt} is degenerate. We just need to select $p_1$ and $x_1$ from the constraint
\eqref{eqreal-conrtrain}.
Lemma \ref{Lemma3.2} in subsection \ref{subsection3.6}  implies that $M_{0}$ is of full row rank, thus the feasibility of
\eqref{eqreal-conrtrain}.
We can select a normalized $p_1$, then
$Q_2^\top x_1=\alpha_1^{-1} \sqrt{\alpha_1^2+\beta_1^2}Q_2^\top Ap_1$ (for $|\alpha_1|\geq|\beta_1|$)
or  $Q_2^\top x_1=\beta_1^{-1} \sqrt{\alpha_1^2+\beta_1^2} Q_2^\top Ep_1$ (for $|\alpha_1|<|\beta_1|$).
In addition, we have  $x_1=Q_1y_1+Q_2(Q_2^\top x_1)$ with $y_1$ to be chosen (as in Section~\ref{subsection3.5}).
Similar comments hold when we have to assign a complex conjugate pair of finite poles first.
\end{note}

We summarize the assignment of the finite real pole $(\alpha_{j+1}, \beta_{j+1})$ in
Algorithm~\ref{algorithm2}, with $\Xi_j \equiv Q_2^\top X_j$.
\begin{algorithm}
\caption{\ Assigning finite real pole  $(\alpha_{j+1}, \beta_{j+1})$  }
\begin{algorithmic}[1]\label{algorithm2}
\REQUIRE ~~\\
$A$, $E$, $Q_2$, $P_j$, $\Xi_j$, $S_j$, $T_j$ and $(\alpha_{j+1}, \beta_{j+1})\in\mathbb{R}\times \mathbb{R}$.
\ENSURE ~~\\
Orthogonal $P_{j+1}\in\mathbb{R}^{n\times (j+1)}$, upper quasi-triangular $S_{j+1}, T_{j+1}\in\mathbb{R}^{(j+1)\times (j+1)}$
and  $\Xi_{j+1}\in\mathbb{R}^{(n-m)\times (j+1)}$.
\IF {$|\alpha_{j+1}| \geq |\beta_{j+1}|$}
    \STATE Find  $\begin{bmatrix}Z_1^\top &Z_3^\top &Z_4^\top\end{bmatrix}^\top$ with
    $Z_1\in\mathbb{R}^{n\times (m+j)}$, $Z_3, Z_4\in\mathbb{R}^{j\times (m+j)}$,
    whose columns form  an orthonormal basis of $\mathcal{N}(\widetilde{M}_{2,j})$ as defined in \eqref{M_2j-real}.
\ELSE
    \STATE Find $\begin{bmatrix}Z_1^\top &Z_3^\top &Z_4^\top\end{bmatrix}^\top$ with
    $Z_1\in\mathbb{R}^{n\times (m+j)}$, $Z_3, Z_4\in\mathbb{R}^{j\times (m+j)}$,
    whose columns form  an orthonormal basis of $\mathcal{N}(\widetilde{M}_{2,j})$ as defined in \eqref{M_2j-real2}.
\ENDIF
\STATE Solve  the optimization sub-problem \eqref{eqreal-opt2} to get $u\in\mathbb{R}^{m+j}$.
\STATE Set $P_{j+1}=\begin{bmatrix}P_j&p_{j+1}\end{bmatrix}$ with $p_{j+1}=Z_1u$.
\STATE Set
 \[
  S_{j+1}=\begin{bmatrix}S_j& v_{j+1,S}\\&\frac{\alpha_{j+1}}{\sqrt{\alpha_{j+1}^2+\beta_{j+1}^2}}\end{bmatrix}, \quad
  T_{j+1}=\begin{bmatrix}T_j& v_{j+1,T}\\&\frac{\beta_{j+1}}{\sqrt{\alpha_{j+1}^2+\beta_{j+1}^2}}\end{bmatrix}
 \]
 with $v_{j+1,S}=Z_3u$ and $v_{j+1,T}=Z_4u$.
 \STATE Set $\Xi_{j+1}=\begin{bmatrix}\Xi_j&\xi_{j+1}\end{bmatrix}$ with \\
 $\xi_{j+1}=\frac{\sqrt{\alpha_{j+1}^2+\beta_{j+1}^2}}{\alpha_{j+1}}(Q_2^\top Ap_{j+1}-\Xi_jv_{j+1,S})$ \
 (if $|\alpha_{j+1}| \geq |\beta_{j+1}|$), \\
\hspace*{0.57cm} $= \frac{\sqrt{\alpha_{j+1}^2+\beta_{j+1}^2}}{\beta_{j+1}}(Q_2^\top Ep_{j+1}-\Xi_jv_{j+1,T})$ \ (otherwise).
\end{algorithmic}
\end{algorithm}

\subsection{Assigning the finite complex pole $(\alpha_{j+1}, \beta_{j+1})$}\label{subsection-complex}

With the $2\times 2$ diagonal blocks in $S$ and $T$ specified in Lemma \ref{Theorem2.1},
assigning the complex conjugate pair $\{(\alpha_{j+1}, \beta_{j+1}), (\bar{\alpha}_{j+1}, \bar{\beta}_{j+1})\}$,
involves two different cases, when $|\alpha_{j+1}|\geq|\beta_{j+1}|$ or otherwise.

\subsubsection{Situation \uppercase\expandafter{\romannumeral1}} ($|\alpha_{j+1}|\geq|\beta_{j+1}|$)
Under such circumstance, we have $S(j+1\step j+2,j+1\step j+2)=I_2$ and $T(j+1\step j+2,j+1\step j+2)=D_{\delta_{j+1}}(\sigma_{j+1}+i\tau_{j+1})$.
The $(j+1)$-th and  $(j+2)$-th columns of \eqref{condition} can be expanded to
\begin{equation}\label{pj12}
\begin{split}
&Q_2^\top A\begin{row}p_{j+1}&p_{j+2}\end{row}-Q_2^\top X_j\begin{row}v_{j+1,S}&v_{j+2,S}\end{row}
-Q_2^\top\begin{row}x_{j+1}&x_{j+2}\end{row}=0,\\
&Q_2^\top E\begin{row}p_{j+1}&p_{j+2}\end{row}-Q_2^\top X_j\begin{row}v_{j+1,T}&v_{j+2,T}\end{row} \\
& \qquad \qquad -Q_2^\top\begin{row}x_{j+1}&x_{j+2}\end{row}D_{\delta_{j+1}}(\sigma_{j+1}+i\tau_{j+1})=0.
\end{split}
\end{equation}
Let  $\delta_{j+1}=\varsigma_1/\varsigma_2$ with $\varsigma_1, \varsigma_2 \in\mathbb{R}$ and $\varsigma_2\neq0$.
We can verify that
\begin{align*}
&D_{\delta_{j+1}}(\sigma_{j+1}+i\tau_{j+1})\\
=&\frac{1}{2}\begin{bsmallmatrix}\varsigma_1 & \\ & \varsigma_2\end{bsmallmatrix}
\begin{bsmallmatrix}1&1\\i&-i\end{bsmallmatrix} \begin{bsmallmatrix}\sigma_{j+1}+i\tau_{j+1}& \\ & \sigma_{j+1}-i\tau_{j+1}\end{bsmallmatrix}
\begin{bsmallmatrix}1&-i\\1&i\end{bsmallmatrix}\begin{bsmallmatrix}\varsigma_1^{-1} & \\ &\varsigma_2^{-1} \end{bsmallmatrix}.
\end{align*}
Defining  $\tilde{p}_{j+l}=\varsigma_lp_{j+l}$,
$\tilde{x}_{j+l}=\varsigma_lx_{j+l}$, $\tilde{v}_{j+l,S}=\varsigma_lv_{j+l,S}$ and
$\tilde{v}_{j+l,T}=\varsigma_lv_{j+l,T}$ for $l=1, 2$,
\eqref{pj12} is equivalent to
\begin{equation}\label{pj12-equ1}
	\left\{\begin{aligned}
 Q_2^{\top}A(\tilde{p}_{j+1}+i\tilde{p}_{j+2}) - Q_2^\top X_j(\tilde{v}_{j+1,S}+i\tilde{v}_{j+2,S})\qquad\\
  \qquad \qquad \quad -Q_2^\top(\tilde{x}_{j+1}+i\tilde{x}_{j+2})=0,\\
  Q_2^{\top}E(\tilde{p}_{j+1}+i\tilde{p}_{j+2}) - Q_2^\top X_j(\tilde{v}_{j+1,T}+i\tilde{v}_{j+2,T})\qquad\\
  \qquad \qquad \quad -(\sigma_{j+1}+i\tau_{j+1})Q_2^\top(\tilde{x}_{j+1}+i\tilde{x}_{j+2})=0.
\end{aligned}\right.
\end{equation}
Consequently, $\tilde{p}_{j+l}$, $\tilde{x}_{j+l}$, $\tilde{v}_{j+l,S}$ and
$\tilde{v}_{j+l,T}$ ($l=1, 2$) have to be selected satisfying \eqref{pj12-equ1}, in addition to the constraint $\begin{row}\tilde{p}_{j+1}&\tilde{p}_{j+2}\end{row}^\top
\begin{row}\tilde{p}_{j+1}&\tilde{p}_{j+2}\end{row}=\diag(\varsigma_1^2, \varsigma_2^2)$
(so that $\begin{row}p_{j+1}&p_{j+2}\end{row}^\top \begin{row}p_{j+1}&p_{j+2}\end{row}=I_2$).


Recalling the definition $\Delta_F^2(A,E)$ in \eqref{dep},
we then select $(j+1)$-th and  $(j+2)$-th columns of $P_{G,F}$, $X_{G,F}$, $S$ and $T$ while minimizing
their contributions to $\Delta_F^2(A+BF,E+BG)$.
In other words, we solve the optimization sub-problem:
\begin{subequations}\label{eqcom-opt}
\begin{align}
\min_{\begin{subarray}{c}\varsigma_1, \varsigma_2,\\ \tilde{v}_{j+1,S},
\tilde{v}_{j+2,S},\\ \tilde{v}_{j+1,T}, \tilde{v}_{j+2,T}\end{subarray}} &
\sum_{l=1}^2 \left[
\tfrac{\|\tilde{v}_{j+l,S}\|_2^2}{\varsigma_l^2} + \tfrac{\|\tilde{v}_{j+l,T}\|_2^2}{\varsigma_l^2}
\right] +
\tau_{j+1}^2 \left(\tfrac{\varsigma_1}{\varsigma_2} -\tfrac{\varsigma_2}{\varsigma_1} \right)^2
 \label{eqcom-opt-obj}\\
\mbox{s.t.} \quad & [\ \tilde{p}_{j+1}^{\top}+i\tilde{p}_{j+2}^{\top}, \ \ \tilde{x}_{j+1}^{\top}+i\tilde{x}_{j+2}^{\top}, \notag \\
& \tilde{v}_{j+1,S}^{\top}+i\tilde{v}_{j+2,S}^{\top}, \ \ \tilde{v}_{j+1,T}^{\top}+i\tilde{v}_{j+2,T}^{\top}\ ]^{\top}
\in \mathcal{N} (M_{j}), \label{eqcom-opt-constraina}\\
&\begin{row}\tilde{p}_{j+1}&\tilde{p}_{j+2}\end{row}^\top \begin{row}\tilde{p}_{j+1}&\tilde{p}_{j+2}\end{row}=\diag(\varsigma_1^2, \varsigma_2^2),\label{eqcom-opt-constrainb}
\end{align}
\end{subequations}
where
\begin{align}\label{M_j_complex}
 M_{j}=\begin{bmatrix}Q_2^\top A & -Q_2^\top& -Q_2^\top X_j&0\\
 Q_2^\top E  & -(\sigma_{j+1}+i\tau_{j+1})Q_2^\top&0& -Q_2^\top X_j\\
 P_j^\top&0&0&0\end{bmatrix}.
\end{align}

Once a solution to  the optimization sub-problem \eqref{eqcom-opt} is acquired, then the
$(j+1)$-th and $(j+2)$-th columns of $P_{G,F}$, $X_{G,F}$, $S$ and $T$ can be retrieved via normalization: (for  $l=1, 2$)
\begin{align*}
&p_{j+l}=\frac{\tilde{p}_{j+l}}{\|\tilde{p}_{j+l}\|_2},     &\qquad
x_{j+l}=\frac{\tilde{x}_{j+l}}{\|\tilde{p}_{j+l}\|_2}, &\\
&v_{j+l,S}=\frac{\tilde{v}_{j+l,S}}{\|\tilde{p}_{j+l}\|_2}, & \qquad
v_{j+l,T}=\frac{\tilde{v}_{j+l,T}}{\|\tilde{p}_{j+l}\|_2}. &\qquad
\end{align*}

To solve the constrained  optimization sub-problem \eqref{eqcom-opt}, we firstly consider the constraint \eqref{eqcom-opt-constraina}. Analogous to the previous section for finite real poles,
define $\gamma_{j+1} \equiv \sigma_{j+1}+i\tau_{j+1}$,
\begin{align}\label{M_2j-com}
	\widetilde{M}_{2,j}=\begin{bmatrix}
Q_2^\top (E-\gamma_{j+1} A) &  \gamma_{j+1} Q_2^\top X_j  &-Q_2^\top X_j\\
P_j^\top& 0& 0\end{bmatrix},
\\\nonumber
\widetilde{M}_j=\begin{bgapmatrix}-Q_2^\top &\vrule&
Q_2^\top A  &  -Q_2^\top X_j &  0 \\
\hline
0 &\vrule& &\widetilde{M}_{2,j}&\end{bgapmatrix},
\end{align}
then we have $M_{j}\begin{row}z_1^\top &z_2^\top &z_3^\top &z_4^\top\end{row}^\top=0$
if and only if $\widetilde{M}_{j}\begin{row}z_2^\top &z_1^\top &z_3^\top &z_4^\top\end{row}^\top=0$
with $z_1, z_2 \in\mathbb{C}^n$, $z_3, z_4 \in\mathbb{C}^j$.
Furthermore, it follows from Theorem~\ref{theorem-real} in Section~\ref{subsection3.6} that
$M_j$, thus $\widetilde{M}_{2,j}$, are of full row rank, or
$\dim(\mathcal{N}(\widetilde{M}_{2,j}))=m+j$.
Now let the columns of $\begin{bmatrix}Z_1^\top&Z_3^\top&Z_4^\top\end{bmatrix}^\top$
be   an orthonormal basis of $\mathcal{N}(\widetilde{M}_{2,j})$, where $Z_1\in\mathbb{C}^{n\times (m+j)}$,
$Z_3, Z_4\in\mathbb{C}^{j\times (m+j)}$, then the columns of
\[
\begin{bmatrix}0 & Q_1^\top &\quad  0&\quad  0\\
Z_1^\top \qquad& \qquad (AZ_1-X_jZ_3)^\top Q_2Q_2^\top \qquad &\quad Z_3^\top\qquad &\quad Z_4^\top\\
\end{bmatrix}^\top
\]
constitute a   basis of $\mathcal{N}(M_j)$. We can then select
\begin{align*}
\tilde{p}_{j+1}+i\tilde{p}_{j+2}&=Z_1b, \\
\tilde{x}_{j+1}+i\tilde{x}_{j+2}&=Q_1y + Q_2 Q_2^\top(AZ_1-X_jZ_3)b, \\
\tilde{v}_{j+1,S}+i\tilde{v}_{j+2,S}&=Z_3b,\\
\tilde{v}_{j+1,T}+i\tilde{v}_{j+2,T}&=Z_4b,
\end{align*}
for some $0 \neq b\in\mathbb{C}^{m+j}$ and $y\in\mathbb{C}^{m}$.
Accordingly, the optimization problem \eqref{eqcom-opt} is reduced
to choosing some suitable nonzero $b\in\mathbb{C}^{m+j}$ such that $\tilde{p}_{j+1}^\top \tilde{p}_{j+2}=0$, while
minimizing the objective function formulated in \eqref{eqcom-opt-obj}.
It is worthwhile to point out that Theorem \ref{theorem-real} in Section~\ref{subsection3.6} guarantees that
$Z_1\neq0$ and there exist some nontrivial $b\in\mathbb{C}^{m+j}$ such that
$\tilde{p}_{j+1}$ and $\tilde{p}_{j+2}$ are linearly independent,
which is necessary for $\{\tilde{p}_{j+1}, \tilde{p}_{j+2}\}$ to be orthogonal.
In what follows, we consider how $b\in\mathbb{C}^{m+j}$ is selected, in two distinct cases.

Denote  $Z_1=U_{Z_1} \Sigma_{Z_1} V_{Z_1}^{\ast}$  as the Singular Value Decomposition (SVD) of $Z_1$, where
its nonzero singular values $\nu_j$ ($j=1,\cdots,r_{Z_1}$) satisfy
$\nu_1 \geq \nu_2 \geq \ldots \geq \nu_{r_{Z_1}}>0$.
Note that  $Z_1^*Z_1 + Z_3^*Z_3 + Z_4^*Z_4 = I_{m+j}$,
implying that $Z_3^*Z_3 + Z_4^*Z_4= V_{Z_1}(I_{m+j}-\Sigma_{Z_1}^\top\Sigma_{Z_1})V_{Z_1}^*$.

\begin{enumerate}[label={\bf Case \roman*},leftmargin=0em,itemindent=5em]
  \item ($\rank(Z_1)=1$) In this case, there exists a unique nonzero singular
  value $\nu_1$ for $Z_1$, with the corresponding left-singular vector $\psi_1$.   Define
  \begin{align*}
  &\mathcal{N}_1(\widetilde{M}_{2,j}) \equiv \{  \begin{bmatrix}\psi_1^\top &z_3^\top &z_4^\top \end{bmatrix}^\top: \quad
  z_3=Z_3b, \quad  z_4=Z_4b, \\
  &\quad b=V_{Z_1}\begin{bmatrix}\frac{1}{\nu_1} &\eta_2&\cdots &\eta_{m+j-1}\end{bmatrix}^\top,
   \quad  \eta_2,\ \ldots, \ \eta_{m+j-1}\in\mathbb{C} \},
  \end{align*}
  then with proper scaling, one can show that $\mathcal{N}_1(\widetilde{M}_{2,j})$
  is the unique  subset of $\mathcal{N}(\widetilde{M}_{2,j})$, which contains $\begin{row}z_1^\top&z_3^\top&z_4^\top\end{row}^\top$ with $0\neq z_1\in \mathbb{C}^n$, $z_3, z_4\in \mathbb{C}^j$.
  Furthermore,  it follows from Theorem~\ref{theorem-real} in Section~\ref{subsection3.6} that
  $\Re(\psi_1)$ and $\Im(\psi_1)$ are linearly independent.

  We then select $\tilde{p}_{j+1}$ and $\tilde{p}_{j+2}$ as the vectors generated by
  the Jacobi orthogonal transformation on  $\Re(\psi_1)$ and $\Im(\psi_1)$:
  \begin{align}\label{jacobi}
  \begin{row}\tilde{p}_{j+1}&\tilde{p}_{j+2}\end{row}=
  \begin{row}\Re(\psi_1)&\Im(\psi_1)\end{row}\begin{bsmallmatrix}c&s\\-s&c\end{bsmallmatrix},
  \end{align}
 with $c$ and $s$ selected to enforce $\tilde{p}_{j+1}^\top\tilde{p}_{j+2}=0$.
  Notice that $D_{\delta_{j+1}}(\sigma_{j+1}+i\tau_{j+1})$  is already determined
  with $\delta_{j+1}=\|\tilde{p}_{j+1}\|_2/\|\tilde{p}_{j+2}\|_2$
  (with $\varsigma_1=\|\tilde{p}_{j+1}\|_2$, $\varsigma_2=\|\tilde{p}_{j+2}\|_2$).
  Accordingly, $\tilde{v}_{j+1,S}$, $\tilde{v}_{j+2,S}$, $\tilde{v}_{j+1,T}$ and  $\tilde{v}_{j+2,T}$ will be  selected from
  \begin{equation} \label{v1v2}
  \begin{aligned}
  \begin{row}\tilde{v}_{j+1,S}&\tilde{v}_{j+2,S}\end{row}=\begin{row}\Re(Z_3b)&\Im(Z_3b)\end{row}
  \begin{bsmallmatrix}c&s\\-s&c\end{bsmallmatrix}, \\
  \begin{row}\tilde{v}_{j+1,T}&\tilde{v}_{j+2,T}\end{row}=\begin{row}\Re(Z_4b)&\Im(Z_4b)\end{row}
  \begin{bsmallmatrix}c&s\\-s&c\end{bsmallmatrix},
  \end{aligned}
  \end{equation}
  where $b=V_{Z_1}\begin{bmatrix}\frac{1}{\nu_1} &\eta_2&\cdots &\eta_{m+j-1}\end{bmatrix}^\top$,
  with $\eta_2, \ldots, \eta_{m+j-1} \in\mathbb{C}$ to be determined.
 Our goal will then be  to choose  some appropriate $\eta$'s to minimize
 the first term in \eqref{eqcom-opt-obj}.

  Define $\begin{row}w&W\end{row} \equiv
  \begin{row}Z_3^\top &Z_4^\top\end{row}^\top V_{Z_1}$  with $w\in\mathbb{C}^{2j}$, $W\in\mathbb{C}^{2j\times (m+j-1)}$,
  $K_1\equiv\begin{row}\Re(W)&-\Im(W)\end{row}$,
  $K_2\equiv\begin{row}\Im(W)&\Re(W)\end{row}$,
  and
  \[ g=\Re(g)+i\Im(g)\equiv\begin{row}\eta_2&\cdots&\eta_{m+j-1}\end{row}^\top,\]
  then a simple manipulation shows that the first term in \eqref{eqcom-opt-obj} equals
  \begin{equation}\label{opt}
  \begin{aligned}
   &\sum_{l=1}^2 \left[ \frac{\|\tilde{v}_{j+l,S}\|_2^2}{\varsigma_l^2} +
   \frac{\|\tilde{v}_{j+l,T}\|_2^2}{\varsigma_l^2} \right] \\
  =& \left\{ \begin{row}\Re(g)^\top &\Im(g)^\top\end{row}H
      + h^\top \right\}
     \begin{row}\Re(g)^\top & \Im(g)^\top\end{row}^\top + \zeta,
  \end{aligned}
  \end{equation}
  where
  \begin{multline*}
  H=\tfrac{1}{\|\tilde{p}_{j+1}\|_2^2}(cK_1-sK_2)^\top(cK_1-sK_2)\\
      +\tfrac{1}{\|\tilde{p}_{j+2}\|_2^2}(sK_1+cK_2)^\top(sK_1+cK_2),
  \end{multline*}
  \begin{multline*}
  h=\tfrac{2}{\nu_1}\left(\tfrac{c^2}{\|\tilde{p}_{j+1}\|_2^2}+\tfrac{s^2}{\|\tilde{p}_{j+2}\|_2^2}\right)K_1^\top\Re(w)\\
  +\tfrac{2}{\nu_1}\left(\tfrac{s^2}{\|\tilde{p}_{j+1}\|_2^2}+\tfrac{c^2}{\|\tilde{p}_{j+2}\|_2^2}\right)K_2^\top\Im(w)\\
   +\tfrac{2cs}{\nu_1}\left(\tfrac{1}{\|\tilde{p}_{j+2}\|_2^2}
  -\tfrac{1}{\|\tilde{p}_{j+1}\|_2^2}\right)(K_2^\top\Re(w)+ K_1^\top \Im(w)),
  \end{multline*}
  \begin{multline*}
  \zeta= \left(\tfrac{c^2}{\|\tilde{p}_{j+1}\|_2^2}+\tfrac{s^2}{\|\tilde{p}_{j+2}\|_2^2}\right)\tfrac{\|\Re(w)\|_2^2}{\nu_1^2}\\
  +\left(\tfrac{s^2}{\|\tilde{p}_{j+1}\|_2^2}+\tfrac{c^2}{\|\tilde{p}_{j+2}\|_2^2}\right)\tfrac{\|\Im(w)\|_2^2}{\nu_1^2}\\
    +\tfrac{ 2cs}{\nu_1^2}\left(\tfrac{1}{\|\tilde{p}_{j+2}\|_2^2}
  -\tfrac{1}{\|\tilde{p}_{j+1}\|_2^2}\right)\Re(w)^\top \Im(w).
  \end{multline*}
  Obviously, $H$ is symmetric semipositive definite. In fact, $H$ is symmetric positive definite. For if  $Hf=0$ with $f\in\mathbb{R}^{2m+2j-2}$,
  then $K_1f=K_2f=0$ by the definition of $H$. On the other hand, it follows from
  the definitions of $K_1,K_2$ and $W$ that $K_1^\top K_1+K_2^\top K_2=I_{2(m+j-1)}$. Hence $f=0$,
  proving  that $H$ is nonsingular.
  Consequently, the minimizer of  \eqref{opt} is given by
 \begin{align}\label{getg}
 \begin{row}\Re(g)^\top &\Im(g)^\top\end{row}^\top=-\frac{1}{2}H^{-1}h.
 \end{align}

  Once we obtain $g\in\mathbb{C}^{m+j-1}$,  $\tilde{v}_{j+1,S}$,  $\tilde{v}_{j+2,S}$,  $\tilde{v}_{j+1,T}$ and
  $\tilde{v}_{j+2,T}$ can be computed via \eqref{v1v2}.
  Also, we observe that
  $b=(c+is)V_{Z_1}\begin{bmatrix}\frac{1}{\nu_1}&g^\top\end{bmatrix}^\top$ here, with
  $c$ and $s$ from the Jacobi orthogonal transformation on
  $[ \Re(\psi_1), \Im(\psi_1) ]$.

  We still need to determine $\tilde{x}_{j+1}$
  and $\tilde{x}_{j+2}$, where
  $\tilde{x}_{j+1}+i\tilde{x}_{j+2}=Q_1y+Q_2Q_2^\top A(\tilde{p}_{j+1}+i\tilde{p}_{j+2}) -
  Q_2(Q_2^\top X_j)(\tilde{v}_{j+1,S}+i\tilde{v}_{j+2,S})$ for some $y\in\mathbb{C}^m$, which is equivalent to
  \begin{align*}
  &x_{j+1}=\tfrac{1}{\|\tilde{p}_{j+1}\|_2}\left[ Q_1\Re(y)+Q_2Q_2^\top A\tilde{p}_{j+1}-Q_2(Q_2^\top X_j)\tilde{v}_{j+1,S} \right], \\
  &x_{j+2}=\tfrac{1}{\|\tilde{p}_{j+2}\|_2}\left[ Q_1\Im(y)+Q_2Q_2^\top A\tilde{p}_{j+2}-Q_2(Q_2^\top X_j)\tilde{v}_{j+2,S} \right].
  \end{align*}
  This implies that  $Q_2^\top x_{j+l}=Q_2^\top A p_{j+l}-(Q_2^\top X_j) v_{j+l,S}$ ($l=1, 2$).
  Again, as pointed out previously, only
  $Q_2^\top X_{j+2}=Q_2^\top \begin{row}X_{j}&x_{j+1}&x_{j+2}\end{row}$ is required
  for the assigning procedure to continue. When computing $x_{j+1}$ and $x_{j+2}$, we may rewrite
  $y_{j+1}=\Re(y)/\|\tilde{p}_{j+1}\|_2$ and  $y_{j+2}=\Im(y)/\|\tilde{p}_{j+2}\|_2$, which will be selected in Section~\ref{subsection3.5}.

  \item ($\rank(Z_1)\geq2$) In this case, we shall
  employ the strategy for placing complex conjugate  pairs  in \cite{GCQX}.
  It produces reasonably good suboptimal feasible points for \eqref{eqcom-opt}. We shall sketch the placement process;
  for details, please consult \cite{GCQX}.

  We set
  \[
  b=V_{Z_1}\begin{bmatrix}\frac{e_1}{\nu_1}&\frac{e_2}{\nu_2}\end{bmatrix}
  \begin{bsmallmatrix}\gamma_1+i\zeta_1\\\gamma_2+i\zeta_2 \end{bsmallmatrix},
  \]
  where $\gamma_1, \gamma_2, \zeta_1, \zeta_2\in\mathbb{R}$ are to be determined with
  $\gamma_1^2+\gamma_2^2+\zeta_1^2+\zeta_2^2=1$, and let  $\psi_1$ and $\psi_2$ denote
  the left singular vectors of $Z_1$ corresponding to its two largest singular values $\nu_1$
  and $\nu_2$, respectively. It then follows that
  \begin{equation}\label{pv12}
  \begin{aligned}
  &\tilde{p}_{j+1}+i\tilde{p}_{j+2}=\begin{bmatrix}\psi_1&\psi_2\end{bmatrix}
  \begin{bmatrix}\gamma_1+i\zeta_1\\\gamma_2+i\zeta_2 \end{bmatrix},\\
  &\begin{bmatrix}\tilde{v}_{j+1,S}+i\tilde{v}_{j+2,S}\\ \tilde{v}_{j+1,T}+i\tilde{v}_{j+2,T}\end{bmatrix}
  =\begin{bmatrix}w_1 &w_2\end{bmatrix}\begin{bmatrix}\gamma_1+i\zeta_1\\\gamma_2+i\zeta_2 \end{bmatrix},
  \end{aligned}
  \end{equation}
  where $w_l=\frac{1}{\nu_l}\begin{bmatrix}V^\top_{Z_1}Z_3^\top &V^\top_{Z_1}Z_4^\top\end{bmatrix}^\top e_l$
  for $l=1, 2$. In the case of $\Re(\psi_1)^\top\Im(\psi_1)=0$
  and $\|\Re(\psi_1)\|_2=\|\Im(\psi_1)\|_2=\frac{1}{\sqrt{2}}$,
  we simply take $\gamma_1=1$, $\zeta_1=\gamma_2=\zeta_2=0$,
  yielding  $\tilde{p}_{j+1}=\Re(\psi_1)$ and $\tilde{p}_{j+2}=\Im(\psi_1)$.
  This actually gives the objective function in \eqref{eqcom-opt-obj}
  its minimum $2(1-\nu_1^2)/\nu_1^2$.
  In general, there are two simple possibilities.
  One is to  apply  the Jacobi orthogonal transformation on
  $\left[ \Re(\psi_1), \Im(\psi_1) \right]$ to produce  $\left[ \tilde{p}_{j+1}, \tilde{p}_{j+2} \right]$. This postulates that
  $\Re(\psi_1)$ and $\Im(\psi_1)$ are linearly independent,
  and  the value of the objective function in \eqref{eqcom-opt-obj} equals
  \begin{align*}
  \varrho_1=&\frac{\|c\Re(w_1)-s\Im(w_1)\|_2^2}{\|\tilde{p}_{j+1}\|_2^2}+
  \frac{\|s\Re(w_1)+c\Im(w_1)\|_2^2}{\|\tilde{p}_{j+2}\|_2^2}\\
  &+\tau_{j+1}^2 \left( \frac{\|\tilde{p}_{j+1}\|_2}{\|\tilde{p}_{j+2}\|_2}-\frac{\|\tilde{p}_{j+2}\|_2}{\|\tilde{p}_{j+1}\|_2} \right)^2 \\
  \leq & \frac{1}{\min\{\|\tilde{p}_{j+1}\|^2_2, \ \|\tilde{p}_{j+2}\|_2^2\}}  \left( \frac{1-\nu_1^2}{\nu_1^2}+\tau_{j+1}^2 \right).
  \end{align*}
  The other possibility makes use of the following  spectral decomposition of the Hamiltonian matrix \cite{GCQX}:
  \begin{multline*}
  \begin{bmatrix}K_R^{\top}K_R-K_I^{\top}K_I&-(K_R^{\top}K_I+K_I^{\top}K_R)\\
  -(K_R^{\top}K_I+K_I^{\top}K_R)&K_I^{\top}K_I-K_R^{\top}K_R \end{bmatrix}
  \\=\Omega 
  \diag(\phi_1, \phi_2,  -\phi_1, -\phi_2) \Omega^\top,
  \end{multline*}
  with  $K_R=\begin{row}\Re(\psi_1)&\Re(\psi_2)\end{row}$,
  $K_I=\begin{row}\Im(\psi_1)&\Im(\psi_2)\end{row}$, $\phi_1\geq \phi_2>0$.
  Some  $\gamma_1, \gamma_2, \zeta_1, \zeta_2$ are chosen (essentially determined by  $\phi_1, \phi_2$) such that
  $\tilde{p}_{j+1}^\top \tilde{p}_{j+2}=0$ and $\|\tilde{p}_{j+1}\|_2=\|\tilde{p}_{j+2}\|_2=\frac{1}{\sqrt{2}}$.
  This eventually gives the  objective function in \eqref{eqcom-opt-obj} the value
  \[
  \varrho_2= 2\sum_{l=1}^2 \frac{1-\nu_l^2}{\nu_l^2}(\gamma_l^2+\zeta_l^2)
  \leq \frac{2(1-\nu_2^2)}{\nu_2^2}.
  \]
  Then we take the possibility corresponding to the minimum of  $\varrho_1$ and $\varrho_2$,
  choosing the $(j+1)$-th and $(j+2)$-th columns of $P_{G,F}$, $S$ and $T$ accordingly.

 \smallskip

   As in {\bf Case \expandafter{\romannumeral1}},  we also need to determine, for $l=1,2$:
\begin{align*}
x_{j+l} &= Q_1y_{j+l}+Q_2Q_2^\top Ap_{j+l}-Q_2(Q_2^\top X_j)v_{j+l,S}, \\
Q_2^\top x_{j+l} &= Q_2^\top Ap_{j+l}-(Q_2^\top X_j)v_{j+l,S},
\end{align*}
for some $y_{j+l}\in\mathbb{R}^m$ to be determined in Section~\ref{subsection3.5}.
\end{enumerate}

\subsubsection{Situation \uppercase\expandafter{\romannumeral2}} ($|\alpha_{j+1}|<|\beta_{j+1}|$)
Contrasting {\em Situation \uppercase\expandafter{\romannumeral1}}, here we have
$S(j+1\step j+2,j+1\step j+2)=D_{\delta_{j+1}}(\tilde{\sigma}_{j+1}+i\tilde{\tau}_{j+1})$ and $T(j+1\step j+2,j+1\step j+2)=I_2$.
Similarly to  \eqref{pj12}, the $(j+1)$-th and  $(j+2)$-th columns of
$P_{G,F}$, $X_{G,F}$, $S$ and $T$ satisfy
\begin{align*}
Q_2^\top A\begin{row}p_{j+1}&p_{j+2}\end{row}-Q_2^\top X_j\begin{row}v_{j+1,S}&v_{j+2,S}\end{row}
\qquad\qquad\\
-Q_2^\top\begin{row}x_{j+1}&x_{j+2}\end{row}D_{\delta_{j+1}}(\tilde{\sigma}_{j+1}+i\tilde{\tau}_{j+1})=0,\\
Q_2^\top E\begin{row}p_{j+1}&p_{j+2}\end{row}-Q_2^\top X_j\begin{row}v_{j+1,T}&v_{j+2,T}\end{row}
\qquad\qquad\\-Q_2^\top\begin{row}x_{j+1}&x_{j+2}\end{row}=0.
\end{align*}
Similar to {\em Situation \uppercase\expandafter{\romannumeral1}},
after defining $\delta_{j+1}= \varsigma_1/\varsigma_2$ with
$\varsigma_1, \varsigma_2\in\mathbb{R}$ and $\varsigma_2\neq0$, we need
to solve the constrained  optimization sub-problem:
\begin{equation}\label{eqcom-opt-b}
\begin{aligned}
\min_{\begin{subarray}{c}\varsigma_1, \varsigma_2,\\ \tilde{v}_{j+1,S},
\tilde{v}_{j+2,S},\\ \tilde{v}_{j+1,T}, \tilde{v}_{j+2,T}\end{subarray}} &
\sum_{l=1}^2 \left\{
\tfrac{\|\tilde{v}_{j+l,S}\|_2^2}{\varsigma_l^2} + 
\tfrac{\|\tilde{v}_{j+l,T}\|_2^2}{\varsigma_l^2} 
\right\}
+\tilde{\tau}_{j+1}^2 \left(\tfrac{\varsigma_1}{\varsigma_2} -\tfrac{\varsigma_2}{\varsigma_1} \right)^2\\
\mbox{s.t.}\quad & [\ \tilde{p}_{j+1}^\top+i\tilde{p}_{j+2}^\top, \ \
\tilde{x}_{j+1}^\top+i\tilde{x}_{j+2}^\top, \\
& \tilde{v}_{j+1,S}^\top+i\tilde{v}_{j+2,S}^\top, \ \  \tilde{v}_{j+1,T}^\top+i\tilde{v}_{j+2,T}^\top \ ]^{\top} \in \mathcal{N}(M_{j}), \\
&\begin{row}\tilde{p}_{j+1}&\tilde{p}_{j+2}\end{row}^\top \begin{row}\tilde{p}_{j+1}&\tilde{p}_{j+2}\end{row}=\diag(\varsigma_1^2, \varsigma_2^2),
\end{aligned}
\end{equation}
where
\[
 M_{j}=\begin{bmatrix}Q_2^\top A & -(\tilde{\sigma}_{j+1}+i\tilde{\tau}_{j+1})Q_2^\top& -Q_2^\top X_j&0\\
 Q_2^\top E&-Q_2^\top &0& -Q_2^\top X_j\\
 P_j^\top&0&0&0\end{bmatrix}.
\]
The above  optimization problem \eqref{eqcom-opt-b} can be treated similarly as the optimization problem \eqref{eqcom-opt},
We skip the details here.

\begin{note}
Analogously to last section when assigning the finite real poles,
we need to pay some attention when  $j=0$.
Suppose that no infinite poles exist  and the
first finite poles to be assigned are $(\alpha_1, \beta_1)$ and $(\bar{\alpha}_1, \bar{\beta}_1)$
with $\Im(\alpha_1)\Im(\beta_1)\neq0$.
It follows from the structure of $S$ and $T$ that  we just need to compute
the first two columns of $P_{G,F}$, $Q_2^\top X_{G,F}$ and
the $2\times 2$ leading principal submatrices of $S$ and $T$.

When $|\alpha_1| \geq |\beta_1|$ (and neglecting the  complementary  case),
the first two columns $p_1, p_2\in\mathbb{R}^n$ of $P_{G,F}$ and $x_1, x_2\in\mathbb{R}^{n}$ of $X_{G,F}$
should be chosen to satisfy
\begin{align}\label{p1p2condition}
	\left\{\begin{aligned}
Q_2^{\top}A\begin{row}p_1 &p_2\end{row} - Q_2^\top\begin{row}x_1 &x_2\end{row}&=0,\\
Q_2^{\top}E\begin{row}p_1 &p_2\end{row} - Q_2^\top\begin{row}x_1 &x_2\end{row}D_{\delta_1}(\sigma_1+i\tau_1)&=0,\\
\begin{row}p_1 &p_2 \end{row}^\top \begin{row}p_1 &p_2 \end{row}&=I_2,
\end{aligned}\right.
\end{align}
so that $(\delta_1-\delta_1^{-1})^2\tau_1^2$ is minimized. This is obviously achieved when $\delta_1=1$,
where \eqref{p1p2condition} is reduced to \eqref{eqcom-opt-constraina} and \eqref{eqcom-opt-constrainb}
with $j=0$, $\tilde{v}_{j+l, S}, \tilde{v}_{j+l, T}$ vanished,
 $\tilde{p}_{j+l}, \tilde{x}_{j+l}$ replaced by $p_{l}, x_{l}$ respectively,
 and $\varsigma_l=1$ ($l=1,2$).
Let  the columns of $Z\in\mathbb{C}^{n\times m}$  be an  orthonormal basis of
$\mathcal{N}(M)$, where $M \equiv Q_2^\top[E-(\sigma_1+i\tau_1)A]$ is of
full row rank follows by Lemma~\ref{Lemma3.2} in Section~\ref{subsection3.6}.
Then the columns of
\[
\begin{bmatrix}0& Z\\ Q_1 & Q_2Q_2^\top AZ\end{bmatrix}
\]
construct a basis of $\mathcal{N}(M_0)$. Furthermore,
 there exist $p\in\mathbb{C}^n$ and $x\in\mathbb{C}^{n}$ with $\{ \Re(p), \Im(p) \}$
being linearly  independent such that $\begin{row}p^\top &x^\top\end{row}^\top\in\mathcal{N}(M_0)$.
Obviously,
$p_1=\begin{row}\Re(Z)&-\Im(Z)\end{row}\begin{row}u_1^\top &u_2^\top\end{row}^\top$,
$p_2=\begin{row}\Im(Z)&\Re(Z)\end{row}\begin{row}u_1^\top &u_2^\top\end{row}^\top$
for some $u_1, u_2\in\mathbb{R}^{m}$. Adopting the method in \cite{GCQX},
where two Hamiltonian matrices would be constructed and their
spectral decompositions lead to
$p_1^\top p_2=0$ and $\|p_1\|_2=\|p_2\|_2=1$.
We also have
$x_l=Q_1y_{l}+Q_2Q_2^\top Ap_l$ for  some $y_l\in\mathbb{R}^m$ (from Section~\ref{subsection3.5}),
leading  to $Q_2^\top x_l=Q_2^\top Ap_l$ for $l=1, 2$. This is  sufficient for the process to continue.
\end{note}

We recapitulate the assigning process
of the  complex  conjugate poles $\{(\alpha_{j+1}, \beta_{j+1}), (\bar{\alpha}_{j+1}, \bar{\beta}_{j+1})\}$
in Algorithm \ref{algorithm3}, where  $\Xi_j \equiv Q_2^\top X_j$.
Again, the suboptimal strategies  in \cite{GCQX} for complex conjugate poles
is employed.

\begin{algorithm*}
\caption{\ Assigning complex conjugate poles $\{(\alpha_{j+1}, \beta_{j+1}), (\bar{\alpha}_{j+1}, \bar{\beta}_{j+1})\}$}
\begin{algorithmic}[1]\label{algorithm3}
\REQUIRE ~~\\
$A$, $E$, $Q_2$, $P_j$, $\Xi_j$, $S_j$, $T_j$ and  $(\alpha_{j+1}, \beta_{j+1})\in\mathbb{C}\times \mathbb{C}$.
\ENSURE ~~\\
Orthogonal $P_{j+2}\in\mathbb{R}^{n\times(j+2)}$,  upper quasi-triangular $S_{j+2}, T_{j+2}\in\mathbb{R}^{(j+2)\times(j+2)}$ and
$\Xi_{j+2}\in\mathbb{R}^{(n-m)\times(j+2)}$.
\IF {$|\alpha_{j+1}|\geq|\beta_{j+1}|$}
    \STATE Compute $\sigma_{j+1}$ and $\tau_{j+1}$ as defined in Lemma \ref{Theorem2.1}.
    \STATE Find $\begin{bmatrix}Z_1^\top &Z_3^\top &Z_4^\top\end{bmatrix}^\top$ with $Z_1\in\mathbb{C}^{n\times (m+j)}$,
           $Z_3, Z_4\in\mathbb{C}^{j\times (m+j)}$, whose orthonormal columns span $\mathcal{N}(\widetilde{M}_{2,j})$  from \eqref{M_2j-com}.
\ELSE
    \STATE Compute $\tilde{\sigma}_{j+1}$ and $\tilde{\tau}_{j+1}$ as defined in Lemma \ref{Theorem2.1}.
    \STATE Find $\begin{bmatrix}Z_1^\top &Z_3^\top &Z_4^\top\end{bmatrix}^\top$ with $Z_1\in\mathbb{C}^{n\times (m+j)}$,
           $Z_3, Z_4\in\mathbb{C}^{j\times (m+j)}$, whose orthonormal columns span the null space of
           \[
           \widetilde{M}_{2,j}=\begin{bmatrix}Q_2^\top(A-(\tilde{\sigma}_{j+1}+i\tilde{\tau}_{j+1})E) &-\Xi_j&
           \ (\tilde{\sigma}_{j+1}+i\tilde{\tau}_{j+1})\Xi_{j}\\ P_j^\top &0&0\end{bmatrix}.
           \]
\ENDIF
\IF {$\rank(Z_1)=1$}
    \STATE Perform the Jacobi orthogonal transformation on $\psi_1$ as in \eqref{jacobi} to compute
           $\tilde{p}_{j+1}$ and $\tilde{p}_{j+2}$, where $\psi_1$ is the left singular vector of
           $Z_1$  corresponding to its unique nonzero singular value $\nu_1$.
    \STATE Compute $\tilde{v}_{j+1,S}$, $\tilde{v}_{j+2,S}$, $\tilde{v}_{j+1,T}$ and $\tilde{v}_{j+2,T}$
	by \eqref{v1v2}, with $b=V_{Z_1}\begin{bmatrix}\frac{1}{\nu_1}&g^\top\end{bmatrix}^\top$,
           the SVD $Z_1=U_{Z_1}\Sigma_{Z_1}V_{Z_1}$ and  $g$ from \eqref{getg}.
\ELSE
    \STATE Solve the optimization problem  \eqref{eqcom-opt}  with the constraint \eqref{eqcom-opt-constraina}
    replaced by
    $$\widetilde{M}_{2,j}\begin{bmatrix}\tilde{p}_{j+1}^\top+i\tilde{p}_{j+2}^\top &
    \tilde{v}_{j+1,S}^\top+i\tilde{v}_{j+2,S}^\top &\tilde{v}_{j+1,T}^\top+i\tilde{v}_{j+2,T}^\top\end{bmatrix}^\top=0,$$
    employing  the suboptimal strategies in \cite{GCQX} to compute $\tilde{p}_{j+1}$, $\tilde{p}_{j+2}$, $\tilde{v}_{j+1,S}$, $\tilde{v}_{j+2,S}$, $\tilde{v}_{j+1,T}$ and $\tilde{v}_{j+2,T}$.
\ENDIF
\STATE Set $\delta_{j+1}=\frac{\|\tilde{p}_{j+1}\|_2}{\|\tilde{p}_{j+2}\|_2}$.
\STATE Set $P_{j+2}=\begin{bmatrix}P_j&p_{j+1}&p_{j+2}\end{bmatrix}$ with $p_{j+1}$ and
           $p_{j+1}$ being the normalized vectors of $\tilde{p}_{j+1}$ and $\tilde{p}_{j+2}$, respectively.
\STATE Set $v_{j+1,S}$, $v_{j+2,S}$, $v_{j+1,T}$, $v_{j+2,T}$ be the normalized vectors of $\tilde{v}_{j+1,S}$,
    $\tilde{v}_{j+2,S}$, $\tilde{v}_{j+1,T}$, $\tilde{v}_{j+2,T}$, respectively.
\IF{$|\alpha_{j+1}|\geq|\beta_{j+1}|$}
    \STATE Set
    \begin{equation*}
    S_{j+2}=
	\begin{bgapmatrix}
		S_j  &\vrule& \begin{matrix} v_{j+1,S}&v_{j+2,S}\end{matrix}\\
		\hline &\vrule& I_2
	\end{bgapmatrix} ,\quad
    T_{j+2}=
	\begin{bgapmatrix}
		T_j &\vrule& 	\begin{matrix}v_{j+1,T}&v_{j+2,T}\end{matrix}\\
	\hline &\vrule& D_{\delta_{j+1}(\sigma_{j+1}+i\tau_{j+1})} \\
	\end{bgapmatrix}.
    \end{equation*}
    \STATE  Set $\Xi_{j+2}=\begin{bmatrix}\Xi_j&\xi_{j+1}&\xi_{j+2}\end{bmatrix}$ with
    $\xi_{j+1}=Q_2^\top Ap_{j+1}-\Xi_jv_{j+1,S}$ and  $\xi_{j+2}=Q_2^\top Ap_{j+2}-\Xi_jv_{j+2,S}$.
\ELSE
    \STATE Set
    \begin{equation*}
    S_{j+2}=
	\begin{bgapmatrix}
		S_j  &\vrule& \begin{matrix} v_{j+1,S}&v_{j+2,S}\end{matrix}\\
    \hline &\vrule& D_{\delta_{j+1}(\tilde{\sigma}_{j+1}+i\tilde{\tau}_{j+1})}
	\end{bgapmatrix}, \quad
	T_{j+2}=
	\begin{bgapmatrix}
		T_j &\vrule& 	\begin{matrix}v_{j+1,T}&v_{j+2,T}\end{matrix}\\
	\hline &\vrule& I_2
	\end{bgapmatrix} .
    \end{equation*}
	\STATE  Set $\Xi_{j+2}=\begin{bmatrix}\Xi_j&\xi_{j+1}&\xi_{j+2}\end{bmatrix}$ with
    $\xi_{j+1}=Q_2^\top Ep_{j+1}-\Xi_jv_{j+1,T}$ and  $\xi_{j+2}=Q_2^\top Ep_{j+2}-\Xi_jv_{j+2,T}$.
\ENDIF
\end{algorithmic}
\end{algorithm*}

\subsection{Determining $X_{G,F}$}\label{subsection3.5}

We have $X_{G,F}=Q_1Y+Q_2 (Q_2^\top X_{G,F})$,
where $Q_2^\top X_{G,F}$ has been computed and $Y=\begin{row}y_1&\cdots&y_n\end{row}
\in\mathbb{R}^{m\times n}$ is to be determined. This last gap is to be filled in this section.

\begin{lem}
The computed matrix $Q_2^\top X_{G,F}$ is of full row rank.
\end{lem}
\begin{pf}
Since the descriptor system $(E, A, B)$  is S-controllable, then
$\begin{row}E&AN_{\infty}&B\end{row}$ is of full row rank
with $\mathcal{R}(N_{\infty})$ being the null space of $E$.
This consequently leads to $\rank(\begin{row}Q_2^\top E&Q_2^\top A N_{\infty}\end{row})=n-m$.
It hence holds that
$\begin{row}Q_2^\top E &Q_2^\top A\end{row}$ is of full row rank, which is equivalent to
$\begin{row}Q_2^\top EP_{G,F} &Q_2^\top A P_{G,F}\end{row}$ of having full row rank.
Also, it follows from
$Q_2^\top AP_{G,F}=Q_2^\top X_{G,F}S$ and $Q_2^\top EP_{G,F}=Q_2^\top X_{G,F}T$ that
$\begin{row}Q_2^\top X_{G,F}T&Q_2^\top X_{G,F}S\end{row}$ is of full row rank,
yielding  the same for $Q_2^\top X_{G,F}$.
\end{pf}

Now rewrite
\[
X_{G,F}=\begin{bmatrix}Q_1&Q_2\end{bmatrix}\begin{bmatrix}Y\\ Q_2^\top X_{G,F}\end{bmatrix},
\]
it is nonsingular with $Y^\top$ not deficient in
the complementary subspace  of $\mathcal{R}(X_{G,F}^\top Q_2)$.
Furthermore, we wish $X_{G,F}$ to be as well conditioned as possible. Thus we should choose  $Y^\top$
whose orthonormal columns span the complementary subspace  of $\mathcal{R}(X_{G,F}^\top Q_2)$.
From the QR factorization  $X_{G,F}^\top Q_2 =Q_X\begin{row}R_X^\top &0\end{row}^\top
=Q_{1,X}R_X$
with  $Q_X = \begin{row}Q_{1,X} &Q_{2,X}\end{row} \in\mathbb{R}^{n\times n}$ being orthogonal, $Q_{1,X}\in\mathbb{R}^{n\times(n-m)}$  and
$R_X\in\mathbb{R}^{(n-m)\times(n-m)}$ being nonsingular upper triangular,
we then select $Y=Q_{2,X}^\top$, leading to
$X_{G,F}=Q_1Q_{2,X}^\top+Q_2(Q_2^\top X_{G,F})$.

\subsection{Supporting Theorems }\label{subsection3.6}

\begin{lem}\label{Lemma3.2}
For any $\lambda\in\mathbb{C}$, $Q_2^\top(\lambda E-A)$ is of full row rank.
\end{lem}
\begin{pf}
Since $\begin{bmatrix}\lambda E-A &B\end{bmatrix}$ is of full row rank for all $\lambda \in\mathbb{C}$ and
\[
\begin{bmatrix}Q_1 &Q_2\end{bmatrix}^\top\begin{bmatrix}\lambda E-A &B\end{bmatrix}=
\begin{bmatrix}Q_1^\top(\lambda E-A) & R\\Q_2^\top(\lambda E-A)&0 \end{bmatrix},
\]
the result holds because  $R$ is nonsingular.
\end{pf}

\begin{thm}\label{theorem-infinite}
For an S-controllable descriptor system $(E, A, B)$, assume
$j$ infinite poles  ($0\leq j \leq (n-r-1)$) have already been assigned with
$P_{j}=\begin{row}p_1&\cdots&p_{j}\end{row}=ZW_{j}$  and
$Q_2^\top X_{j}=Q_2^\top \begin{row}x_1&\cdots &x_{j}\end{row}$ computed,
where the orthonormal columns of  $Z\in\mathbb{R}^{n\times l}$ span $Q_2^\top E$, $l=n-\rank(Q_2^\top E)$ and $W_{j}\in\mathbb{R}^{l\times j}$ satisfies
$W_{j}^\top W_{j}=I_{j}$. Define
\[
M=\begin{bmatrix} Q_2^\top AZ & -Q_2^\top& -Q_2^\top X_{j} \\  W_{j}^\top &0 &0\end{bmatrix},
\]
then we have
\begin{enumerate}
 \item[(a)] $\dim(\mathcal{N}(M))=m+l$; and
 \item[(b)] there exist a nonzero $w\in\mathbb{R}^l$, $x\in\mathbb{R}^{n}$ and $v_S\in\mathbb{R}^{j}$ such that
            $\begin{row}w^\top &x^\top  &v_S^\top\end{row}^\top \in\mathcal{N}(M)$.
\end{enumerate}
\end{thm}

\begin{pf}
We firstly consider {\em (a)}, which is equivalent to $M$ possessing full row rank.
Let $f\in\mathbb{R}^{n-m}$, $h\in\mathbb{R}^{j}$ be vectors satisfying
$\begin{row}f^\top &h^\top\end{row} M=0$,
we have $f^\top Q_2^\top AZ + h^\top W_{j}^\top=0$ and $f^\top Q_2^\top=0$. Hence
$f$ and $h$ vanish for $Q_2$ and $W_{j}$ are of full column rank, implying the result.

For {\em (b)}, assume the contrary and we have $\rank(Q_2^\top)=\rank(\begin{row}Q_2^\top &Q_2^\top X_{j}\end{row})= (n-m)+(j-l)$,
implying   $j=l$.
Since $l=n-\rank(Q_2^\top E)$ and $\rank(\begin{row}E&B\end{row})=m+\rank(Q_2^\top E)$,
then $l=j=n-(\rank(\begin{row}E&B\end{row})-m)\geq n-r$ since
$\rank(\begin{row}E&B\end{row})-m\leq r$.
On the other hand, $j\leq n-r-1$ since there exists at least one infinite pole that  is not placed.
Thus we get a contradiction and {\em (b)} holds.
\end{pf}

\begin{thm}\label{theorem-real}
For an S-controllable descriptor system $(E, A, B)$,
assume all infinite poles
\[\{(\alpha_1, \beta_1), \ldots, (\alpha_{n-r}, \beta_{n-r})\}\]
and  $j$  finite poles
\[\{(\alpha_{n-r+1}, \beta_{n-r+1}), \ldots, (\alpha_{n-r+j}, \beta_{n-r+j})\}\subseteq\mathfrak{L}\]
have already been assigned as described in Section~\ref{section3},
where  $j<r$ if there is at least one  unassigned  finite real pole,
or  $(j+1)<r$ when there is  at least a pair of  unassigned complex conjugate   poles.
Let $P_{n-r+j}=\begin{row}p_1&\cdots&p_{n-r+j}\end{row}$ contain the first $n-r+j$ columns
of $P_{G,F}$, $S_{n-r+j}$ and $T_{n-r+j}$ be the
$(n-r+j)\times (n-r+j)$ principal submatrices of $S$ and $T$,  respectively, and
$Q_2^\top X_{n-r+j}=Q_2^\top \begin{row}x_1&\cdots&x_{n-r+j}\end{row}$, all
computed in the assigning process. Thus we have
\begin{align}\label{Q_2AE}
\left\{\begin{aligned}
Q_2^\top AP_{n-r+j}&=Q_2^\top X_{n-r+j}S_{n-r+j},\\
Q_2^\top EP_{n-r+j}&=Q_2^\top X_{n-r+j}T_{n-r+j}.
\end{aligned}
\right.
\end{align}
Assume that $(\alpha, \beta)\in\mathfrak{L}$ is  the  finite real pole or
$\{(\alpha, \beta), (\bar{\alpha}, \bar{\beta})\}\subseteq\mathfrak{L}$
are the complex conjugate poles to be assigned.
Denote
\begin{align*}
M=\begin{bmatrix}Q_2^\top A&-\epsilon_1Q_2^\top&-Q_2^\top X_{n-r+j}&0\\
Q_2^\top E&-\epsilon_2Q_2^\top&0&-Q_2^\top X_{n-r+j}\\
P_{n-r+j}^\top&0&0&0\end{bmatrix},
\end{align*}
where (i) $\epsilon_1=\frac{\alpha }{\sqrt{|\alpha|^2+|\beta|^2}}$,
$\epsilon_2=\frac{\beta}{\sqrt{|\alpha|^2+|\beta|^2}}$ for $(\alpha, \beta)\in\mathbb{R}\times \mathbb{R}$;
(ii) $\epsilon_1=1$ and $\epsilon_2=\sigma+i\tau$ for $(\alpha, \beta)\in\mathbb{C}\times \mathbb{C}$ and $|\alpha|\geq|\beta|$,
with $\sigma=\frac{\Re(\bar{\alpha}\beta)}{|\alpha|^2}$, $\tau=\frac{\Im(\bar{\alpha}\beta)}{|\alpha|^2}$;
or (iii) $\epsilon_1=\tilde{\sigma}+i\tilde{\tau}$, $\epsilon_2=1$
for $(\alpha, \beta)\in\mathbb{C}\times \mathbb{C}$ and $|\alpha|<|\beta|$,
with $\sigma=\frac{\Re(\bar{\beta}\alpha)}{|\beta|^2}$ and $\tau=\frac{\Im(\bar{\beta}\alpha)}{|\beta|^2}$.
Let the columns of
\begin{align*}
Z=\begin{bmatrix} Z_{1}\\Z_{2}\end{bmatrix}
\begin{array}{>{\scriptstyle}l}n\\n+2(n-r+j)\end{array}
\end{align*}
be an orthonormal basis of $\mathcal{N}(M)$,  then we have:
\begin{enumerate}
 \item[(a)] $\dim(\mathcal{R}(Z))= 2m+(n-r+j)$;
 \item[(b)] $Z_1\neq 0$; and
 \item[(c)] for $(\alpha, \beta)\in\mathbb{C}\times \mathbb{C}$, there exist $0\neq p=\Re(p)+i\Im(p)\in\mathbb{C}^n$ with $\{ \Re(p), \Im(p) \}$ being linearly independent,
 $x\in\mathbb{C}^{n}$, $v_S\in\mathbb{C}^{n-r+j}$ and $v_T\in\mathbb{C}^{n-r+j}$
 such that
 $
 \begin{row}p^\top &x^\top &v_S^\top&v_T^\top\end{row}^\top\in\mathcal{R}(Z).
 $
\end{enumerate}
\end{thm}

\begin{pf}
Obviously, $\dim(\mathcal{R}(Z))= 2m+(n-r+j)$ if and only if $M$ has full row rank.
Suppose $z\in\mathbb{C}^{n-m}$, $y\in\mathbb{C}^{n-m}$ and $w\in\mathbb{C}^{n-r+j}$
satisfy $\begin{row}z^\top&y^\top&w^\top\end{row}M=0$, which is equivalent to
\begin{subequations}\label{uvw}
\begin{align}
z^\top Q_2^\top A + y^\top Q_2^\top E + w^\top P_{n-r+j}^\top&=0,\label{uvw1}\\
\epsilon_1 z^\top Q_2^\top + \epsilon_2 y^\top Q_2^\top&=0,\label{uvw2}\\
z^\top Q_2^\top X_{n-r+j}=y^\top Q_2^\top X_{n-r+j}&=0. \label{uvw3}
\end{align}
\end{subequations}
Post-multiplying $P_{n-r+j}$ on both sides of \eqref{uvw1} gives
$z^\top Q_2^\top AP_{n-r+j} + y^\top Q_2^\top EP_{n-r+j} + w^\top=0$. Together with
\eqref{Q_2AE} and  \eqref{uvw3}, we get
\begin{align*}
z^\top Q_2^\top AP_{n-r+j}&=z^\top Q_2^\top X_{n-r+j}S_{n-r+j}=0,\\
y^\top Q_2^\top EP_{n-r+j}&=y^\top Q_2^\top X_{n-r+j}T_{n-r+j}=0,
\end{align*}
leading to $w=0$. Thus $z^\top Q_2^\top A + y^\top Q_2^\top E=
z^\top(Q_2^\top A - \frac{\epsilon_1}{\epsilon_2}Q_2^\top E)=0$ follows from \eqref{uvw1} and $\epsilon_2\neq0$ for  all  three cases.
Furthermore from Lemma \ref{Lemma3.2}, $Q_2^\top (A - \frac{\epsilon_1}{\epsilon_2}E)$ is
of full row rank, implying that $y=z=0$.
So $M$ is of full row rank, hence {\em (a)} holds. 

To prove {\em (b)}, we assume the contrary that $Z_1=0$. Then we have
\begin{multline*}
\rank\left(\begin{bmatrix}-\epsilon_1Q_2^\top&-Q_2^\top X_{n-r+j}&0\\
-\epsilon_2Q_2^\top&0&-Q_2^\top X_{n-r+j}\end{bmatrix}\right)\\
=(n-m)+(n-r+j-m).
\end{multline*}
Since $Q_2$ is of full column rank and
\begin{equation}\label{Y_j^2}
\begin{aligned}
&\begin{bmatrix}0&-Q_2^\top X_{n-r+j}&\frac{\epsilon_1}{\epsilon_2}Q_2^\top X_{n-r+j}\\
-\epsilon_2Q_2^\top&0&0\end{bmatrix}\\
=&\begin{bmatrix}I_{n-m}&-\frac{\epsilon_1}{\epsilon_2}I_{n-m}\\ &I_{n-m}\end{bmatrix}
\begin{bmatrix}-\epsilon_1Q_2^\top&-Q_2^\top X_{n-r+j}&0\\
-\epsilon_2Q_2^\top&0&-Q_2^\top X_{n-r+j}\end{bmatrix}\\
&\qquad \qquad \qquad \qquad \cdot
\begin{bmatrix}I_{n}&0&-\frac{1}{\epsilon_2} X_{n-r+j}\\ &I_{n-r+j}&0\\ &&I_{n-r+j}\end{bmatrix},
\end{aligned}
\end{equation}
we deduce that
$ \rank(\begin{row}-Q_2^\top X_{n-r+j}&\frac{\epsilon_1}{\epsilon_2}Q_2^\top X_{n-r+j}\end{row})
=\rank(Q_2^\top  X_{n-r+j})=(n-r+j)-m$,
thus requiring $(n-r+j)\geq m$.
When $(n-r+j)<m$, we have $Z_1\neq0$ and we need to consider the complementary case when $(n-r+j)\geq m$.  Assume that $Q_2^\top X_{n-r+j}H=0$ with
$H\in\mathbb{R}^{(n-r+j)\times m}$ and $H^\top H=I_m$,
then $\begin{row}H&X_{n-r+j}^\top Q_2\end{row}^\top$ is of full column rank.
Define $Y_{n-r+j}=Q\begin{row}H&X_{n-r+j}^\top Q_2\end{row}^\top$,
we have $Q_2^\top Y_{n-r+j}=Q_2^\top X_{n-r+j}$.
From \eqref{Q_2AE},  there exist
$W_A, W_E\in\mathbb{R}^{m\times (n-r+j)}$ such that
\begin{align*}
\left\{\begin{aligned}
AP_{n-r+j}&=Y_{n-r+j} S_{n-r+j}+BW_A, \\
EP_{n-r+j}&=Y_{n-r+j} T_{n-r+j}+BW_E.
\end{aligned}
\right.
\end{align*}
Moreover, it follows from
\[
	Y_{n-r+j}H=Q\begin{row}H&X_{n-r+j}^\top Q_2\end{row}^\top H
=Q_1=BR^{-1}
\]
that $B=Y_{n-r+j} HR$. Consequently, it holds that
\begin{align}\label{HRW}
	\left\{\begin{aligned}
AP_{n-r+j}&=Y_{n-r+j}(S_{n-r+j}+HRW_A), \\
EP_{n-r+j}&=Y_{n-r+j}(T_{n-r+j}+HRW_E).
\end{aligned}
\right.
\end{align}
Define  $K\in\mathbb{R}^{n\times (r-j)}$ satisfying
$K^\top Y_{n-r+j}=0$ and $K^\top K=I_{r-j}$,
then pre-multiplying
\[
	L=\begin{bmatrix}(Y_{n-r+j}^\top Y_{n-r+j})^{-1}&\\&I_{r-j}\end{bmatrix}
\begin{bmatrix}Y_{n-r+j}^\top \\ K^\top\end{bmatrix}
\] on both sides of \eqref{HRW}  yields
\begin{align*}
L AP_{n-r+j}&=\begin{bmatrix}S_{n-r+j}+HRW_A\\ 0\end{bmatrix}, \\
L EP_{n-r+j}&=\begin{bmatrix}T_{n-r+j}+HRW_E\\ 0\end{bmatrix}.
\end{align*}
Let $P_{\bot}\in\mathbb{R}^{n\times (r-j)}$
satisfying $P_{\bot}^\top P_{n-r+j}=0$ and $P_{\bot}^\top P_{\bot}=I_{r-j}$,
we have that
\begin{equation}
\begin{aligned}
L A\begin{bmatrix}P_{n-r+j}&P_{\bot}\end{bmatrix}&=\begin{bmatrix}S_{n-r+j}+HRW_A&A_{12}\\0&A_{22}\end{bmatrix}, \\
L E\begin{bmatrix}P_{n-r+j}&P_{\bot}\end{bmatrix}&=\begin{bmatrix}T_{n-r+j}+HRW_E&E_{12}\\0&E_{22}\end{bmatrix},
\end{aligned}
\end{equation}
where $\begin{row}A_{12}^\top&A_{22}^\top\end{row}^\top=LAP_{\bot}$ and
$\begin{row}E_{12}^\top&E_{22}^\top\end{row}^\top=LEP_{\bot}$.
Since $LB=\begin{row}R^\top H^\top&\ 0\end{row}^\top$, then for
the system $(LE\begin{row}P_{n-r+j}&\ P_{\bot}\end{row}, \quad
LA\begin{row}P_{n-r+j}&\ P_{\bot}\end{row}, \quad LB)$ (which is equivalent to the
descriptor system $(E, A, B)$), there are at most $j$ finite poles
assignable with the PD-SF.
Such result obviously contradicts with Lemma  \ref{Theorem2.2} since $j<r$,
hence {\em (b)} holds.

Regarding {\em (c)}, we just give the proof when $\epsilon_1=1$ and $\epsilon_2=\sigma+i\tau$.
(When $\epsilon_1=\tilde{\sigma}+i\tilde{\tau}$ and $\epsilon_2=1$, the proof is similar and ignored.)

If
\begin{multline*}
\rank\left(\begin{bmatrix}-\epsilon_1Q_2^\top&-Q_2^\top X_{n-r+j}&0\\
-\epsilon_2Q_2^\top&0&-Q_2^\top X_{n-r+j}\end{bmatrix}\right)\\
\geq (n-m)+(n-r+j-m)+2,
\end{multline*}
then there exist
\[
\begin{bmatrix}\tilde{p}_1^\top &\tilde{x}_1^\top &\tilde{v}_{1,S}^\top&\tilde{v}_{1,T}^\top\end{bmatrix}^\top, \quad
\begin{bmatrix}\tilde{p}_2^\top &\tilde{x}_2^\top  &\tilde{v}_{2,S}^\top&\tilde{v}_{2,T}^\top\end{bmatrix}^\top\in\mathcal{R}(Z),
\]
where $\tilde{p}_1, \tilde{p}_2\in\mathbb{C}^n$, $\tilde{x}_1, \tilde{x}_2\in\mathbb{C}^{n}$,
$\tilde{v}_{1,S}, \tilde{v}_{2,S}\in\mathbb{C}^{n-r+j}$ and  $\tilde{v}_{1,T}, \tilde{v}_{2,T}\in\mathbb{C}^{n-r+j}$,
such that $\tilde{p}_1$ and $\tilde{p}_2$ are linearly independent.
Let
\begin{align*}
\begin{bmatrix}p^\top &x^\top &v_{S}^\top&v_{T}^\top\end{bmatrix}^\top
&=(\xi_1+i\eta_1)\begin{bmatrix}\tilde{p}_1^\top &\tilde{x}_1^\top &\tilde{v}_{1,S}^\top &\tilde{v}_{1,T}^\top
\end{bmatrix}^\top\\
&\qquad+(\xi_2+i\eta_2)\begin{bmatrix}\tilde{p}_2^\top &\tilde{x}_2^\top  &\tilde{v}_{2,S}^\top &\tilde{v}_{2,T}^\top
\end{bmatrix}^\top,
\end{align*}
then we can always find suitable $\xi_1, \xi_2, \eta_1, \eta_2 \in\mathbb{R}$ such that the real and the imaginary parts of
the resulting $p\in\mathbb{C}^n$ are linearly independent.

While
\begin{multline*}
\rank\left(\begin{bmatrix}-\epsilon_1Q_2^\top&-Q_2^\top X_{n-r+j}&0\\
-\epsilon_2Q_2^\top&0&-Q_2^\top X_{n-r+j}\end{bmatrix}\right)\\
= (n-m)+(n-r+j-m)+1,
\end{multline*}
it follows from \eqref{Y_j^2} that
\begin{align*}
	\rank(Q_2^\top X_{n-r+j})&=\rank(\begin{row}-Q_2^\top X_{n-r+j}&\frac{\epsilon_1}{\epsilon_2}Q_2^\top X_{n-r+j}\end{row})
	\\&=(n-r+j)-(m-1).
\end{align*}
Thus it is necessary that $(n-r+j)\geq(m-1)$,
which we assume from now on. 
If {\em (c)} does not hold, then there exist vectors  $0\neq p\in\mathbb{R}^n$,
$\Re(x)+i\Im(x)=x\in\mathbb{C}^{n}$,
$\Re(v_{S})+i\Im(v_{S})=v_{S}\in\mathbb{C}^{n-r+j}$,
$\Re(v_{T})+i\Im(v_{T})=v_{T}\in\mathbb{C}^{n-r+j}$ such that
$\begin{row}p^\top &x^\top &v_{S}^\top&v_{T}^\top\end{row}^\top\in\mathcal{N}(M)$, which
is equivalent to
\begin{subequations}\label{real-z}
\begin{align}
&Q_2^\top Ap=Q_2^\top \Re(x)+Q_2^\top X_{n-r+j}\Re(v_S), \label{real-z1}\\
&Q_2^\top Ep=Q_2^\top(\sigma\Re(x)-\tau\Im(x))+Q_2^\top X_{n-r+j}\Re(v_T), \label{real-z2}\\
&Q_2^\top\Im(x)=-Q_2^\top X_{n-r+j}\Im(v_S),\label{real-z3}\\
&Q_2^\top(\sigma \Im(x)+\tau\Re(x))+ Q_2^\top X_{n-r+j}\Im(v_T)=0.\label{real-z4}
\end{align}
\end{subequations}
From \eqref{real-z2} and \eqref{real-z3}, we have
$Q_2^\top Ep=\sigma Q_2^\top\Re(x)+Q_2^\top X_{n-r+j}[\tau\Im(v_S)+\Re(v_T)]$;
and from \eqref{real-z3} and \eqref{real-z4}, we  get
$Q_2^\top\Re(x)=Q_2^\top X_{n-r+j}[\frac{\sigma}{\tau}\Im(v_S)-\frac{1}{\tau}\Im(v_T)]$.
Consquently, writing $\widehat{Q}^\top\equiv\begin{bmatrix}Q_2^\top X_{n-r+j}&Q_2^\top\Re(x)\end{bmatrix}$, we obtain
\begin{equation}\label{Y_j^2w}
\left\{\begin{aligned}
  &Q_2^\top A\begin{bmatrix}P_{n-r+j}&p\end{bmatrix}=
  \widehat{Q}^\top\begin{bgapmatrix}S_{n-r+j}  &\vrule&\Re(v_S)\\  \hline
    &\vrule& 1\end{bgapmatrix}, \\
  &Q_2^\top E\begin{bmatrix}P_{n-r+j}&p\end{bmatrix}=
  \widehat{Q}^\top\begin{bgapmatrix}T_{n-r+j}  &\vrule&\tau\Im(v_S)+\Re(v_T)\\ \hline
	  &\vrule&\sigma\end{bgapmatrix},
  \end{aligned}
  \right.
\end{equation}
and $\rank(\widehat{Q}^\top)=\rank(Q_2^\top X_{n-r+j})=(n-r+j)-m+1$.
Now let  $H\in\mathbb{R}^{(n-r+j+1)\times m}$ satisfy
$\widehat{Q}^\top H=0$ and $H^\top H=I_{m}$.
Define
\[
	Y_{n-r+j+1}=Q\begin{bmatrix}H^\top \\ \widehat{Q}^\top \end{bmatrix},
\]
which is of full column rank,  and let
$K\in\mathbb{R}^{n\times (r-j-1)}$ be the matrix satisfying $K^\top Y_{n-r+j+1}=0$ and $K^\top K=I_{r-j-1}$.
Then it follows from \eqref{Y_j^2w} and
$Q_2^\top Y_{n-r+j+1}=\widehat{Q}^\top$
that there exist $L_A\in\mathbb{R}^{m\times (n-r+j+1)}$,
$L_E\in\mathbb{R}^{m\times (n-r+j+1)}$ such that
\begin{align*}
A\begin{row}P_{n-r+j}&p\end{row}&=Y_{n-r+j+1}\begin{bgapmatrix}S_{n-r+j}&\vrule&\Re(v_S)\\ \hline
 &\vrule&1\end{bgapmatrix}+BL_A,\\
E\begin{row}P_{n-r+j}&p\end{row}&=Y_{n-r+j+1}\begin{bgapmatrix}T_{n-r+j}&\vrule&\tau\Im(v_S)+\Re(v_T)\\ \hline
 &\vrule&\sigma\end{bgapmatrix}\\
 &\hphantom{=Y_{n-r+j+1}\begin{bgapmatrix}S_{n-r+j}&\vrule&\Re(v_S)\\ \hline &\vrule&1\end{bgapmatrix}}
 +BL_E.
\end{align*}
Furthermore, it is easy to verify  that $B=Y_{n-r+j+1}HR$.
Now let  $P_{\bot}\in\mathbb{R}^{n\times (r-j-1)}$ satisfy
$P_{\bot}^\top \begin{row}P_{n-r+j}&p\end{row}=0$ and
$P_{\bot}^\top P_{\bot}=I_{r-j-1}$. Denoting
\begin{align*}
&A_1=\begin{bgapmatrix}S_{n-r+j}   &\vrule&\Re(v_S)\\   \hline
  &\vrule&1\end{bgapmatrix}+HRL_A, \\
&E_1=\begin{bgapmatrix}T_{n-r+j}  &\vrule&\tau\Im(v_S)+\Re(v_T)\\ \hline
 &\vrule&\sigma\end{bgapmatrix}+HRL_E,
\end{align*}
and writing
\[
	Y_{n-r+j+1}^\dagger = (Y_{n-r+j+1}^\top Y_{n-r+j+1})^{-1}Y_{n-r+j+1}^\top,
\]
then simple manipulations show  that
\begin{align*}
 \begin{bmatrix}Y_{n-r+j+1}^\dagger\\ K^\top\end{bmatrix}A\begin{bmatrix}P_{n-r+j}&p&\vrule&P_{\bot}\end{bmatrix}
&= \begin{bgapmatrix}A_1  &\vrule&A_{12} \\ \hline &\vrule& A_2\end{bgapmatrix},\\
 \begin{bmatrix}Y_{n-r+j+1}^\dagger\\ K^\top\end{bmatrix}E\begin{bmatrix}P_{n-r+j}&p&\vrule&P_{\bot}\end{bmatrix}
&= \begin{bgapmatrix}E_1  &\vrule&E_{12} \\ \hline &\vrule& E_2\end{bgapmatrix},\\
 \begin{bmatrix}Y_{n-r+j+1}^\dagger\\ K^\top\end{bmatrix}B\hphantom{\begin{bmatrix}P_{n-r+j}&p&\vrule&P_{\bot}\end{bmatrix}}
&= \begin{bmatrix} HR\\0	\end{bmatrix},
\end{align*}
where
\begin{align*}
A_{12}=Y_{n-r+j+1}^\dagger AP_{\bot}, &\quad A_2=K^\top AP_{\bot},\\
E_{12}=Y_{n-r+j+1}^\dagger EP_{\bot}, & \quad E_2=K^\top EP_{\bot}.
\end{align*}
Apparently, for the descriptor system
\[
\left(
\begin{bgapmatrix}E_1  &\vrule&E_{12} \\ \hline &\vrule& E_2\end{bgapmatrix}, \
\begin{bgapmatrix}A_1  &\vrule&A_{12} \\ \hline &\vrule& A_2\end{bgapmatrix},\
\begin{bmatrix} HR\\0	\end{bmatrix}
\right),
\]
which is equivalent to $(E, A, B)$, there are at most $(j+1)$ finite poles
assignable. This contradicts the fact that at least $(j+2)$ finite poles are
assignable, hence {\em (c)} holds.
\end{pf}

\subsection{Algorithm}\label{alg:Framwork_new}\label{subsection3.7}
The framework of our algorithm, referred to as ``\verb|DRSchurS|", is given in this section.
We assume that all infinite poles   appear together in $\mathfrak{L}$,
while  complex conjugate poles appear in pairs.
The   \verb|DRSchurS| algorithm below firstly assigns all infinite poles, then the finite ones.

\begin{algorithm}
\caption{Framework of our \texttt{DRSchurS} algorithm.}
\begin{algorithmic}[1]
\REQUIRE ~~\\
$A\in\mathbb{R}^{n\times n}, E\in\mathbb{R}^{n\times n}, B\in\mathbb{R}^{n\times m}$ and
$\mathfrak{L}=\{(\alpha_1, 0), \ldots, (\alpha_{n-r}, 0), (\alpha_{n-r+1}, \beta_{n-r+1}), \ldots, (\alpha_n, \beta_n)\}$.
\ENSURE ~~\\
The feedback matrix $F\in\mathbb{R}^{m\times n}$ and $G\in\mathbb{R}^{m\times n}$.
\STATE Compute the QR factorization: \\
\ \ \ $B=Q\begin{bmatrix}R^{\top}&0\end{bmatrix}^{\top}=\begin{bmatrix}Q_1&Q_2\end{bmatrix}
\begin{bmatrix}R^{\top}&0\end{bmatrix}^{\top}=Q_1R$.
\STATE Compute $P_{n-r}\in\mathbb{R}^{n\times (n-r)}$, $S_{n-r}\in\mathbb{R}^{(n-r)\times (n-r)}$
and $\Xi_{n-r}\in\mathbb{R}^{(n-m)\times (n-r)}$ by
Algorithm \ref{algorithm1};
       set $j=n-r$.
\WHILE{$j<n$}
\IF {$(\alpha_{j+1}, \beta_{j+1})\in\mathbb{R}\times \mathbb{R}$  }
     \STATE Compute $P_{j+1}\in\mathbb{R}^{n\times (j+1)}$, $S_{j+1}, T_{j+1}\in\mathbb{R}^{(j+1)\times (j+1)}$
     and $\Xi_{j+1}\in\mathbb{R}^{(n-m)\times(j+1)}$ by  Algorithm \ref{algorithm2}; set $j=j+1$.
\ELSE
    \STATE $P_{j+2}\in\mathbb{R}^{n\times (j+2)}$,  $S_{j+2}, T_{j+2}\in\mathbb{R}^{(j+2)\times (j+2)}$
    and  $\Xi_{j+2}\in\mathbb{R}^{(n-m)\times(j+2)}$ by  Algorithm \ref{algorithm3}; set $j=j+2$..
\ENDIF
\ENDWHILE
\STATE Compute the QR factorization: \\
\ \ \ $\Xi_n^\top=Q_X\begin{bmatrix}R_X^{\top}&0\end{bmatrix}^{\top}=
\begin{bmatrix}Q_{1,X}&Q_{2,X}\end{bmatrix}\begin{bmatrix}R_X^{\top}&0\end{bmatrix}^{\top}=Q_{1,X}R_X$.
\STATE Compute $X_{G,F}=Q_1Q_{2,X}^\top+Q_2\Xi_n$.
\STATE Compute $F=R^{-1}Q_1^{\top}(X_{G,F}S_nP_n^{\top}-A)$ and  \\
      \hspace*{1.1cm}  $G=R^{-1}Q_1^{\top}(X_{G,F}T_nP_n^{\top}-E)$.
\end{algorithmic}
\end{algorithm}

\section{Numerical examples}\label{section4}

In this section, we illustrate the performance of our  method  \verb|DRSchurS|
by applying it to  several examples, some from various references and
others generated randomly.

Similar to the definition of the precision of the assigned poles in \cite{GCQX}, we define
\[
precs=\max_{n-r+1\le j\le n}\left( \log \left| \frac{\lambda_j-\hat{\lambda}_j}{\lambda_j} \right| \right)
\]
to characterize the precision of the assigned finite  poles,
where $\lambda_j=\alpha_j/\beta_j$ are the chosen finite poles
and $\hat{\lambda}_j$ are  the computed ones from $(A+BF, E+BG)$.
Implicitly, we expect the number of  computed finite  eigenvalues
to be identical to that of those to be placed. Apparently,
smaller $``precs"$ indicates  more accurately computed finite eigenvalues.
To reveal  the robustness of the closed-loop system,
in addition to $\Delta_F^2(A+BF,E+BG)$ in \eqref{dep} (abbreviated as $\Delta_F^2$ in the following Tables and Figures),  the condition number of the closed-loop generalized eigenvectors matrix will also be displayed.
Specifically, assume that
\begin{align*}
&A+BF=Y\diag(\alpha_1, \ldots, \alpha_n)X,\\
&E+BG=Y\diag(\beta_1, \ldots, \beta_n)X,
\end{align*}
where $\{(\alpha_1, \beta_1), \ldots, (\alpha_n, \beta_n)\}$ are the  eigenvalues,
$Y$ and $X$ are nonsingular,  the Bauer-Fike type theorem in \cite{SSun} then shows that
$\kappa_F(X)=\|X\|_F\|X^{-1}\|_F$ would measure the sensitivity of the  eigenvalues
relative to  perturbations on the matrix pencil $(A+BF, E+BG)$.
When determining the nonsingular $X_{G,F}$ in Section~\ref{subsection3.5}, it is hoped
that $X_{G,F}$ would  be  well-conditioned. Accordingly, $\kappa_F(X_{G,F})$
is given explicitly for all examples.
In addition, $\|F\|_F$ and $\|G\|_F$, representing the energy  involved in the  feedback control,
are also  displayed.

All calculations are carried out  in MATLAB R2012a,
with the machine accuracy represented by $\epsilon\approx 2.2\times 10^{-16}$,
on an Intel\textregistered    Core\texttrademark i3 dual core 2.27 GHz machine with $2.00$ GB RAM.


\begin{exmp}\label{example1}
This illustrative set includes the  examples from
\cite{KNC, DP2, Duan, FLE, RZ, DP3,  DP1, MI, LC, CBS, SL, ZB, ZB1},
some of which  are  employed to compare the efficiency of the method proposed in \cite{Var2} and \verb|DRSchurS| here.
When assigning the finite poles, we firstly place the finite real poles in
non-descendent order, then the finite complex conjugate poles.
Tables~\ref{table_0} and \ref{table_1} present the numerical results,
with $\alpha(k)$ representing $\alpha \times 10^{k}$ for space saving.

Table~\ref{table_0} shows that \verb|DRSchurS| produces  comparable or  occasionally better results
than  the method proposed in \cite{Var2}.
For Example~$5$, the relative accuracy $``precs"$ produced by  \verb|DRSchurS| is
not that high, probably because some poles are close to the imaginary axis.
This is possibly a weakness of our algorithm.
 Notice also that the numerical results corresponding to Example~$6$ is not that satisfactory,
probably due to the difference in magnitudes of the entries in $A$.
\verb|DRSchurS| produces nice numerical results for Example~$4$ except for $\kappa_F(X)$,
indicating that it may  not be  wise to access an algorithm on only one criterion.

For all examples in Table~\ref{table_1}, the method put forward in \cite{Var2} is not tested
since it would fail. The reason may be one of the followings ---
$\lambda(A, E)\cap\mathfrak{L}\neq \emptyset$ or the sizes of the Jordan blocks corresponding
to some repeated poles cannot be determined.
Note that there are two more rows, in $9^\prime$ and $12^\prime$,
in Table~\ref{table_1}, which correspond to the inputs in rows~$1$ and $5$,  respectively,
but with all finite poles are placed before the infinite ones. Though we cannot prove  the
 feasibility of \verb|DRSchurS| when all infinite poles are assigned last,
 numerical results demonstrate better performance for certain examples.
 Our method  \verb|DRSchurS| produces  fairly low  relative accuracy $``precs"$  and very large $\kappa_F(X)$
for  Examples ~$13$ and  $18$, both possessing repeated  finite poles with algebraic multiplicities greater than $m$.
The treatment of repeated finite  poles deserves further investigation.

\tabcolsep 0.02in
\begin{table*}
\footnotesize
\centering
\caption{ Numerical results for Example  \ref{example1} (compare with the method in \cite{Var2} )}
\begin{tabular}{c|cccccc|cccccc}
\Xhline{1.0pt}
\multirow{2}*{} &  \multicolumn{6}{c|}{\text{Method in \cite{Var2}}}&\multicolumn{6}{c}{\texttt{DRSchurS}}\\
\Xcline{2-13}{0.5pt}
& $precs$ & $\Delta_F^2$ & $\|F\|_F$& $\|G\|_F$& $\kappa_F(Y)$& $\kappa_F(X)$
  &$precs$ & $\Delta_F^2$ & $\|F\|_F$& $\|G\|_F$& $\kappa_F(X_{G,F})$& $\kappa_F(X)$\\
\Xhline{0.7pt}
$1$\cite{KNC}&$ -4.1$ & $7.80(0)$& $1.43(0)$ &$1.09(0)$ &$1.57(0)$ &$3.75(0)$ &$-15.65$&$2.18(0)$&$1.71(0)$&$7.28(-1)$&$6.11(0)$&$1.21(1)$\\ 
$2$\cite{DP2}&$-12.19$&$7.55(5)$&$5.26(2)$&$1.34(2)$&$2.19(1)$&$6.18(1)$ &$-13.49$&$5.42(0)$&$3.00(0)$&$9.95(-1)$&$6.90(0)$&$5.80(2)$\\ 
$3$\cite{RZ}&$-14.26$&$6.13(1)$&$6.24(0)$&$6.27(0)$&$1.49(0)$&$1.15(1)$ &$-14.61$&$4.23(0)$&$5.11(0)$&$2.59(0)$&$3.12(0)$&$6.23(0)$\\  
$4$\cite{MI}&$-9.09$&$9.52(0)$&$1.96(-2)$&$9.20(-3)$&$1.17(6)$&$4.88(6)$&$-11.08$&$1.75(1)$&$3.23(0)$&$2.34(0)$&$7.95(0)$&$2.71(4)$ \\ 
$5$\cite{LC}&$-6.55$&$1.78(6)$&$2.48(1)$&$6.41(1)$&$9.72(2)$&$5.74(0)$&$-7.34$&$1.62(1)$&$9.89(0)$&$1.59(1)$&$1.56(2)$&$4.40(1)$ \\
$6$\cite{LC}&$-10.63$&$1.13(8)$&$5.49(3)$&$9.65(3)$&$1.14(1)$&$3.79(1)$&$-10.53$&$3.40(7)$&$4.67(3)$&$1.41(3)$&$4.31(1)$&$1.47(5)$ \\
$7$\cite{LC}&$-12.53$&$8.38(4)$&$5.54(2)$&$1.14(3)$&$5.72(-1)$&$5.16(-1)$&$-14.23$&$5.65(1)$&$1.18(2)$&$2.97(1)$&$1.29(1)$&$1.11(2)$ \\
$8$\cite{LC}&$-12.45$&$1.16(4)$&$6.77(0)$&$2.76(1)$&$2.48(2)$&$3.34(2)$&$-14.19$&$2.31(0)$&$6.18(0)$&$5.08(1)$&$9.22(1)$&$2.60(1)$ \\ 
\Xhline{1.0pt}
\end{tabular}
\label{table_0}
\end{table*}

\tabcolsep 0.08in
\begin{table*}
\normalsize
\centering
\caption{ Numerical results for Example  \ref{example1}}
\begin{tabular}{c|c|c|c|c|c|c}
\Xhline{1.0pt}
Num.& $precs$ & $\Delta_F^2$ & $\|F\|_F$& $\|G\|_F$& $\kappa_F(X_{G,F})$& $\kappa_F(X)$\\
\Xhline{0.5pt}
$9$\cite{Duan} &$-15.86$&$7.80(-1)$&$2.16(0)$&$1.45(0)$&$3.00(0)$&$8.38(15)$\\
$9^\prime$\cite{Duan}&$-17$&$2.47(-32)$&$1.58(0)$&$1.58(0)$&$3.16(0)$&$3.00(0)$\\
$10$\cite{FLE}&$-15.48$&$2.97(-1)$&$1.23(0)$&$1.37(0)$&$4.09(0)$&$6.40(0)$\\ 
$11$\cite{DP3}&$-15.78$&$2.80(0)$&$1.12(0)$&$1.89(0)$&$3.00(0)$&$4.57(0)$\\  
$12$\cite{DP3}&$-8.52$&$5.03(-1)$&$1.54(0)$&$1.12(0)$&$3.05(0)$&$1.86(8)$\\ 
$12^\prime$\cite{DP3}&$-17.00$&$5.00(-1)$&$1.27(0)$&$1.35(0)$&$3.06(0)$&$4.00(0)$\\ 
$13$\cite{DP1}&$-5.09$&$1.18(1)$&$2.28(0)$&$4.41(0)$&$6.76(0)$&$1.29(11)$ \\
$14$\cite{CBS}&$-15.09$&$4.67(-1)$&$4.90(0)$&$2.39(0)$&$3.49(0)$&$3.29(0)$ \\ 
$15$\cite{SL}&$-15.48$&$1.84(0)$&$1.92(0)$&$3.32(-2)$&$4.00(0)$&$1.53(1)$ \\ 
$16$\cite{ZB}&$-15.35$&$2.45(0)$&$1.89(0)$&$1.15(0)$&$4.30(0)$&$2.79(1)$\\ 
$17$\cite{ZB}&$-14.81$&$2.72(0)$&$2.88(0)$&$1.46(0)$&$6.39(0)$&$3.06(1)$ \\ 
$18$\cite{ZB1}&$-5.63$&$9.02(0)$&$4.41(0)$&$1.39(0)$&$6.41(0)$&$2.12(11)$\\ 
\Xhline{1.0pt}
\end{tabular}
\label{table_1}
\end{table*}
\end{exmp}

\begin{exmp}\label{example2}
This test set contains $255$ random examples, where $n = 6, 15, 30$,
the rank of the singular matrix $E$ equals $2, \lfloor\frac{n}{2}\rfloor, n-1$,
and $m = 2, \lfloor\frac{n}{2}\rfloor, n-2$.
For each triple $(n, \rank(E), m)$, the number  of the finite
poles, denoted by $r$,  increases from
$\rank(\begin{row}E&B\end{row} )-m$ to $\rank(\begin{row}E&B\end{row})$ in increment of $1$.
Note that for randomly generated examples, we usually have
$\rank(\begin{row}E&B\end{row})=\min\{n, \rank(E)+m\}$, bringing
$r=\rank(E), \rank(E)+1, \ldots, \rank(E)+m$ or $r=(n-m), (n-m)+1, \ldots, n$.
All examples are generated randomly as follows.
For a  fixed  4-turple $(n, \rank(E), m, r)$, we first randomly generate five  matrices
$A, E\in\mathbb{R}^{n\times n}$, $B\in\mathbb{R}^{n\times m}$,
$W\in\mathbb{R}^{r\times r}$, $Y\in\mathbb{R}^{r\times r}$ by the MATLAB function \verb|randn|,
and set the finite poles as $\mathfrak{L}_1=\verb|eig|(W,Y)$ and $\mathfrak{L}=\{\underbrace{\infty, \ldots, \infty}_{n-r}, \mathfrak{L}_1\}$.
Then  compute the QR factorization  $E=Q_ER_E$,  reset the
$(n-\rank(E))\times (n-\rank(E))$ principal submatrix of $R_E$  to $0$,
and reassign $E=Q_ER_EQ_E^{-1}$.
Taking the resulting $A, E, B$ and $\mathfrak{L}$ as the inputs, we apply \verb|DRSchurS|.

All  numerical results are  plotted in the following figures.
Specifically, with the triple  $(n, \rank(E), m)$ fixed,
the precision of the computed finite poles $precs$, $\Delta_F^2$,
the norms of $F$ and $G$ and the condition numbers of $X_{G,F}$ and the
generalized eigenvectors matrix $X$, with respect to $r$, are
displayed in Fig.~\ref{figprecs} to  Fig.~\ref{figcondV},  respectively.
For  each fixed $n$, the three subfigures correspond to $m=2, \lfloor\frac{n}{2}\rfloor$ and $n-2$, respectively,
where the three different lines match the three distinct
$\rank(E)=2, \lfloor\frac{n}{2}\rfloor, n-1$.
The $x$-axis represents $r$, which varies from $\rank(E)$ to $(\rank(E)+m)$ or $(n-m)$ to  $n$,
and the values on the $y$-axis are mean values over $50$ trials for a certain 4-turple $(n, \rank(E), m, r)$.

Fig.~\ref{figprecs}  reveals that \verb|DRSchurS| can produce high relative accuracy of the
assigned finite poles for all the examples  except  the special  case when
$n=30, \rank(E)=29, m=2$.  In that case, $r = 28, 29, 30$,
and the decline of the relative accuracy  can be attributed to the differences between
the number of the  finite poles and $m$.
In addition, when $\rank(E)=2$, the value of $precs$ exhibits an ascending trend   with
respect to $r$, probably due to the low rank of $E$. As for
$\Delta_F^2$, $\|F\|_F$ and   $\|G\|_F$,
they all  display an ascending trend  for $\rank(E)=2$,
but an oscillatory or a downward trend for $\rank(E)=\lfloor\frac{n}{2}\rfloor, \ n-1$.
It shows that $\kappa_F(X_{G,F})$ increases with respect to  $r$, probably due to
the greater contributions of $S$ and $T$ to $X_{G,F}$ when $r$ gets larger;
it decreases with respect to  $m$ since the freedom in $X_{G,F}$ increases with respect to $m$.

\begin{figure*}[!t]
\centering
\begin{tabular}{ccc}
    \begin{minipage}[t]{0.3\textwidth}
    \includegraphics[width=1.8in]{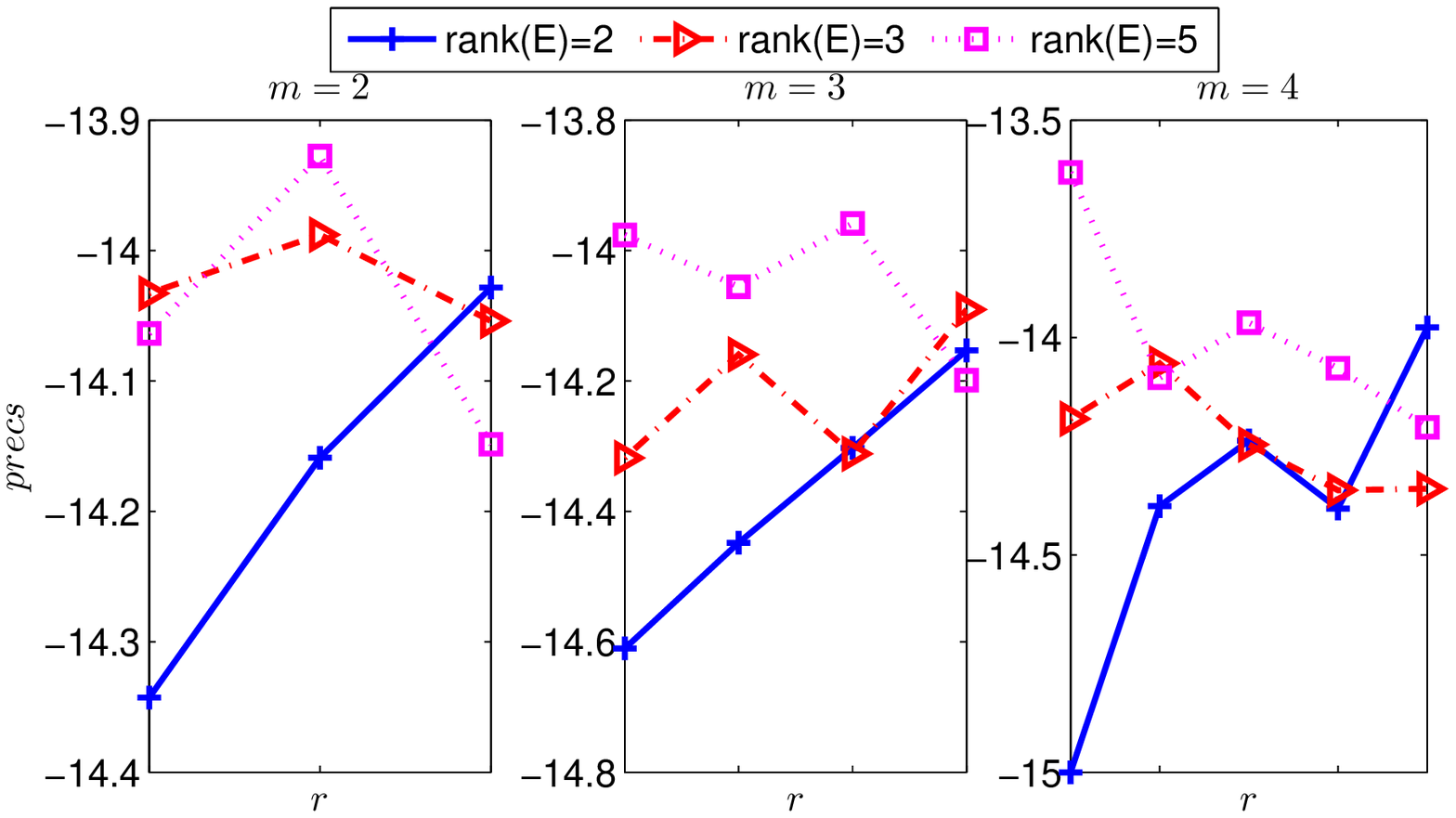}
    \caption*{\uppercase\expandafter{\romannumeral1}: $n=6$}
    \end{minipage}&
    \begin{minipage}[t]{0.3\textwidth}
    \includegraphics[width=1.8in]{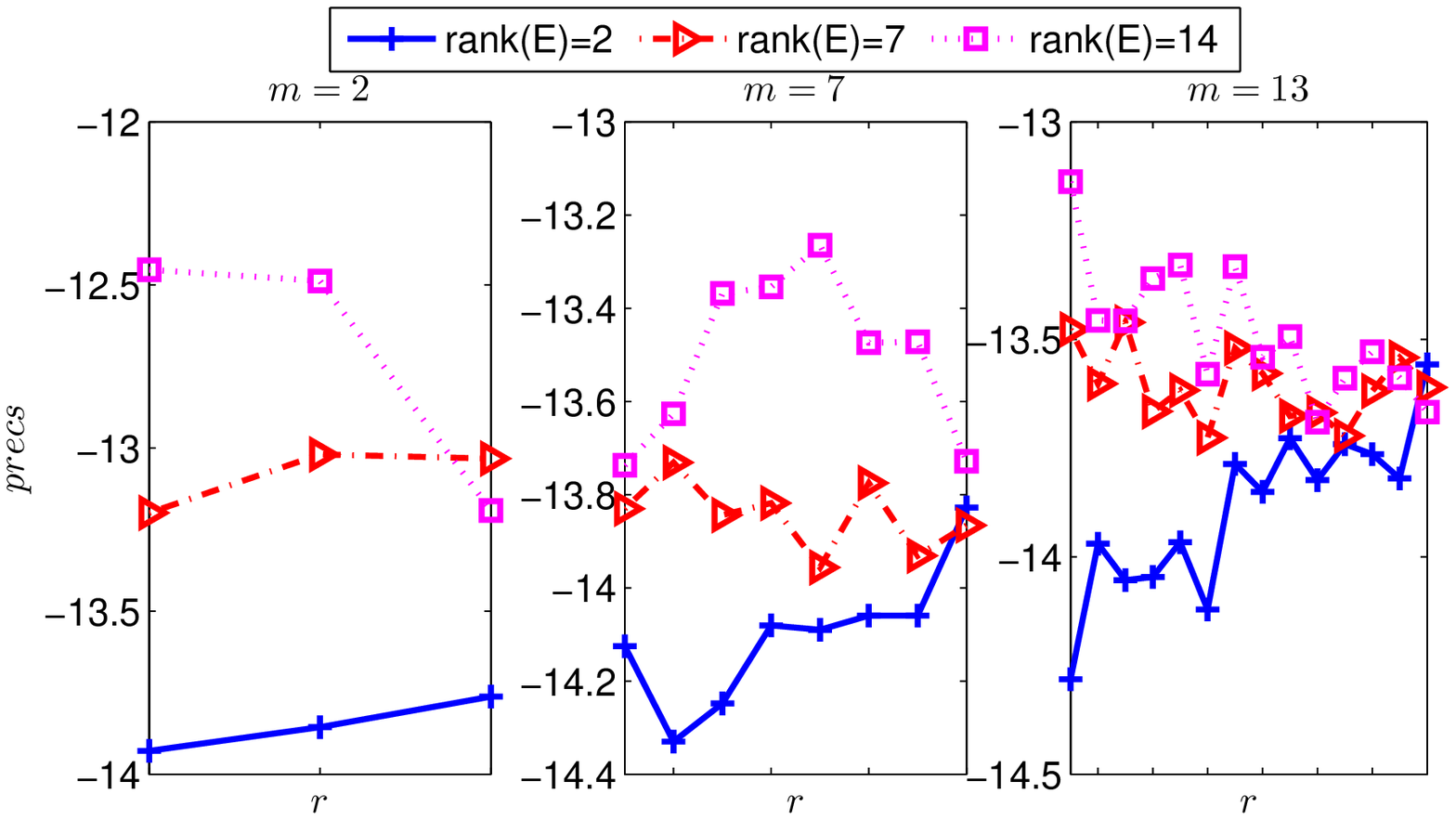}
    \caption*{\uppercase\expandafter{\romannumeral2}: $n=15$}
    \end{minipage}&
    \begin{minipage}[t]{0.3\textwidth}
    \includegraphics[width=1.8in]{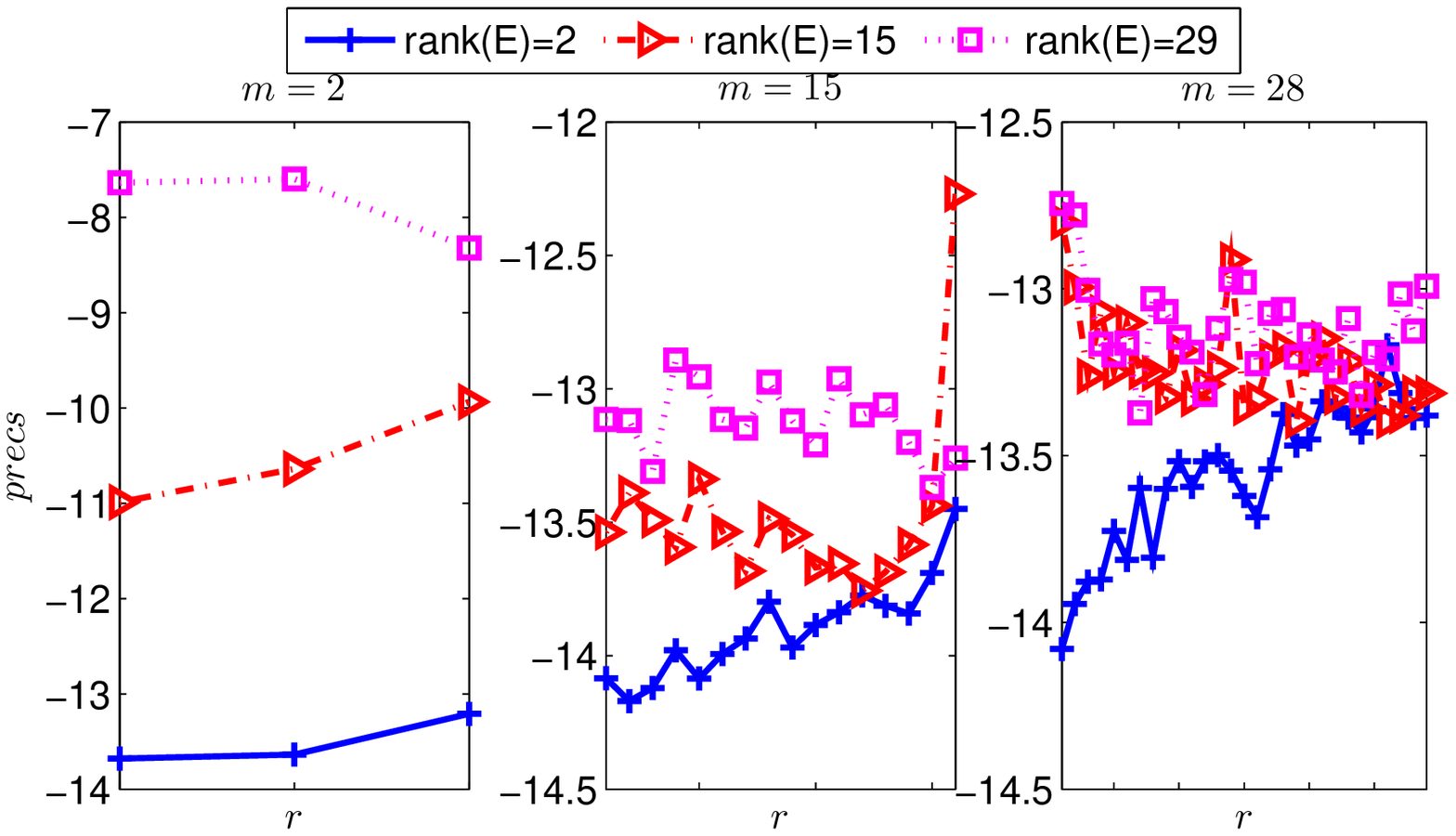}
    \caption*{\uppercase\expandafter{\romannumeral3}: $n=30$}
    \end{minipage}
\end{tabular}
\caption{Precision of the assigned finite poles}\label{figprecs}
\end{figure*}

\begin{figure*}[!t]
\centering
\begin{tabular}{ccc}
    \begin{minipage}[t]{0.3\textwidth}
    \includegraphics[width=1.8in]{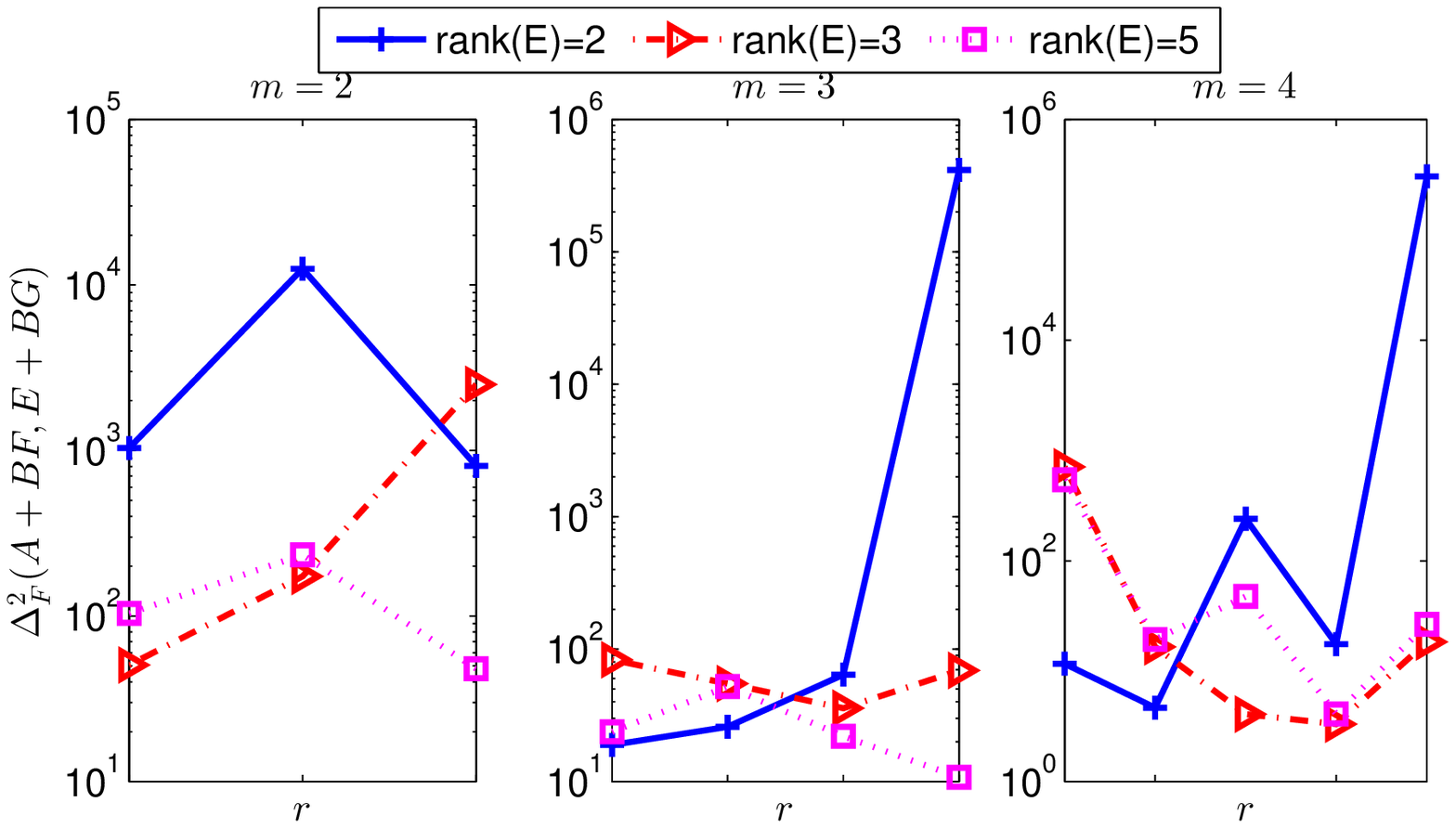}
    \caption*{\uppercase\expandafter{\romannumeral1}: $n=6$}
    \end{minipage}&
    \begin{minipage}[t]{0.3\textwidth}
    \includegraphics[width=1.8in]{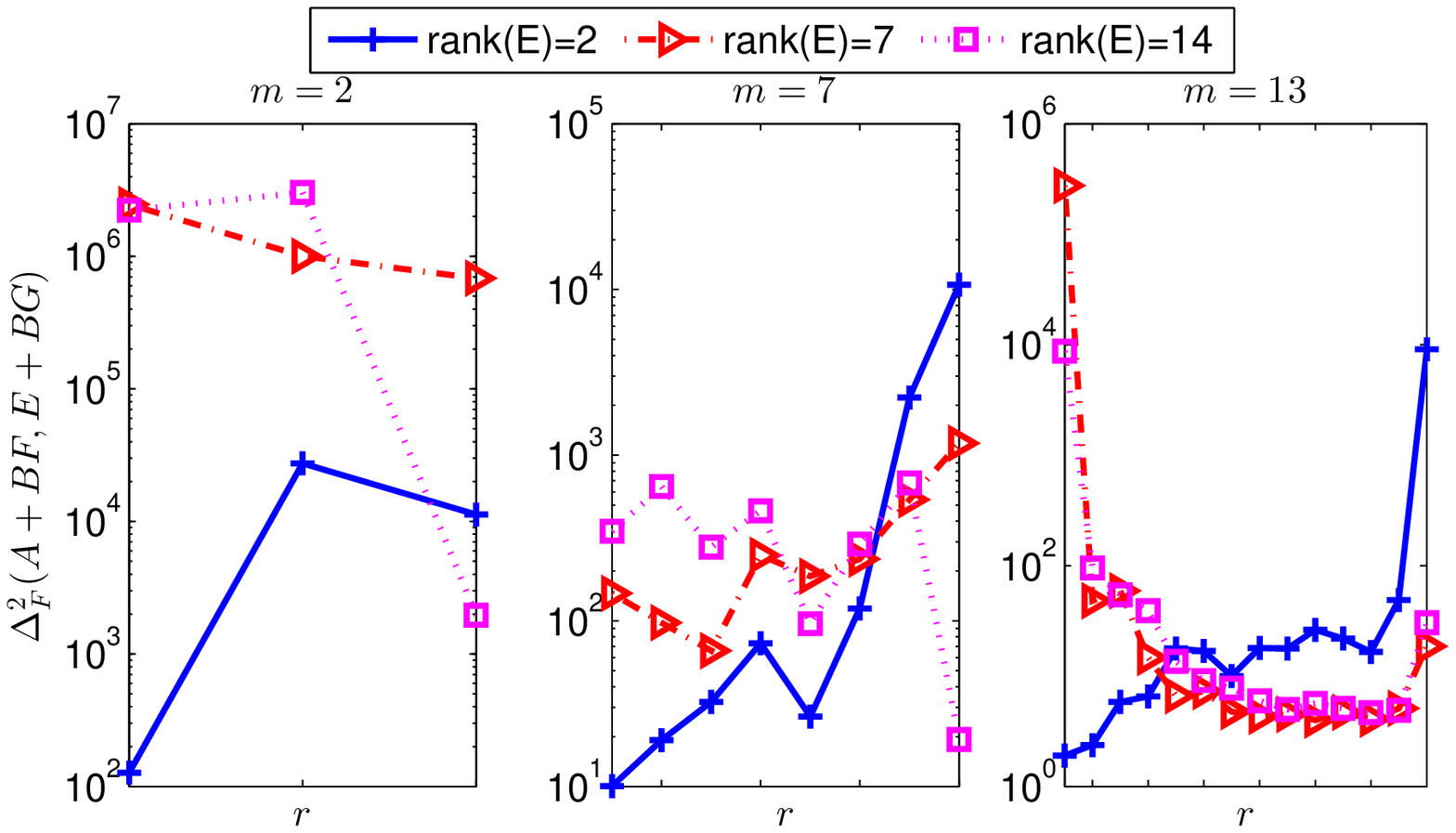}
    \caption*{\uppercase\expandafter{\romannumeral2}: $n=15$}
    \end{minipage}&
    \begin{minipage}[t]{0.3\textwidth}
    \includegraphics[width=1.8in]{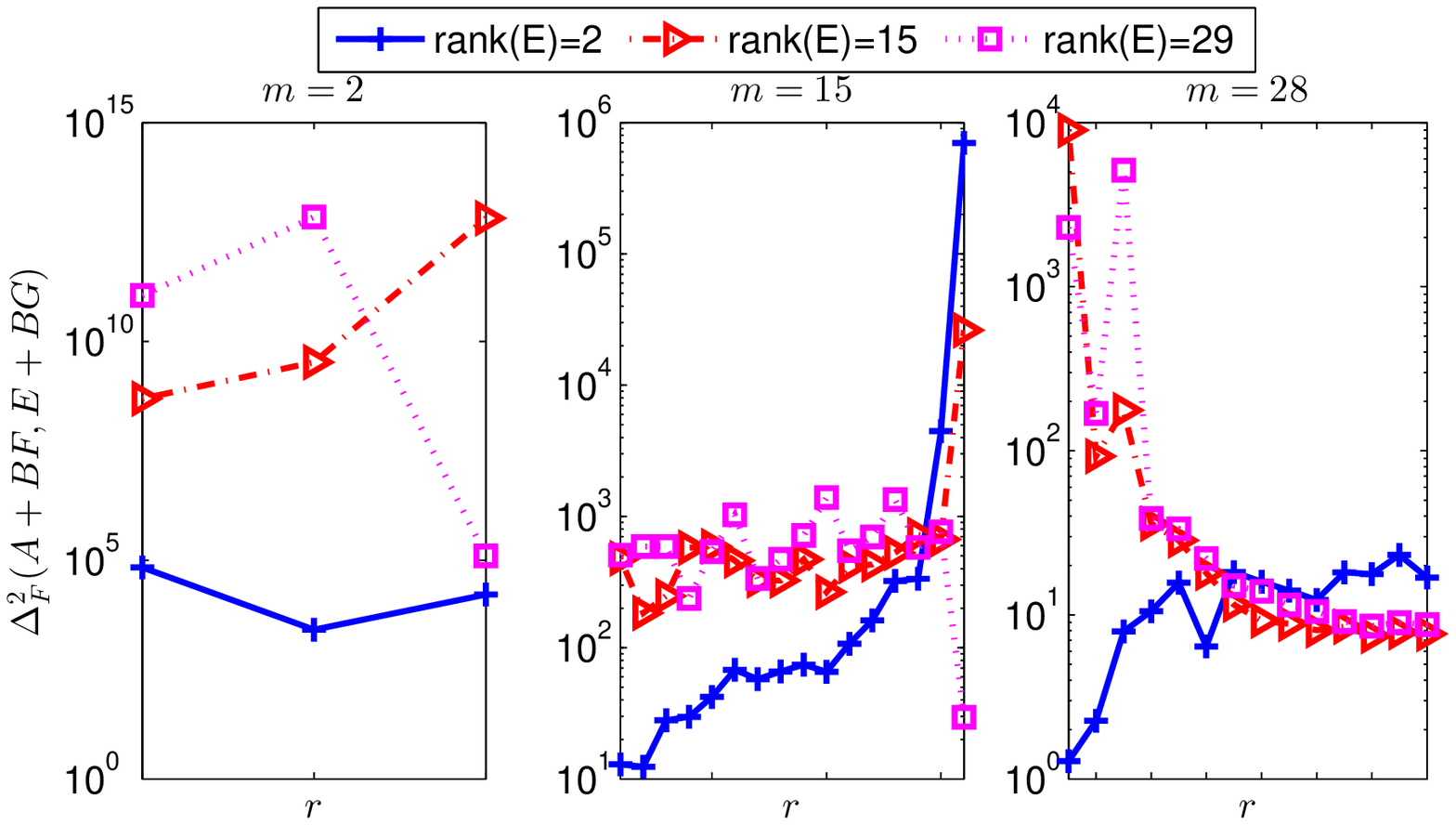}
    \caption*{\uppercase\expandafter{\romannumeral3}: $n=30$}
    \end{minipage}
\end{tabular}
\caption{Numerical results for $\Delta_F^2$}\label{figdep}
\end{figure*}

\begin{figure*}[!t]
\centering
\begin{tabular}{ccc}
    \begin{minipage}[t]{0.3\textwidth}
    \includegraphics[width=1.8in]{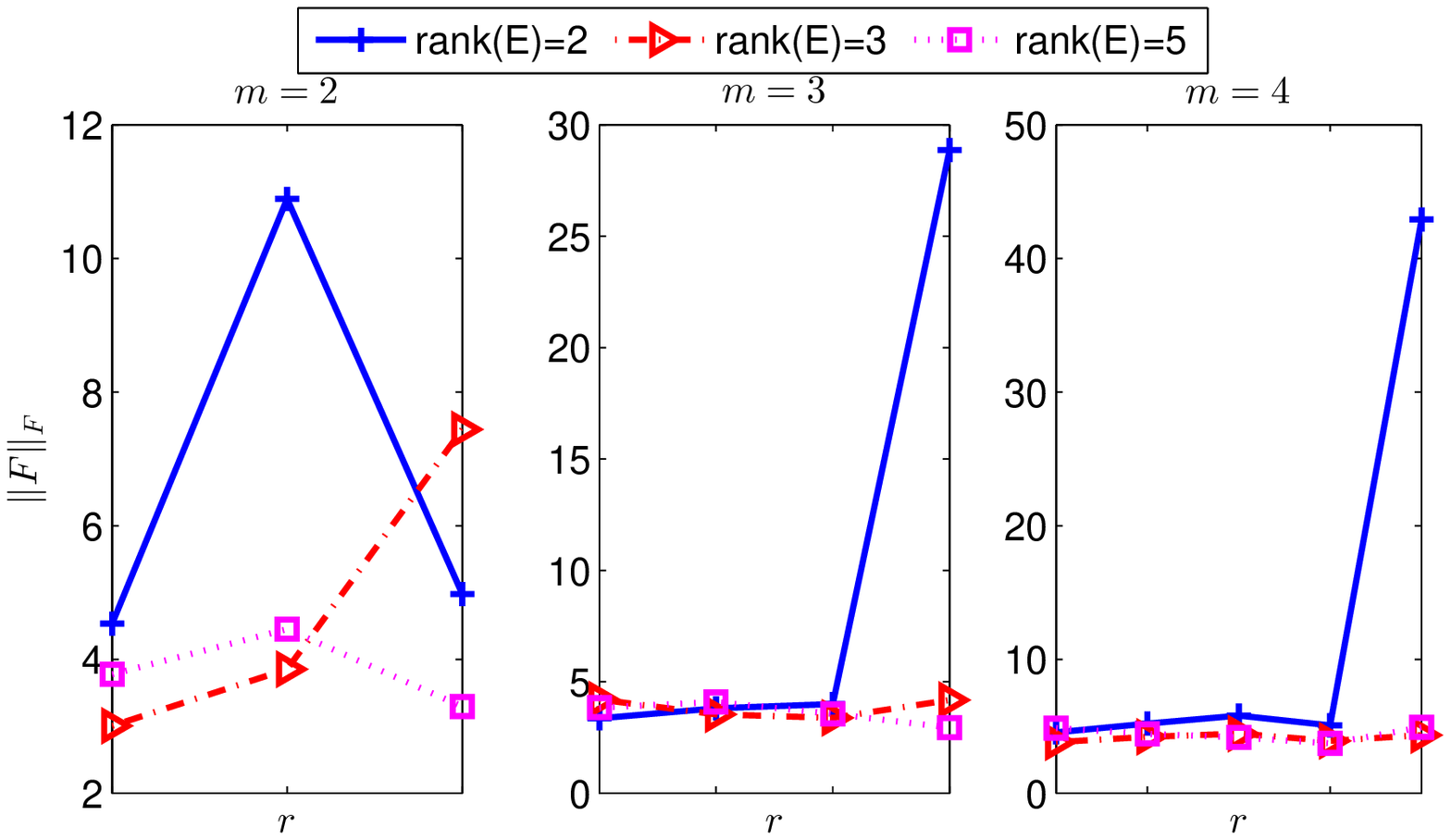}
    \caption*{\uppercase\expandafter{\romannumeral1}: $n=6$}
    \end{minipage}&
    \begin{minipage}[t]{0.3\textwidth}
    \includegraphics[width=1.8in]{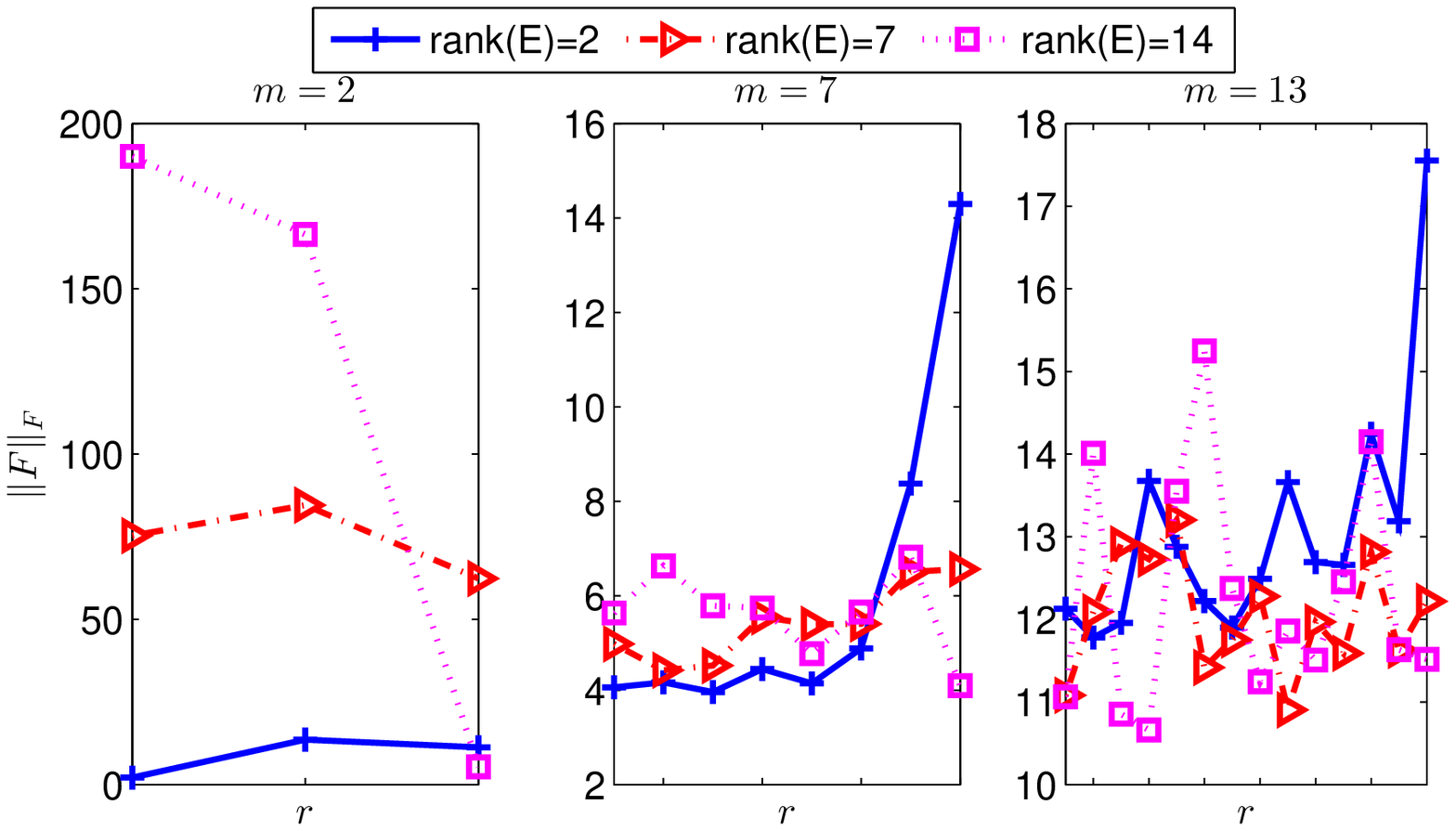}
    \caption*{\uppercase\expandafter{\romannumeral2}: $n=15$}
    \end{minipage}&
    \begin{minipage}[t]{0.3\textwidth}
    \includegraphics[width=1.8in]{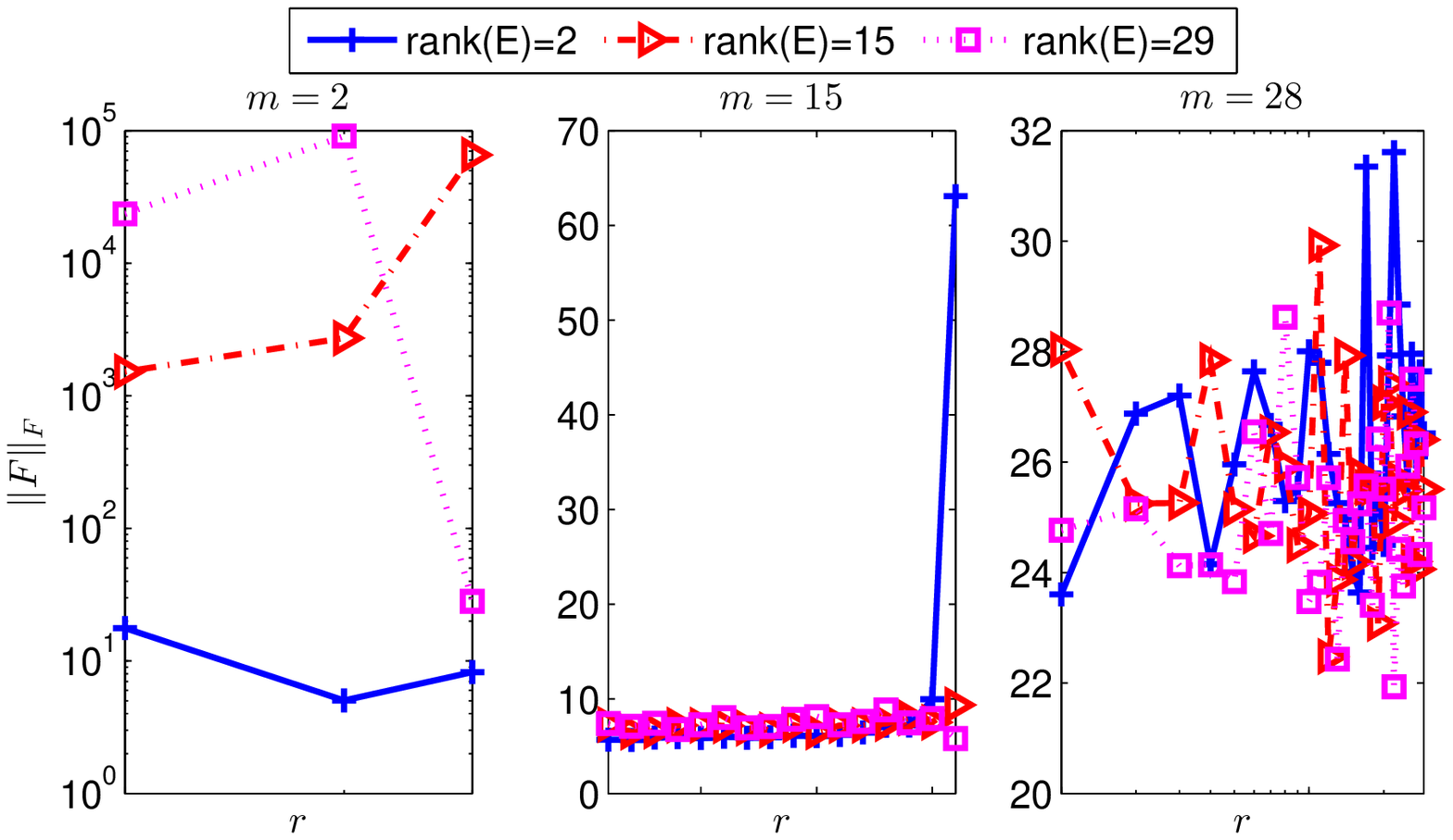}
    \caption*{\uppercase\expandafter{\romannumeral3}: $n=30$}
    \end{minipage}
\end{tabular}
\caption{Norm of the proportional part $F$ }\label{fignormF}
\end{figure*}

\begin{figure*}[!t]
\centering
\begin{tabular}{ccc}
    \begin{minipage}[t]{0.3\textwidth}
    \includegraphics[width=1.8in]{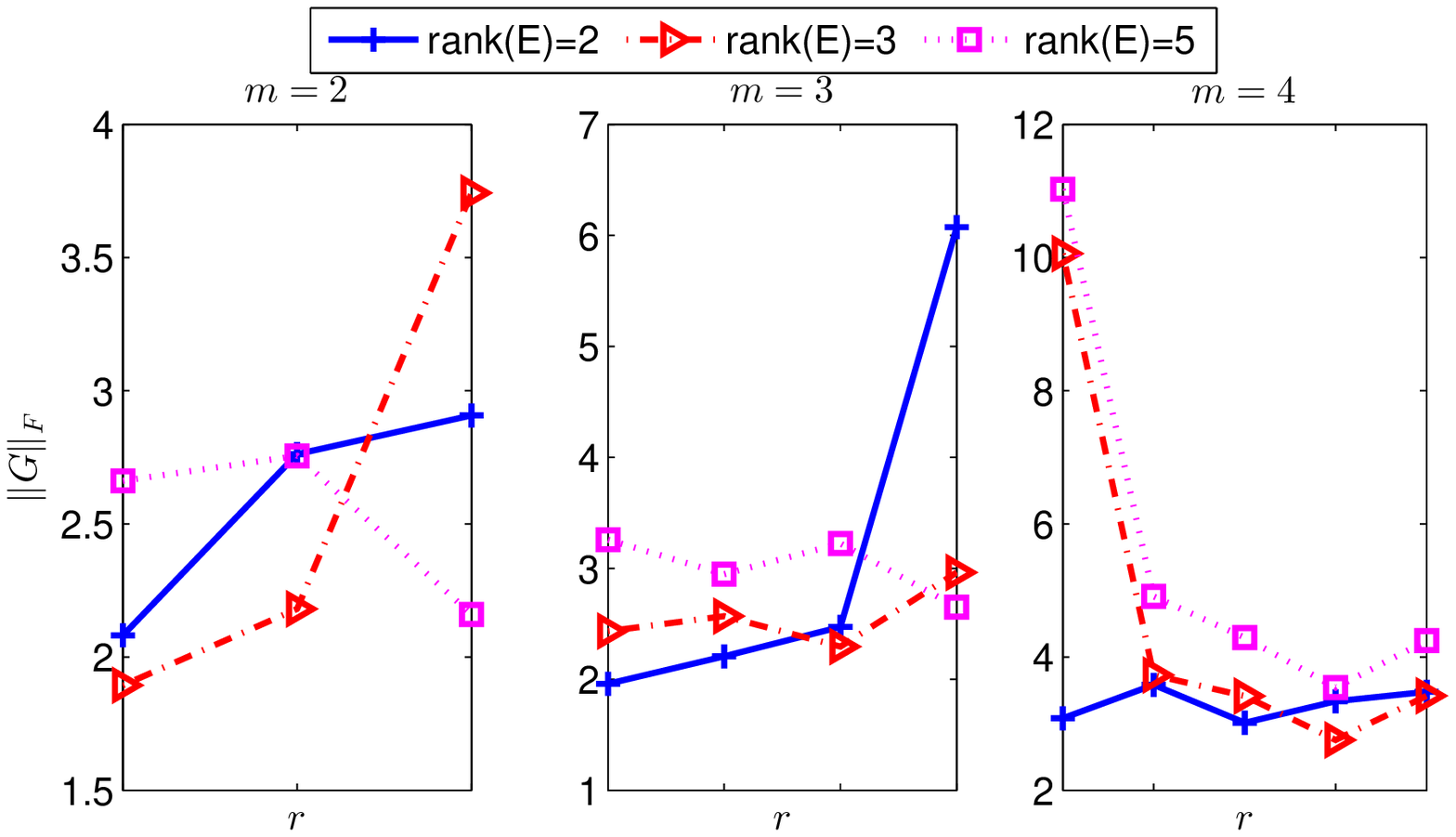}
    \caption*{\uppercase\expandafter{\romannumeral1}: $n=6$}
    \end{minipage}&
    \begin{minipage}[t]{0.3\textwidth}
    \includegraphics[width=1.8in]{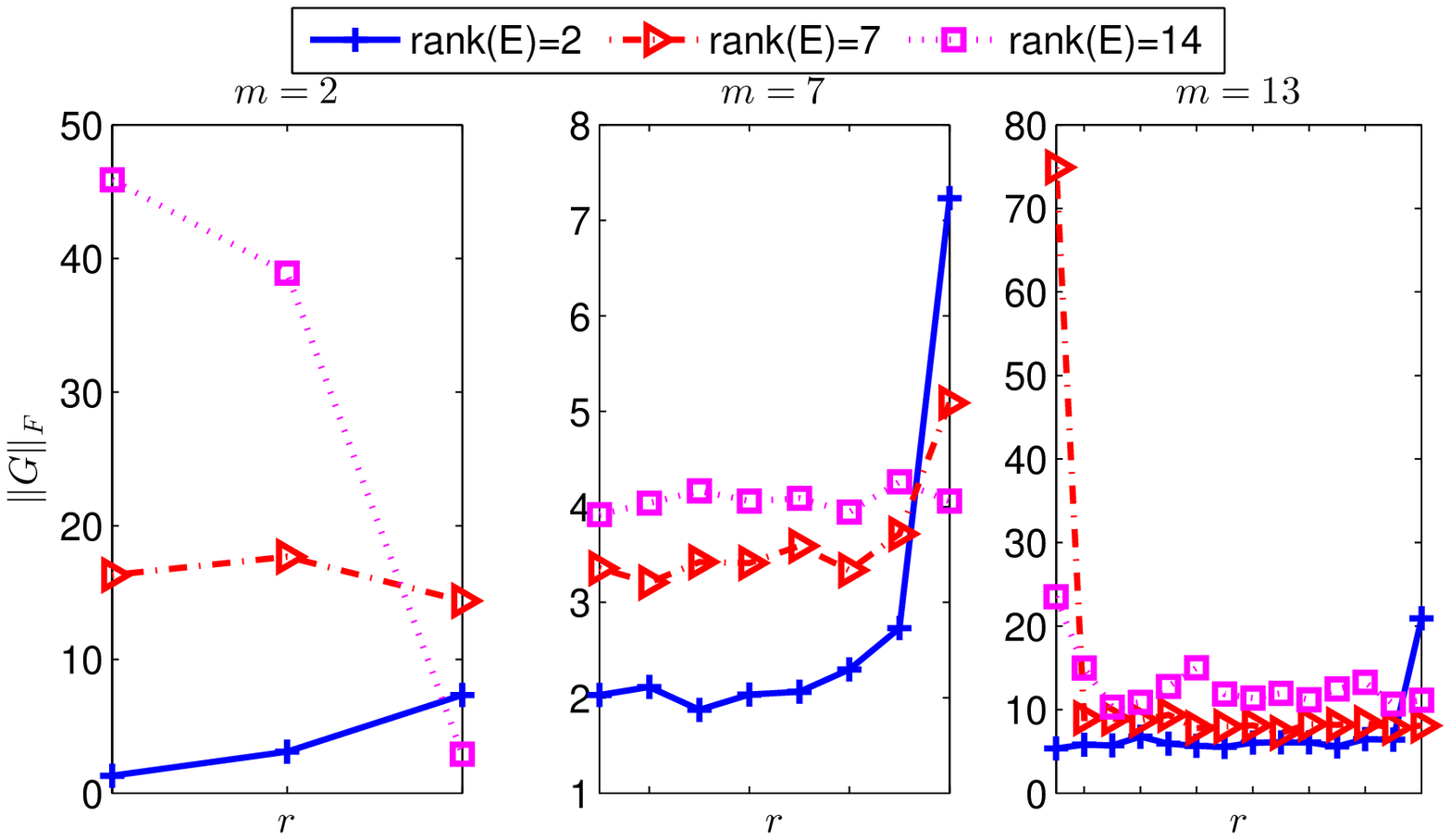}
    \caption*{\uppercase\expandafter{\romannumeral2}: $n=15$}
    \end{minipage}&
    \begin{minipage}[t]{0.3\textwidth}
    \includegraphics[width=1.8in]{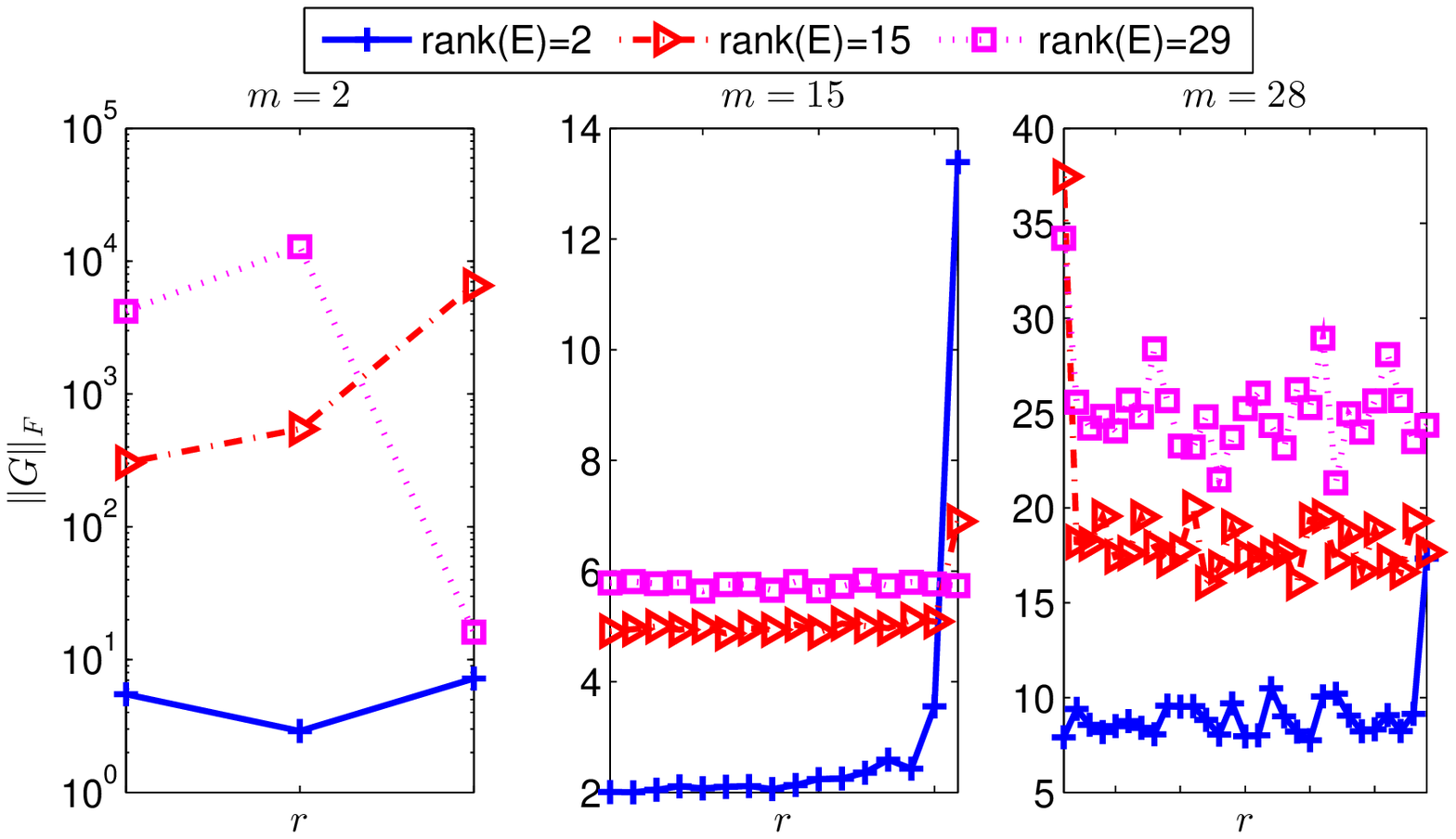}
    \caption*{\uppercase\expandafter{\romannumeral3}: $n=30$}
    \end{minipage}
\end{tabular}
\caption{Norm of the derivative part $G$ }\label{fignormG}
\end{figure*}

\begin{figure*}[!t]
\centering
\begin{tabular}{ccc}
    \begin{minipage}[t]{0.3\textwidth}
    \includegraphics[width=1.8in]{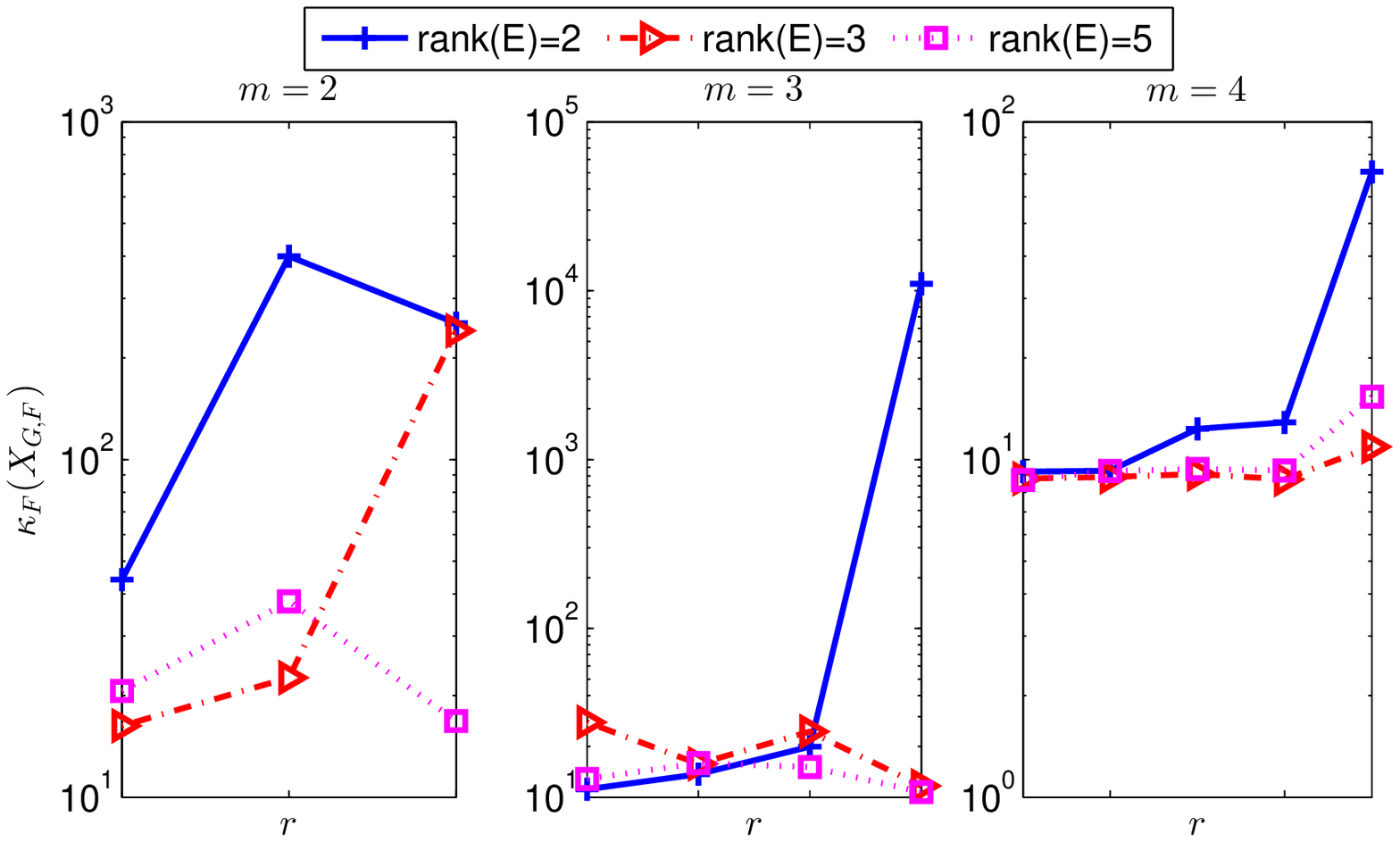}
    \caption*{\uppercase\expandafter{\romannumeral1}: $n=6$}
    \end{minipage}&
    \begin{minipage}[t]{0.3\textwidth}
    \includegraphics[width=1.8in]{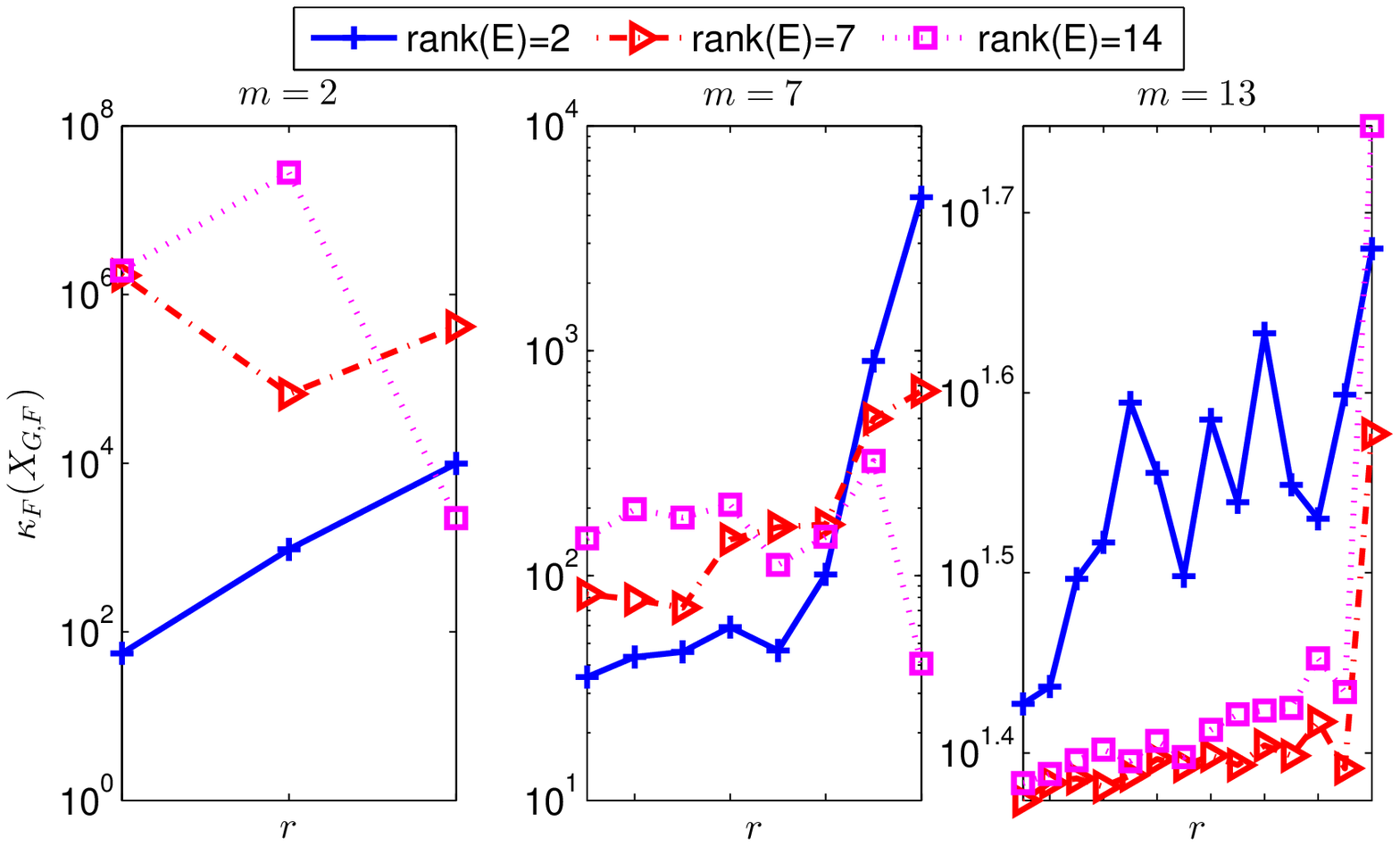}
    \caption*{\uppercase\expandafter{\romannumeral2}: $n=15$}
    \end{minipage}&
    \begin{minipage}[t]{0.3\textwidth}
    \includegraphics[width=1.8in]{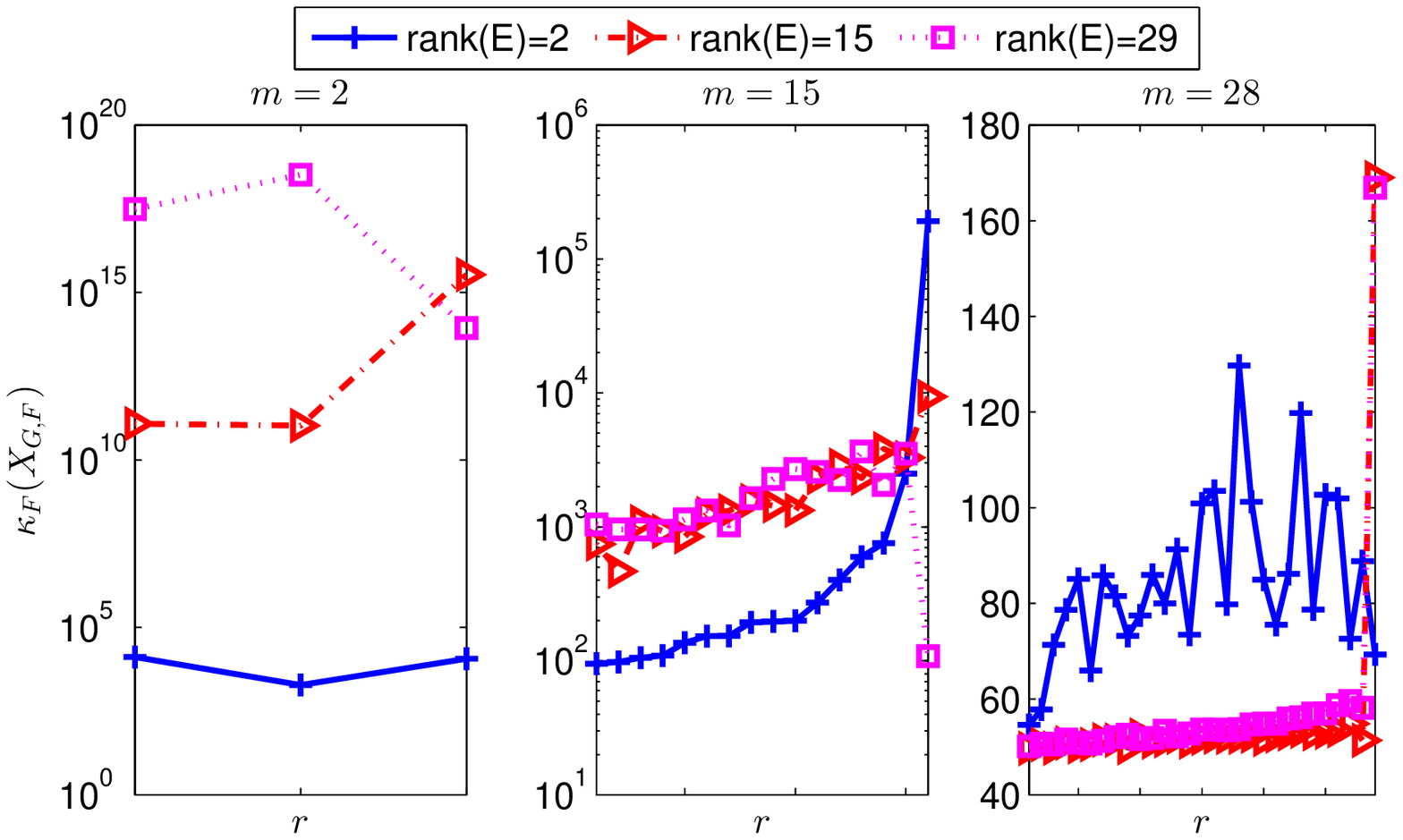}
    \caption*{\uppercase\expandafter{\romannumeral3}: $n=30$}
    \end{minipage}
\end{tabular}
\caption{Condition number of $X_{G,F}$ }\label{figcondX}
\end{figure*}

\begin{figure*}[!t]
\centering
\begin{tabular}{ccc}
    \begin{minipage}[t]{0.3\textwidth}
    \includegraphics[width=1.8in]{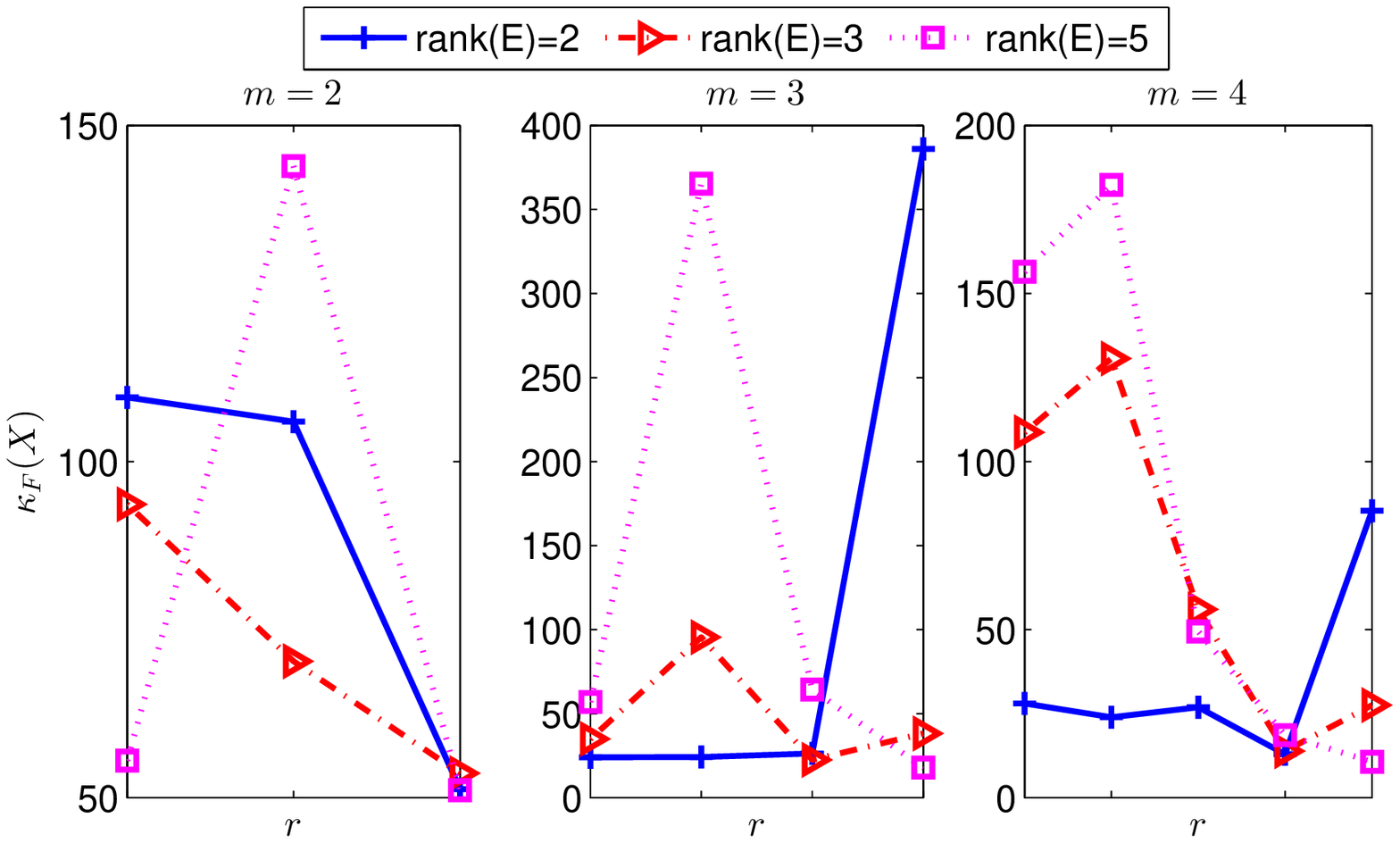}
    \caption*{\uppercase\expandafter{\romannumeral1}: $n=6$}
    \end{minipage}&
    \begin{minipage}[t]{0.3\textwidth}
    \includegraphics[width=1.8in]{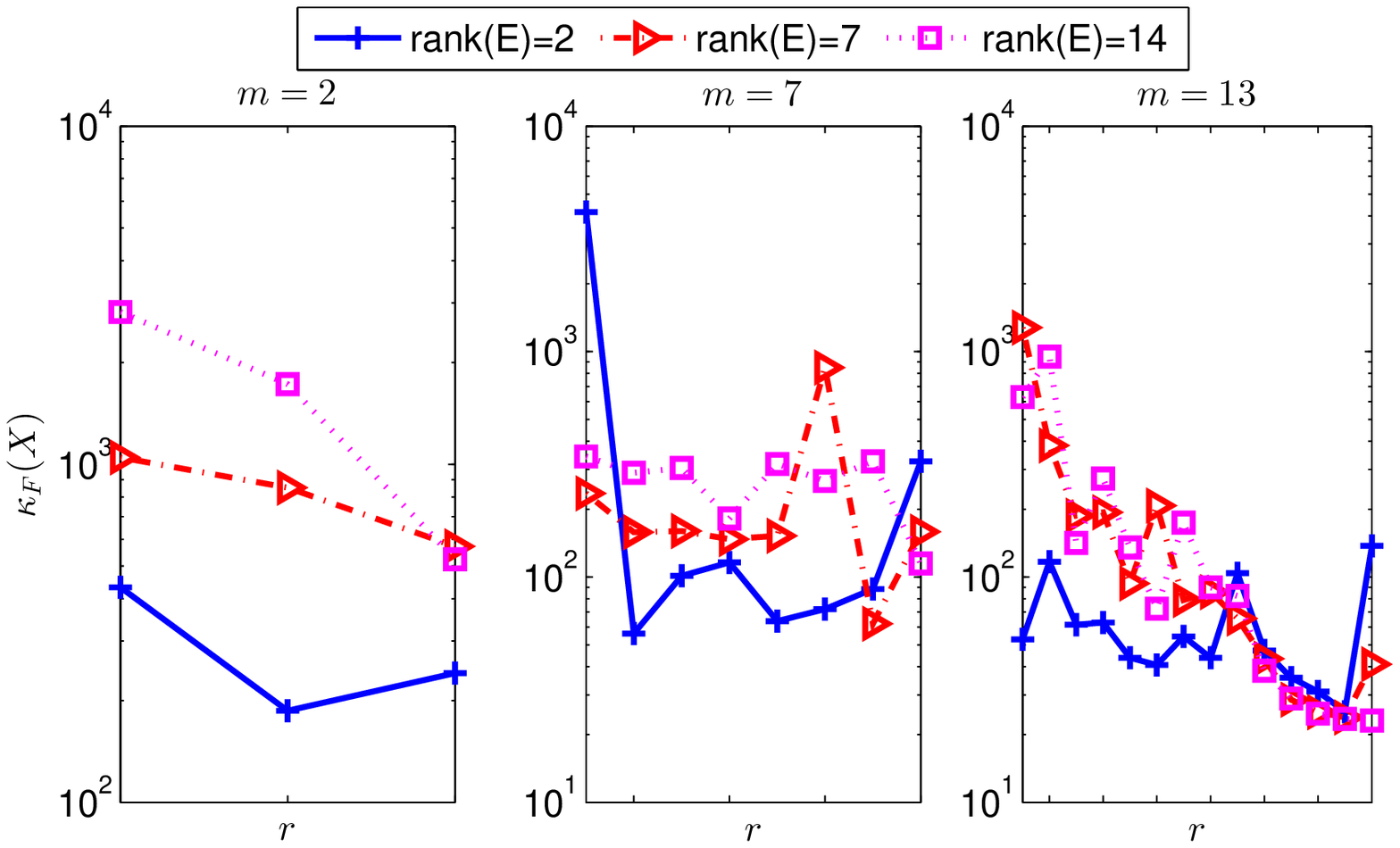}
    \caption*{\uppercase\expandafter{\romannumeral2}: $n=15$}
    \end{minipage}&
    \begin{minipage}[t]{0.3\textwidth}
    \includegraphics[width=1.8in]{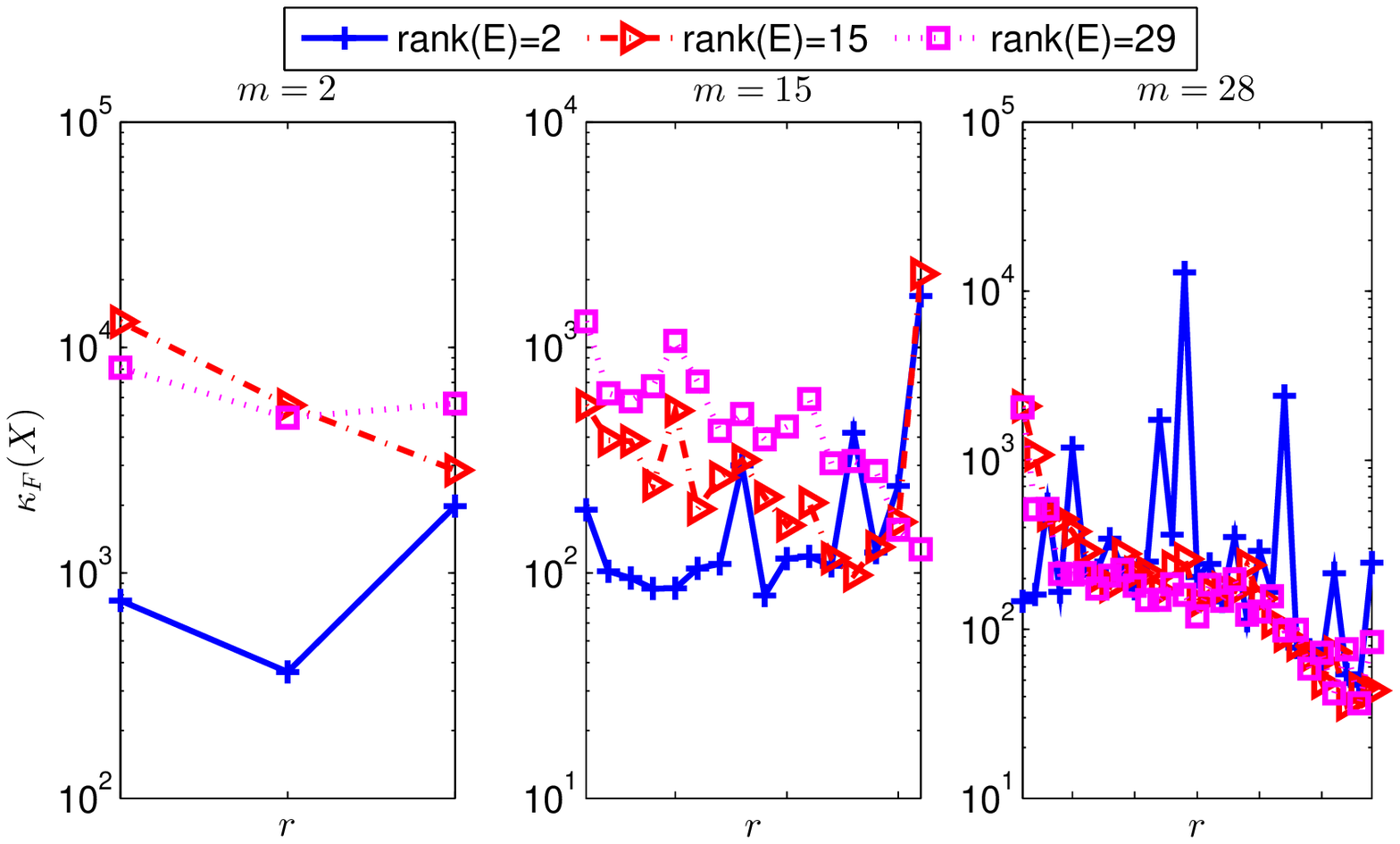}
    \caption*{\uppercase\expandafter{\romannumeral3}: $n=30$}
    \end{minipage}
\end{tabular}
\caption{Condition number of $X$ }\label{figcondV}
\end{figure*}
\end{exmp}

\section{Conclusions}\label{section5}
Based on the remarkable results in \cite{BMN},  a new direct method \verb|DRSchurS| for the \RPAPDSF\ is proposed in this paper.  Using   the generalized real Schur form of the closed-loop system matrix
pencil, \verb|DRSchurS| is  capable  of minimizing a robust measure, which is  closely related to
the  departure from normality of the closed-loop system matrix pencil, via some standard eigen-problems.
Several numerical examples demonstrate that \verb|DRSchurS| solves \RPAPDSF, producing robust closed-loop systems with highly accurate closed-loop finite poles.

For future work, we may further investigate the assignment of repeated finite poles, as well as how the freedom in the first eigenvector for the finite poles and the order of poles in $\mathcal{L}$ can be best exploited.


\bibliographystyle{plain}        


\end{document}